\documentclass[a4paper, british]{amsart}

\usepackage{libertine}
\usepackage[libertine]{newtxmath}

\usepackage{babel}
\usepackage{enumitem}
\usepackage{hyperref}
\usepackage[utf8]{inputenc}
\usepackage{newunicodechar}
\usepackage{mathtools}
\usepackage{varioref}
\usepackage[arrow,curve,matrix]{xy}

\usepackage{colortbl}
\usepackage{graphicx}
\usepackage{tikz}

\usepackage{mathrsfs}

\definecolor{linkred}{rgb}{0.7,0.2,0.2}
\definecolor{linkblue}{rgb}{0,0.2,0.6}

\numberwithin{figure}{section}

\usepackage[hyperpageref]{backref}

\setdescription{labelindent=\parindent, leftmargin=2\parindent}
\setitemize[1]{labelindent=\parindent, leftmargin=2\parindent}
\setenumerate[1]{labelindent=0cm, leftmargin=*, widest=iiii}

\newunicodechar{α}{\ensuremath{\alpha}}
\newunicodechar{β}{\ensuremath{\beta}}
\newunicodechar{χ}{\ensuremath{\chi}}
\newunicodechar{δ}{\ensuremath{\delta}}
\newunicodechar{ε}{\ensuremath{\varepsilon}}
\newunicodechar{Δ}{\ensuremath{\Delta}}
\newunicodechar{η}{\ensuremath{\eta}}
\newunicodechar{γ}{\ensuremath{\gamma}}
\newunicodechar{Γ}{\ensuremath{\Gamma}}
\newunicodechar{ι}{\ensuremath{\iota}}
\newunicodechar{κ}{\ensuremath{\kappa}}
\newunicodechar{λ}{\ensuremath{\lambda}}
\newunicodechar{Λ}{\ensuremath{\Lambda}}
\newunicodechar{ν}{\ensuremath{\nu}}
\newunicodechar{μ}{\ensuremath{\mu}}
\newunicodechar{ω}{\ensuremath{\omega}}
\newunicodechar{Ω}{\ensuremath{\Omega}}
\newunicodechar{π}{\ensuremath{\pi}}
\newunicodechar{Π}{\ensuremath{\Pi}}
\newunicodechar{φ}{\ensuremath{\phi}}
\newunicodechar{Φ}{\ensuremath{\Phi}}
\newunicodechar{ψ}{\ensuremath{\psi}}
\newunicodechar{Ψ}{\ensuremath{\Psi}}
\newunicodechar{ρ}{\ensuremath{\rho}}
\newunicodechar{σ}{\ensuremath{\sigma}}
\newunicodechar{Σ}{\ensuremath{\Sigma}}
\newunicodechar{τ}{\ensuremath{\tau}}
\newunicodechar{θ}{\ensuremath{\theta}}
\newunicodechar{Θ}{\ensuremath{\Theta}}
\newunicodechar{ξ}{\ensuremath{\xi}}
\newunicodechar{Ξ}{\ensuremath{\Xi}}
\newunicodechar{ζ}{\ensuremath{\zeta}}

\newunicodechar{ℓ}{\ensuremath{\ell}}
\newunicodechar{ï}{\"{\i}}

\newunicodechar{⁰}{\ensuremath{^0}}
\newunicodechar{¹}{\ensuremath{^1}}
\newunicodechar{²}{\ensuremath{^2}}
\newunicodechar{³}{\ensuremath{^3}}
\newunicodechar{⁴}{\ensuremath{^4}}
\newunicodechar{⁵}{\ensuremath{^5}}
\newunicodechar{⁶}{\ensuremath{^6}}
\newunicodechar{⁷}{\ensuremath{^7}}
\newunicodechar{⁸}{\ensuremath{^8}}
\newunicodechar{⁹}{\ensuremath{^9}}
\newunicodechar{ⁱ}{\ensuremath{^i}}

\newunicodechar{⌈}{\ensuremath{\lceil}}
\newunicodechar{⌉}{\ensuremath{\rceil}}
\newunicodechar{⌊}{\ensuremath{\lfloor}}
\newunicodechar{⌋}{\ensuremath{\rfloor}}

\newunicodechar{≅}{\ensuremath{\cong}}
\newunicodechar{⇔}{\ensuremath{\Leftrightarrow}}
\newunicodechar{∃}{\ensuremath{\exists}}
\newunicodechar{±}{\ensuremath{\pm}}

\DeclareFontFamily{OMS}{rsfs}{\skewchar\font'60}
\DeclareFontShape{OMS}{rsfs}{m}{n}{<-5>rsfs5 <5-7>rsfs7 <7->rsfs10 }{}
\DeclareSymbolFont{rsfs}{OMS}{rsfs}{m}{n}
\DeclareSymbolFontAlphabet{\scr}{rsfs}
\DeclareSymbolFontAlphabet{\scr}{rsfs}

\DeclareFontFamily{U}{mathx}{\hyphenchar\font45}
\DeclareFontShape{U}{mathx}{m}{n}{
      <5> <6> <7> <8> <9> <10>
      <10.95> <12> <14.4> <17.28> <20.74> <24.88>
      mathx10
      }{}
\DeclareSymbolFont{mathx}{U}{mathx}{m}{n}
\DeclareFontSubstitution{U}{mathx}{m}{n}
\DeclareMathAccent{\wcheck}{0}{mathx}{"71}

\newunicodechar{∂}{\ensuremath{\partial}}
\newunicodechar{∇}{\ensuremath{\nabla}}

\newunicodechar{↺}{\ensuremath{\circlearrowleft}}
\newunicodechar{∞}{\ensuremath{\infty}}
\newunicodechar{⊕}{\ensuremath{\oplus}}
\newunicodechar{⊗}{\ensuremath{\otimes}}
\newunicodechar{•}{\ensuremath{\bullet}}
\newunicodechar{Λ}{\ensuremath{\wedge}}
\newunicodechar{↪}{\ensuremath{\into}}
\newunicodechar{→}{\ensuremath{\to}}
\newunicodechar{↦}{\ensuremath{\mapsto}}
\newunicodechar{⨯}{\ensuremath{\times}}
\newunicodechar{∪}{\ensuremath{\cup}}
\newunicodechar{∩}{\ensuremath{\cap}}
\newunicodechar{⊋}{\ensuremath{\supsetneq}}
\newunicodechar{⊇}{\ensuremath{\supseteq}}
\newunicodechar{⊃}{\ensuremath{\supset}}
\newunicodechar{⊊}{\ensuremath{\subsetneq}}
\newunicodechar{⊆}{\ensuremath{\subseteq}}
\newunicodechar{⊂}{\ensuremath{\subset}}
\newunicodechar{⊄}{\ensuremath{\not \subset}}
\newunicodechar{≥}{\ensuremath{\geq}}
\newunicodechar{≠}{\ensuremath{\neq}}
\newunicodechar{≫}{\ensuremath{\gg}}
\newunicodechar{≪}{\ensuremath{\ll}}

\newunicodechar{≤}{\ensuremath{\leq}}
\newunicodechar{∈}{\ensuremath{\in}}
\newunicodechar{∉}{\ensuremath{\not \in}}
\newunicodechar{∖}{\ensuremath{\setminus}}
\newunicodechar{◦}{\ensuremath{\circ}}
\newunicodechar{°}{\ensuremath{^\circ}}
\newunicodechar{…}{\ifmmode\mathellipsis\else\textellipsis\fi}
\newunicodechar{·}{\ensuremath{\cdot}}
\newunicodechar{⋯}{\ensuremath{\cdots}}
\newunicodechar{∅}{\ensuremath{\emptyset}}
\newunicodechar{⇒}{\ensuremath{\Rightarrow}}

\newunicodechar{⁰}{\ensuremath{^0}}
\newunicodechar{¹}{\ensuremath{^1}}
\newunicodechar{²}{\ensuremath{^2}}
\newunicodechar{³}{\ensuremath{^3}}
\newunicodechar{⁴}{\ensuremath{^4}}
\newunicodechar{⁵}{\ensuremath{^5}}
\newunicodechar{⁶}{\ensuremath{^6}}
\newunicodechar{⁷}{\ensuremath{^7}}
\newunicodechar{⁸}{\ensuremath{^8}}
\newunicodechar{⁹}{\ensuremath{^9}}
\newunicodechar{ⁱ}{\ensuremath{^i}}

\newunicodechar{⌈}{\ensuremath{\lceil}}
\newunicodechar{⌉}{\ensuremath{\rceil}}
\newunicodechar{⌊}{\ensuremath{\lfloor}}
\newunicodechar{⌋}{\ensuremath{\rfloor}}

\newunicodechar{≅}{\ensuremath{\cong}}
\newunicodechar{⇔}{\ensuremath{\Leftrightarrow}}
\newunicodechar{∃}{\ensuremath{\exists}}
\newunicodechar{±}{\ensuremath{\pm}}

\DeclareFontFamily{OMS}{rsfs}{\skewchar\font'60}
\DeclareFontShape{OMS}{rsfs}{m}{n}{<-5>rsfs5 <5-7>rsfs7 <7->rsfs10 }{}
\DeclareSymbolFont{rsfs}{OMS}{rsfs}{m}{n}
\DeclareSymbolFontAlphabet{\scr}{rsfs}
\DeclareSymbolFontAlphabet{\scr}{rsfs}

\DeclareFontFamily{U}{mathx}{\hyphenchar\font45}
\DeclareFontShape{U}{mathx}{m}{n}{
      <5> <6> <7> <8> <9> <10>
      <10.95> <12> <14.4> <17.28> <20.74> <24.88>
      mathx10
      }{}
\DeclareSymbolFont{mathx}{U}{mathx}{m}{n}
\DeclareFontSubstitution{U}{mathx}{m}{n}
\DeclareMathAccent{\wcheck}{0}{mathx}{"71}

\theoremstyle{plain}
\newtheorem{Th}{Théorème}[section]

\newtheorem{cor}[Th]{Corollaire}
\newtheorem{defn}[Th]{Définition}

\newtheorem{lem}[Th]{Lemme}

\newtheorem{Prop}[Th]{Proposition}

\theoremstyle{remark}

\newtheorem{c-n-d}[Th]{Claim and Definition}

\newtheorem{rem}[Th]{\color{blue}{Remarques}}

\newtheorem*{rem-nonumber}{Remark}

\numberwithin{equation}{Th}
\setlist[enumerate]{label=(\thethm.\arabic*), before={\setcounter{enumi}{\value{equation}}}, after={\setcounter{equation}{\value{enumi}}}}
\newcommand{\into}{\hookrightarrow}

\newcommand{\factor}[2]{\left. \raise 2pt\hbox{$#1$} \right/\hskip -2pt\raise -2pt\hbox{$#2$}}

\newcommand{\Publication}[1]{}
\newcommand{\subversionInfo}{}
\newcommand{\svnid}[1]{}
\newcommand{\approvals}[2][Approval]{}
\renewcommand{\phi}{\varphi}
\usepackage{tikz}
\usepackage{tikz-cd}
\tikzset{commutative diagrams/arrow style=Latin Modern}
\author{Mohamed Kaddar} %
\email{\href{mailto:mohamed.kaddar@univ-lorraine.fr}{mohamed.kaddar@univ-lorraine.fr}} %
\title[Sur certains faisceaux de formes méromorphes en géométrie analytique complexe I.]{Sur certains faisceaux de formes méromorphes en géométrie analytique complexe I.}%

\date{\today}
\makeatletter
\hypersetup{
  pdfauthor={\authors},
  pdftitle={\@title},
  pdfsubject={\@subjclass},
  pdfkeywords={\@keywords},
  pdfstartview={Fit},
  pdfpagelayout={TwoColumnRight},
  pdfpagemode={UseOutlines},
  bookmarks,
  colorlinks,
  linkcolor=linkblue,
  citecolor=linkred,
  urlcolor=linkred}
\makeatother
\newcommand{\chapref}[1]{\hyperref[#1]{Chapter~\ref*{#1}}}
\newcommand{\lemmaref}[1]{\hyperref[#1]{Lemme~\ref*{#1}}}
\newcommand{\parref}[1]{\hyperref[#1]{Section~\ref*{#1}}}
\newcommand{\theoremref}[1]{\hyperref[#1]{Théorème~\ref*{#1}}}
\newcommand{\definitionref}[1]{\hyperref[#1]{Définition~\ref*{#1}}}
\newcommand{\propositionref}[1]{\hyperref[#1]{Proposition~\ref*{#1}}}
\newcommand{\conjectureref}[1]{\hyperref[#1]{Conjecture~\ref*{#1}}}
\newcommand{\corollaryref}[1]{\hyperref[#1]{Corollaire~\ref*{#1}}}
\newcommand{\exampleref}[1]{\hyperref[#1]{Exemple~\ref*{#1}}}
\newcommand{\exerciseref}[1]{\hyperref[#1]{Exercise~\ref*{#1}}}
\newcommand{\factref}[1]{\hyperref[#1]{Fact~\ref*{#1}}}
\newcommand{\claimref}[1]{\hyperref[#1]{Claim~\ref*{#1}}}
\newcommand{\remarkref}[1]{\hyperref[#1]{Remark~\ref*{#1}}}
\newcommand{\settingref}[1]{\hyperref[#1]{Setting~\ref*{#1}}}
\newcommand{\appendixref}[1]{\hyperref[#1]{Appendix~\ref*{#1}}}

\theoremstyle{plain}

\theoremstyle{remark}

\newcommand{\eq}[1][r]
   {\ar@<-3pt>@{-}[#1]
    \ar@<-1pt>@{}[#1]|<{}="gauche"
    \ar@<+0pt>@{}[#1]|-{}="milieu"
    \ar@<+1pt>@{}[#1]|>{}="droite"
    \ar@/^2pt/@{-}"gauche";"milieu"
    \ar@/_2pt/@{-}"milieu";"droite"}

\begin{document}

\maketitle
\approvals[Approval for Abstract]{Mohamed & yes}
\begin{abstract}
We study some functorial properties of certain sheaves of meromorphic forms on reduced complex space. Particulary, their behavior under pull back and higher direct image (and in some case by integration  on fibres of equidimensional or open morphism).
\end{abstract}
\maketitle
\tableofcontents

\noindent
\phantomsection\addcontentsline{toc}{part}{Introduction}
%
%

Le présent travail a pour objectif l'étude, dans le cas absolu, de certains faisceaux de formes méromorphes sur un espace analytique complexe réduit; le cas relatif sera traité dans un article indépendant pour ne pas alourdir  la teneur du texte. Le thème principal  tourne autour de la notion de formes méromorphes régulières dont les origines se trouvent dans la théorie de la  dualité dans la catégorie des faisceaux cohérents. On sait qu'en présence de singularités, bon nombre de résultats classiques de la géométrie algébrique ou analytique ne sont jamais vérifiés; en particulier, l'exactitude de certains complexes de formes holomorphes, d'annulation de cohomologie ou d'homologie et bien d'autres exemples. Un des problèmes majeurs était de trouver, dans le cadre singulier,  un substitut au faisceau des formes holomorphes usuelles de degré maximal qui joue un rôle crucial dans la théorie de la dualité  de type Serre ou autre. Malheureusement (ou heureusement) sur un espace complexe à singularités arbitraires, il n'y'avait aucun espoir de mettre en évidence un substitut naïf permettant d'élaborer une dualité cohérente généralisant celle connue pour une variété lisse hormis pour certains espaces soumis à des contraintes draconiennes imposant aux singularités d'être d'une nature très particulière ( rationnelles, quotients, canoniques ...) et dont nous aurons tout le losir d'en parler dans \cite{K6}. \vspace{1mm}

\noindent  A notre connaissance, c'est à Grothendieck que l'on doit la bonne formulation du problème et les premières profondes investigations dans ce sens. D'ailleurs, la dualité topologique de Verdier montre qu'il faut penser la question en termes de complexes défini dans la catégorie dérivée. Il est appelé, à juste titre, complexe dualisant et possède un certain nombre de propriétés fonctorielles en plus, d'être quasi isomorphe, dans le cas lisse, au complexe constitué du faisceau constant muni de son décalage naturel. \vspace{1mm}

\noindent En ce qui nous concerne,  le complexe dualisant, défini dans la catégorie cohérente, est un complexe  de faisceaux cohérents à cohomologie cohérente construit par Ramis et Ruget dans  \cite{RR70}. Nous entrerons dans les détails dans la partie concernant la dualité relative de type Kleiman.  \vspace{1mm}

\noindent Notre propos porte sur un faisceau particulier qui, par essence, possède quelques propriétés dualisantes puisqu'apparaissant dans certains cas comme l'homologie en degré maximal du complexe dualisant. Si l'espace complexe est de dimension pure $m$, ce faisceau noté $\omega^{m}_{X}$ est le faisceau dualisant de Grothendieck dont on connait parfaitement les propriétés. Comme la dualité ne s'appréhende que difficilement quelque soit le cadre envisagé, il a fallu dégager une propriété carcatéristique plus géométrique pour pouvoir le manipuler avec aisance. Les investigations menées par plusieurs auteurs ont permis de le caractériser par ce que l'on appelle communément, maintenant, la {\it{propriété de la trace}} et le désignant par la suite comme le faisceau des {\it{formes méromorphes régulières}}; nous  renvoyons le lecteur à \cite{L} pour une longue et précise discussion dans le cadre algébrique. 
 On a choisi de rester volontairement vague puisque tout ce qui a été dit va être développé  dans les moindres détails dans \cite{K5} et \cite{K8}.
\vspace{1mm}

\noindent
Pour un espace analytique complexe $X$ de dimension pure $m$, le faisceau des $k$-formes régulières $\omega^{k}_{X}$  caractérisé par la {\emph{propriété de la trace}}, doté d'un certain nombre de propriétés fonctorielles  intéressantes mais présente, malheureusement, un  défaut majeur celui de ne pas être  stable par image réciproque (et pas de restriction au sens usuel!). 
Ils possède néanmoins  des sous faisceaux cette propriété pour une certaine classe de morphismes (par exemple, celle des  morphismes surjectifs et à fibres de dimension constante).  La construction de ces pull-back se veut être une généralisation de ceux des formes holomorphes usuelles et doit donc, tant qu'il se peut, préserver toutes les propriétés fonctorielles de l'image réciproque  usuelle. En règle générale, il est pratiquement impossible de réaliser un tel objectif  si ce n'est au prix de quelques concessions  (et pas des moindres!) tant sur la nature des morphismes que sur celle des espaces eux-mêmes (et donc la nature de leurs singularités).\vspace{1mm}

\noindent En désignant par  ${\mathcal L}^{k}_{X}$  le sous faisceau de $\omega^{k}_{X}$ dont les sections sont des $k$-formes méromorphes sur $X$, se prolongent holomorphiquement sur toute désingularisation de $X$, il est possible de construire une image réciproque pour tout morphisme d'espaces analytiques réduits et  tout sous faisceau de ${\mathcal L}^{k}_{X}$ (ce dernier y compris) (cf\cite{K4}, \cite{K0}). Malheureusement, elle n'est pas généralement pas compatible avec la composition des morphismes à moins de travailler avec des morphismes surjectifs permettant un raisonnement générique (c-à-d sur les parties régulières). \vspace{1mm}

\noindent Ces faisceaux   ${\mathcal L}^{\bullet}_{X}$ ont pratiquement toutes les propriétés des faisceaux des formes holomorphes $\Omega^{\bullet}_{X}$ à savoir ceux sont des modules sur ces derniers stables par cup produit, par différentiation extérieure, par image réciproque, par trace mais ne sont, en général, jamais auto-duaux. En un certain sens, ils jouent un  rôle de médiateur entre les sous faisceaux de $\omega^{\bullet}_{X}$ qui sont stables par trace  et ceux qu'ils le sont par image réciproque (même restrictive!). On s'apperçoit que plus on s'éloigne de ${\mathcal L}^{\bullet}_{X}$ en allant vers $\omega^{\bullet}_{X}$, plus les faisceaux rencontrés perdent en stabilité par image réciproque et gagne en stabilité par trace alors que le phénomène inverse se produit quand on s'éloigne de ${\mathcal L}^{\bullet}_{X}$ en allant vers le faisceau $\Omega^{\bullet}_{X}$.  \vspace{2mm}

\indent 
Signalons que dans \cite{Keb}, la stabilité par image réciproque des faisceaux $\omega^{k}_{X}$ est établie pour tout morphisme d'espaces complexes normaux sur une base $X$ de type Kawamata log-terminale. Ce pull back vérifie toutes les propriétés attendues. \vspace{1mm}

\noindent  Un autre aspect de notre étude concerne le  comportement de ces faisceaux par image directe supérieure (propre ou non). Plus précisément, on associe à tout morphisme $\pi:X\rightarrow S$ équidimensionnel ou ouvert à fibres de dimension pure $n$ d'espaces complexes de dimension pure $m$ et $r$ respectivement   et  un faisceau   $\widetilde{\mathcal F}^{n+q}_{X}$ (resp. $\widetilde{\mathcal F}^{q}_{S}$ s'identifiant à l'un des faisceaux  ${\mathcal L}^{n+q}_{X}$, ${\mathcal H}om({\mathcal L}^{r-q}_{X}, {\mathcal L}^{m}_{X})$,  ${\mathcal H}om({\Omega}^{r-q}_{X}, {\mathcal L}^{m}_{X})$, ${\mathcal H}om({\mathcal L}^{r-q}_{X}, \omega^{m}_{X})$, ${\mathcal H}om({\Omega}^{r-q}_{X}, \omega^{m}_{X})$ (resp. leurs analogues sur $S$)  un morphisme  $\displaystyle{{\mathcal T}^{q}_{\pi}:{\rm I}\!{\rm R}^{n}\pi_{!}\widetilde{\mathcal F}^{n+q}_{X}\rightarrow\widetilde{\mathcal F}^{q}_{S}}$ 
(en remplaçant $\pi_{!}$ par $\pi_{*}$ dans le cas propre) compatibles aux restrictions ouvertes sur $X$ et sur $S$ et aux changements de base:\vspace{1mm}

\indent $\star$ arbitraires (resp. surjectif) si ${\mathcal F}^{n+q}_{X}:={\mathcal L}^{n+q}_{X}$  (resp. $\overline{\Omega}^{n+q}_{X}$),\vspace{1mm}

\indent $\star$  préservant l'équidimensionnalité et la pureté des dimensions pour ${\mathcal F}^{n+q}_{X}:={\mathcal H}om({\mathcal L}^{r-q}_{X}, \omega^{m}_{X})$ ou ${\mathcal H}om({\mathcal L}^{r-q}_{X}, {\mathcal L}^{m}_{X})$,\vspace{1mm}

\indent $\star$  entre espaces complexes de lieu singulier de codimension au moins deux pour  ${\mathcal F}^{n+q}_{X}=\omega^{n+r}_{X}$. \vspace{2mm}

\noindent 
Ces morphismes sont induits par les morphismes de degré maximaux  $\displaystyle{{\mathcal T}^{r}_{\pi,\omega}:{\rm I}\!{\rm R}^{n}\pi_{!}{\omega}^{n+r}_{X}\rightarrow{\omega}^{r}_{S}}$,
$\displaystyle{{\mathcal T}^{r}_{\pi, {\mathcal L}}:{\rm I}\!{\rm R}^{n}\pi_{!}{\mathcal L}^{n+r}_{X}\rightarrow{\mathcal L}^{r}_{S}}$ dont l'existence peut être établie de façon conceptuelle ou par construction explicite à l'aide des techniques d'intégration. \vspace{1mm}

\noindent 
Pour $q\!=\!0$ et $S$ normal, on obtient  les  flèches $\displaystyle{{\mathcal T}^{0}_{\pi}:{\rm I}\!{\rm R}^{n}\pi_{!}{\mathcal F}^{n}_{X}\rightarrow{\mathcal O}_{S}}$ et, en particulier, 
$\displaystyle{{\mathcal T}^{0}_{\pi}:{\rm I}\!{\rm R}^{n}\pi_{!}{\omega}^{n}_{X}\rightarrow{\mathcal O}_{S}}$ que l'on peut interpréter  
 en terme d'intégration de classes de cohomologie méromorphes. \vspace{2mm}

\noindent
Pour compléter l'étude, nous construisons un morphisme 
$$\displaystyle{{\mathcal T}^{r}_{\pi, \overline{\Omega}}:{\rm I}\!{\rm R}^{n}\pi_{!}\overline{\Omega}^{n+r}_{X}\rightarrow\overline{\Omega}^{r}_{S}}$$
ou de façon plus générale
$$\displaystyle{{\mathcal T}^{q}_{\pi, \overline{\Omega}}:{\rm I}\!{\rm R}^{n}\pi_{!}\overline{\Omega}^{n+q}_{X}\rightarrow\overline{\Omega}^{q}_{S}}$$
compatible à une certaine classe (suffisament courante!) de changement de base.\vspace{1mm}

\noindent 
Le faisceau $\overline{\Omega}^{j}_{X}$ est un sous faisceau $\omega^{j}_{X}$ dont 
les sections sections vérifient des équations de dépendance intégrale sur l'algèbre symétrique des formes holomorphes introduit dans \cite{Sp} et étudié en détail dabs \cite{B5}, \cite{B6}. Pour $q=0$, ce n'est rien d'autre que le faisceau des formes faiblement holomorphes ($\overline{\Omega}^{0}_{X}={\mathcal L}^{0}_{X}$) et, pour $q=1$, ceux sont les sections de $\omega^{1}_{X}$ qui correspondent aux fonctions faiblement holomorphes sur le tangent de Zariski de $X$.
\vspace{1mm}

\noindent Dans ce qui suit, on convient des notations suivantes auxquelles on renvoie le lecteur au {\bf{\S 2.5}}:\vspace{1mm}

\indent
$\bullet$ ${\mathcal D}_{\omega_{X}}({\mathcal F}):={\mathcal H}om({\mathcal F}, \omega^{m}_{X})$ (resp. ${\mathcal D}_{{\mathcal L}_{X}}({\mathcal F}):={\mathcal H}om({\mathcal F}, {\mathcal L}^{m}_{X})$) le ${\omega}^{m}_{X}$ (resp. ${\mathcal L}^{m}_{X}$)-dual d'un faisceau cohérent  ${\mathcal F}$,  \vspace{1mm}

\indent
$\bullet$ $\overline{\Omega}^{q}_{X}:={\nu_{q}}_{*}({\nu^{*}_{q}}(\widetilde{\Omega}^{q}_{X})/{\mathcal T}_{\nu_q})$, le faisceau associé à toute modification propre 
${\nu_{q}}:\bar{X}\rightarrow X$ dans laquelle $\bar{X}$ est normal et pour laquelle le faisceau ${\nu_{q}}^{*}(\widetilde{\Omega}^{q}_{X})$ ($\widetilde{\Omega}^{q}_{X}$ qui est le faisceau des formes holomorphes usuelles modulo torsion) quotienté par sa torsion ${\mathcal T}_{\nu_q}$ est localement libre. Nous avons conservé la notation de \cite{Sp} alors que dans l'étude détaillée  de \cite{B5}, \cite{B5}, il est noté $\alpha^{\bullet}_{X}$.\vspace{1mm}

\noindent

On a, alors, 
\Th{}{}\label{T1} 
Soient $X$ et $S$ deux espaces analytiques complexes de dimension pure $m$ et $r$ respectivement, $n$ un entier naturel et $\pi:X\rightarrow S$  un morphisme surjectif $n$-équidimensionnel ( ou ouvert à fibres de dimension constante $n$).  Alors: \vspace{1mm}

\noindent
{\bf(i)} si ${\mathcal F}^{q}_{S}=\omega^{q}_{S}, {\mathcal L}^{q}_{S}$, il existe  un unique morphisme ${\mathcal O}_{X}$ linéaire  de faisceaux cohérents  $${\bf{\pi}}^{*}: \pi^{*}({\mathcal D}_{\omega_{_S}}({\mathcal F}^{q}_{S}))\rightarrow{\mathcal D}_{\omega_{_X}}({\mathcal F}^{q}_{X})$$
(et même d'algèbre différentielle sur ${\Omega}^{\bullet}_{S}$) prolongeant l'image réciproque des formes holomorphes usuelles sur la partie régulière, rendant commutatif le diagramme
$$\xymatrix{{\Omega}^{q}_{S}\ar[r]\ar[d]&\pi_{*}{\Omega}^{q}_{X}\ar[d]\\
{\mathcal D}_{\omega_{_S}}({\mathcal F}^{q}_{S})\ar[r]&\pi_{*}{\mathcal D}_{\omega_{_X}}({\mathcal F}^{q}_{X})}$$
compatible avec la différentielle extérieure et avec  la composition des morphismes équidimensionnels.\vspace{1mm}

\noindent
{\bf(ii)} On a aussi une image réciproque 
$${\bf{\pi}}^{*}: \pi^{*}({\mathcal D}_{{\mathcal L}_{_S}}({\mathcal L}^{q}_{S}))\rightarrow{\mathcal D}_{{\mathcal L}_{_X}}({\mathcal L}^{q}_{X})$$
prolongeant l'image réciproque des formes holomorphes usuelles sur la partie régulière et compatible avec  la composition des morphismes équidimensionnels.\vspace{1mm}

\noindent 
{\bf(iii)} si ${\mathcal F}^{q}_{X}=\overline{\Omega}^{q}_{X}, {\mathcal L}^{q}_{X}$, tout morphisme fini et surjectif $f:X\rightarrow Y$ induit un unique morphisme de faisceaux cohérents $$f_{*}({\mathcal D}_{\omega_{_X}}({\mathcal F}^{q}_{X}))\rightarrow{\mathcal D}_{\omega_{_{Y}}}({\mathcal F}^{q}_{Y})$$
et 
$$f_{*}{\mathcal D}_{{\mathcal L}_{_X}}({\mathcal L}^{q}_{X})\rightarrow {\mathcal D}_{{\mathcal L}_{_Y}}({\mathcal L}^{q}_{Y})$$
 \rm\vspace{1mm}

\noindent
\cor{}{}\label{Cor0'} Avec les notations et hypothèses du \theoremref{T1} et $\pi:X\rightarrow S$ un morphisme surjectif à fibres de dimension constante $n$, alors, s'il existe  un morphisme de faisceaux 
${\bf{\pi}}^{*}:{\omega}^{r}_{S}\rightarrow \pi_{*}{\omega}^{r}_{X}$  prolongeant l'image réciproque des formes holomorphes usuelles et si ${\mathcal O}_{S}\simeq {\mathcal H}om(\omega^{r}_{S}, \omega^{r}_{S})$,  $\pi$ est analytiquement géométriquement plat. En particulier, tout morphisme ouvert à fibres de dimension constante $n$ sur une base de Cohen Macaulay et ayant cette propriété du pull-back est analytiquement géométriquement plat. \rm
\vspace{1mm}

\noindent 
 Comme nous ne manipulons que des morphisme analytiquement géométriquement plats (c-à-d  définissant des familles analytiques de cycles), il nous arrivera d'omettre l'adverbe "analytiquement" puisqu'il n'y'a aucune confusion à craindre avec les autres notions (cf \cite{K1}, p.39).  
\Th{}{}\label{T'1} Soit $\pi:X\rightarrow S$  un morphisme surjectif $n$-équidimensionnel ou ouvert à fibres de dimension constante $n$ d'espaces complexes réduits. On a, pour tout entier $q$,  si ${\mathcal F}^{r-q}_{X}$ (resp. ${\mathcal F}^{r-q}_{S}$) est l'un des faisceaux   $\overline{\Omega}^{r-q}_{X}$, $\widehat{\Omega}^{r-q}_{X}$, ${\mathcal L}^{r-q}_{X}$ (resp. leurs analogues sur $S$), des morphismes canoniques  de faisceaux de ${\mathcal O}_S$-modules (cohérents si $\pi$ est propre ):\vspace{2mm}

\indent
$\bullet$  ${\mathcal T}^{q}_{\pi}:{\rm I}\!{\rm R}^{n}\pi_{!}{\mathcal D}_{\omega_{_X}}({\mathcal F}^{r-q}_{X})\rightarrow{\mathcal D}_{\omega_{_S}}({\mathcal F}^{r-q}_{S})$ (en remplaçant $\pi_{!}$ par $\pi_{*}$ dans le cas propre) et \vspace{1mm}

\indent
$\bullet$  $\widetilde{\mathcal T}^{q}_{\pi}:{\rm I}\!{\rm R}^{n}\pi_{*}{\mathcal D}_{{\mathcal L}_{_X}}({\mathcal F}^{r-q}_{X})\rightarrow{\mathcal D}_{{\mathcal L}_{S}}({\mathcal F}^{r-q}_{S})$ (en remplaçant $\pi_{!}$ par $\pi_{*}$ dans le cas propre). 
 \vspace{2mm}

\noindent   Ils sont de formation  compatible aux restrictions ouvertes sur $X$ et $S$, aux changements de base donnés par les morphismes surjectifs à fibres de dimension constante pour ${\mathcal F}^{r-q}_{X}={\mathcal L}^{r-q}_{X}$. De plus,  ils sont compatibles à la composition des morphismes dans le sens suivant:\vspace{ 1mm}

\noindent
Pour tout  diagramme commutatif d'espaces analytiques complexes
$$\xymatrix{X_{2}\ar[rr]^{\Psi}\ar[rd]_{\pi_{2}}&&X_{1}\ar[ld]^{\pi_{1}}\\
&S&}$$
 avec $\pi_{1}$ (resp. $\pi_{2}$)  ouvert de dimension relative $n_{1}$ (resp. $n_{2}$)  et $\Psi$ propre de dimension relative bornée par l'entier   $d:=n_{2}-n_{1}$,
 muni d'un morphisme canonique $\displaystyle{{\rm I}\!{\rm
 R}^{d}{{\Psi}_{*}}{\mathcal D}_{\omega_{X_{2}}}({\mathcal F}^{r-q}_{X_{2}})\rightarrow{\mathcal D}_{\omega_{X_{1}}}({\mathcal F}^{r-q}_{X_{1}})}$, on a un diagramme   commutatif de faisceaux analytiques (cohérents dans le cas propre)
$$\xymatrix{{\rm I}\!{\rm
 R}^{n_{2}}{\pi_{2}}_{!}{\mathcal D}_{\omega_{X_{2}}}({\mathcal F}^{r-q}_{X_{2}})\ar[rr]\ar[rd]_{{\mathcal T}^{q}_{\pi_{2}}}&&
 {\rm I}\!{\rm
 R}^{n_{1}}{\pi_{1}}_{!}{\mathcal D}_{\omega_{X_{1}}}({\mathcal F}^{r-q}_{X_{1}})\ar[ld]^{{\mathcal T}^{q}_{\pi_{1}}}\\
&{\mathcal D}_{\omega_{S}}({\mathcal F}^{r-q}_{S})&}$$
et un diagramme analogue pour les morphismes $\widetilde{\mathcal T}^{q}_{\pi}$, en remplaçant $\pi_{!}$ par $\pi_{*}$ dans le cas propre.\rm
\cor{}{}\label{Cor2}  Soit   $\pi:X\rightarrow S$ un morphisme propre  $n$-équidimensionnel (resp. universellement $n$-équidimensionnel) d'espaces complexes réduits de dimension pure. Alors, pour tout faisceau ${\mathcal O}_{X}$ cohérent ${\mathcal F}$, le morphisme canonique 
$${\rm I}\!{\rm Hom}(X; {\mathcal F},  \omega^{n+r}_{X})\rightarrow{\rm I}\!{\rm Hom}(X; {\rm I}\!{\rm R}^{n}\pi_{*}{\mathcal F}, \omega^{r}_{S})$$
est un isomorphisme de formation fonctorielle en ${\mathcal F}$, compatible aux changements de base à fibres de dimension constante et aux images directes propres.\rm

\noindent
\Th{}{}\label{T2} Soit $\pi:X\rightarrow S$ un morphisme géométriquement plat à fibres de dimension pure $n$ entre espaces complexes réduits de dimension $m$ et $r$ respectivement. Alors, il existe un morphisme de ${\mathcal O}_{S}$-module (cohérents si le morphisme est propre)
$${\mathcal T}^{j}_{\pi,\overline{\Omega}}: {\rm I}\!{\rm R}^{n}\pi_{!}\overline{\Omega}^{n+j}_{X}\rightarrow\overline{\Omega}^{j}_{S}$$ 
( en remplaçant $\pi_{!}$ par $\pi_{*}$ dans le cas propre) compatible: \vspace{1mm}

\noindent {\bf(i)} aux morphismes d'intégration du \theoremref{T1} de sorte que les diagrammes
$$\xymatrix{{\rm I}\!{\rm R}^{n}{\pi_{!}}{\overline{\Omega}^{n+q}_{X}}\ar[d]_{{\mathcal T}^{j}_{\pi,\overline{\Omega} }}\ar[r]&{\rm I}\!{\rm R}^{n}{\pi_{!}}{\mathcal L}^{n+q}_{X}\ar[d]_{{\mathcal T}^{j}_{\pi,{\mathcal L}}} \ar@{=}[rr]&&{\rm I}\!{\rm R}^{n}{\pi_{!}}{\mathcal L}^{n+q}_{X}\ar[r]\ar[d]_{{\mathcal T}^{j}_{\pi,{\mathcal L}}}&{\rm I}\!{\rm R}^{n}{\pi_{!}}{\omega}^{n+q}_{X}\ar[d]_{{\mathcal T}^{j}_{\pi,\omega}}\\
{\Omega}^{q}_{S}\ar[r]&{\mathcal L}^{q}_{S}\ar@{=}[rr]&&{\mathcal L}^{q}_{S}\ar[r]&{\omega}^{q}_{S}}$$
soient commutatifs, \vspace{1mm}

\noindent 
{\bf(ii)}  aux  changements de base admissibles sur $S$ c'est-à-dire que pour tout diagramme de changement de base
$$\xymatrix{\widetilde{X}\ar[d]_{\tilde\pi}\ar[r]^{\theta}&X\ar[d]^{\pi}\\
\widetilde{S}\ar[r]_{\nu}&S}$$
avec $\nu^{-1}(S)\subsetneq{\rm Sing}({\widetilde S})$, 
induit un diagramme commutatif
$$\xymatrix{\nu^{*}{\rm I}\!{\rm R}^{n}{\pi_{!}}{\overline{\Omega}^{n+q}_{X}}\ar[d]_{\nu^{*}({\mathcal T}^{j}_{\pi,\overline{\Omega} })}\ar[r]&{\rm I}\!{\rm R}^{n}{\tilde\pi_{!}}\theta^{*}{\overline{\Omega}^{n+q}_{X}}\ar[r]&{\rm I}\!{\rm R}^{n}{\tilde\pi_{!}}{\overline{\Omega}^{n+q}_{\widetilde X}}\ar[d]^{{\mathcal T}^{j}_{\tilde\pi,\overline{\Omega}}}\\
\nu^{*}{\overline{\Omega}}^{q}_{S}\ar[rr]&&{\overline{\Omega}}^{q}_{\widetilde S}}$$
{\bf(iii)} Soient $Z$ et $S$  deux espaces  analytiques complexes de dimension pure avec 
$Z$ dénombrable 
à l'infini et $S$ réduit. Soit $(X_{s})_{s\in S}$ une famille analytique de $n$- cycles de $Z$ 
paramétrée par $S$ dont le graphe $X$ est supposé non entièrement inclus dans le lieu singulier de $S\times Z$ et muni du morphisme $\pi: X\rightarrow S$ (induit par la projection canonique sur $S$). Soient $\Phi$ et  $\Psi$ deux familles de supports avec $\Phi$ paracompactifiante et $\Phi \cap \Psi$ contenue dans la famille des compacts de $Z$, $\forall s\in S$, $X_{s} \in \Psi$ et  $X\cap (S\times \Phi)$ contenu dans la famille des ferm\'es $S$- 
propres.\vspace{1mm}

\noindent  
 Alors il existe
un morphisme d'int\'egration d'ordre sup\'erieur sur les cycles 
$$\overline{\sigma}^{q,0}_{{\Phi},X}:{\rm H}^{n}_{\Phi}(Z, \overline{\Omega}^{n+q}_{Z})\longrightarrow 
{\rm H}^{0}(S, \overline{\Omega}^{q}_{S})$$
\noindent v\'erifiant les propri\'et\'es énoncées dans le {\bf{Thm 1.3}} de \cite{K0}, \cite{K4}.\rm 

\cor{}{}\label{C0} On conserve  les notations et  hypothèses du \theoremref{T2}.
Soient $\Phi$ une famille paracompactifiante de supports de $Z$, $\Psi$ une famille de supports de $Z$ contenant les supports des cycles $X_s$,  $\Theta$  une famille paracompactifiante de supports de $S$  telle que la famille $\widetilde{\Psi}:=(S\times (\Phi\cap \Psi))\cap X$ soit paracompactifiante et contenue dans la famille paracompactifiante des fermés $\Theta$-propre de sorte  que le couple $(\widetilde{\Psi}, \Theta)$ soit adapté à $\pi$. Alors, il existe
un morphisme d'int\'egration d'ordre sup\'erieur sur les cycles 
$$\overline{\sigma}^{q,p}_{{\Phi},X}:{\rm H}^{n+p}_{\Phi}(Z, \overline{\Omega}^{n+q}_{Z})\longrightarrow 
{\rm H}^{p}_{\Theta}(S, \overline{\Omega}^{q}_{S})$$
possédant les mêmes propriétés fonctorielles que $\overline{\sigma}^{q,0}_{{\Phi},X}$.\rm

\vfill\eject

\phantomsection\addcontentsline{toc}{part}{Formes méromorphes régulières absolues.}
%
%
\svnid{$Id: S07-proof-setup.tex 269 2020-01-20 11:28:53Z kebekus $}

\section{\color{blue}{Formes méromorphes régulières absolues.}}
\subversionInfo

\subsection{\color{blue}{Préliminaires.}}
\approvals{Mohamed & yes}
\par\vspace{2mm}

\noindent 
\subsubsection{\bf{Formes holomorphes et torsion}.} 
\approvals{Mohamed & yes} 
Rappelons rapidement la définition de Ferrari [F] des formes holomorphes sans torsion.\vspace{1mm}

\noindent
Si $X$ est un sous ensemble analytique d'un ouvert $V$ d'un espace numérique 
${\Bbb C}^{m}$, on note $X^0:=X$, $X^{1}:={\rm Sing}(X^0)$,...., $X^{j+1}:={\rm Sing}(X^j)$ et 
$${\mathcal H}^{p}:=\{\alpha\in {\Omega^{p}_{V}}:\alpha\!\!\mid_{{\rm Reg}(X)}=0\}$$ 
$${\mathcal H}^{p}_{1}:={\bigcap_{j}}\{\alpha\in\Omega^{p}_{V}:\alpha\!\!\mid_{{\rm Reg}(X^{j})}=0\}$$

Si $\frak{M}$ est l'ensemble de tous les morphismes $f:W\rightarrow V$ avec $W$ variété lisse telle que $f(W)\subset X$, 
$${\mathcal H}^{p}_{2}:=\bigcap_{f\in\frak{M}}\{\alpha\in {\Omega^{p}_{V}}: \,f^{*}(\alpha)=0\}$$
Alors, d'après ([F], {\bf lemma(1.1), p.67}:
$${\mathcal H}^{p}={\mathcal H}^{p}_{1}={\mathcal H}^{p}_{2}$$
De plus, le faisceau $\widetilde{\Omega}^{p}_{X}:={{\Omega}^{p}_{V}}/{\mathcal H}^{p}$ est lié au faisceau des formes holomorphes au sens de Grothendieck par la suite exactes courte
$$\xymatrix{0\ar[r]&{\mathcal K}^p\ar[r]&\Omega^{p}_{X}\ar[r]&\widetilde{\Omega}^{p}_{X}\ar[r]&0}$$
et aux faisceaux $\omega^{\bullet}_{X}$ et ${\mathcal L}^{p}_{X}$ par le diagramme
de suite exacte courtes
$$\xymatrix{0\ar[r]&\widetilde{\Omega}^{p}_{X}\ar@{=}[d]\ar[r]&{\mathcal L}^{p}_{X}\ar[r]\ar[d]&{\mathcal K'}^p\ar[r]\ar[d]&0\\
0\ar[r]&\widetilde{\Omega}^{p}_{X}\ar[r]&\omega^{p}_{X}\ar[r]&{\mathcal K"}^p\ar[r]&0}$$
${\mathcal K}^{p}$, ${\mathcal K'}^{p}$ et ${\mathcal K"}^p$ étant des faisceaux de torsion à support dans le lieu singulier (d'ailleurs $K^p$ est le sous faisceau de torsion de  ${\Omega}^{p}_{X}$).\vspace{1mm}

\noindent
Signalons, au passage, qu'un morphisme de faisceaux cohérents envoie toujours la torsion de l'un sur celle de l'autre.\vspace{1mm}

\noindent
\subsubsection{\bf{Module dualisant}.}\vspace{1mm}

\noindent
\rm Rappelons qu' un {\emph{module dualisant}}, sur une variété algébrique projective $X$ de dimension $n$
définie sur un corps algébriquement clos de caractéristique
nulle $\bf{k}$ (resp. un espace analytique complexe compact  $X$
 de dimension $n$), est la donnée d'un ${\mathcal O}_{X}$- module
cohérent ${\mathcal K}_{X}$ et d'un  morphisme $\bf{k}$-linéaire
appelé {\emph{trace}} $\int_{X}: {\rm H}^{n}(X, {\mathcal
K}_{X})\rightarrow {\bf{k}}\,\, ({\rm resp.} {\Bbb C})$ assurant
l'isomorphisme de foncteurs de la catégorie des ${\mathcal O}_{X}$-
modules cohérents sur la catégorie des $\bf{k}$ (resp. ${\Bbb
C}$) -espaces vectoriels donné par $${\rm H}om({\mathcal F}, {\mathcal
K}_{X})\simeq {\rm H}^{n}(X, {\mathcal F})^{*}={\rm H}om({\rm H}^{n}(X,
{\mathcal F}), {\bf{k}})$$ Le couple $({\mathcal K}_{X} , \int_{X})$,
appelé {\emph{ paire dualisante}} (unique à isomorphisme canonique
près), représente le foncteur ${\mathcal F}\rightarrow {\rm
H}^{n}(X, {\mathcal F})^{*}$ et $\int_{X}$ correspond à l'identité
dans ${\rm H}om({\mathcal K}_{X}, {\mathcal K}_{X})$.
\vspace{1mm}

\noindent
{\bf Grothendieck} donne dans \cite{G1} montre comment associer à toute
variété projective (pas nécessairement normale) un faisceau dualisant ${\mathcal K}_{X}$ qu'il note $\omega^{n}_{X}$ et qu'il décrit comme étant  l'unique faisceau cohérent de profondeur au moins deux sur $X$ qui, pour tout plongement $\sigma$, de $X$ dans un espace projectif ${\Bbb P}_{N}$, coincide avec  $\sigma^{*}{\mathcal
E}xt^{N-n}_{{\mathcal O}_{{\Bbb P}_{N}}}({\mathcal O}_{X}, \Omega^{N}_{{\Bbb
P}_{N}})$. Cette construction s'étend sans aucune difficultés au cas
de la geométrie analytique complexe comme on peut le voir, par exemple, en consultant le texte clair et  instructif de Monique
Lejeune Jalabert (\cite{LJ}).\vspace{1mm}

\noindent  
\subsubsection{\bf{Faisceaux dualisants de Andréotti-Kas-Golovin}.}\vspace{2mm}

\noindent
Rappelons que Ramis et Ruget dans  \cite{RR70} ont construit, pour tout espaces analytique complexe $X$, un complexe  ${\mathcal D}^{\bullet}_{X}$ de ${\mathcal O}_{X}$-modules cohérents   vérifiant les propriétés suivantes:\vspace{1mm}

\noindent
{\bf(i)} Si $Z$ est une variété lisse de dimension $m$, les fibres de ses composantes, en chaque point $x$ de $X$, sont des  ${\mathcal O}_{X,x}$-modules injectifs et ${\mathcal D}^{\bullet}_{X}$ est une résolution de $\Omega^{n}_{X}[m]$.\vspace{1mm}

\noindent
{\bf(ii)} Pour toute immersion fermée $\sigma:X\rightarrow Z$ d'espaces complexes, on a  ${\mathcal D}^{\bullet}_{X}\simeq {\mathcal H}om_{{\mathcal O}_{Y}}(\sigma_{*}{\mathcal O}_{X}, {\mathcal D}^{\bullet}_{Y})|_{X}$.\vspace{1mm}

\noindent
{\bf(iii)} Il est à cohomologie cohérente avec
${\mathcal H}^{k}({\mathcal D}^{\bullet}_{X})=0$ pour $k<-{\rm Prof}({\mathcal O}_{X})$.
De plus, si  $X$ est de dimension finie $m$, ${\mathcal D}^{\bullet}_{X}$ est
 d'amplitude $[-m,0]$ et muni d'une application ${\Bbb C}$-linéaire appelée {\emph trace}, ${\rm T}_{X}:{\rm H}^{0}_{c}(X, {\mathcal D}^{\bullet}_{X})\rightarrow {\Bbb C}$
\vspace{2mm}

\noindent
Le $m$-ème faisceau d'homologie de ce complexe est
   exactement le faisceau canonique de Grothendieck\footnote{Le fait qu'il soit dualisant au
  sens de la géométrie analytique complexe relève de résultats
  non triviaux de  Andréotti-Kas  et Golovin (\cite{A.K}, \cite{G}).} $\omega^{m}_{X}={\mathcal H}^{-m}({\mathcal
  D}^{\bullet}_{X})$). De plus, si $X$ est un espace de Cohen Macaulay, le complexe dualisant n'a de cohomologie qu'en degré $-m$ et on a $\omega^{m}_{X}[m]={\mathcal
  D}^{\bullet}_{X})$. \vspace{1mm}

  \noindent
  On
  profite de cette occasion pour rappeler la notion de {\emph{faisceau
  dualisant}}  dans ce cadre.\par\noindent Soit $X$ un espace complexe et ${\mathcal
F}\in{\rm Coh}(X)$. Soit ${\mathfrak U}$ la famille des ouverts de Stein
relativement compacts dans $X$. Alors, à toute inclusion
$U'\subset U$ d'ouverts de ${\mathfrak U}$ est associé, en tout
degré $k$,  un morphisme continu d'espaces vectoriels topologiques
de type {\bf D.F.S} ${\rho^{k}_{U',U}:{\rm H}^{k}_{c}(U', {\mathcal
F})\rightarrow {\rm H}^{k}_{c}(U, {\mathcal F})}$. Pour tout entier $k$,
on désigne par ${\mathcal D}^{k}({\mathcal F})$ le faisceau associé au
préfaisceau \[U\in{\mathfrak U}\rightarrow{\rm {le\,
dual\, topologique\, de}}\,{\rm H}^{k}_{c}(U, {\mathcal F})\]
\[U'\subset U\rightarrow {\rm {le\, transpos\acute{e}\,
du\, morphisme}}\,\rho^{k}_{U',U}\]
\smallskip\noindent Andréotti
et Kas montrent dans \cite{AK} qu'à tout faisceau cohérent ${\mathcal
F}$ sur l'espace complexe $X$, est associé, de façon
fonctorielle, un faisceau analytique ${\mathcal D}^{k}({\mathcal F})$
(appelée  {\emph faisceau  dualisant} de ${\mathcal F}$) qui est  ${\mathcal
O}_{X}$- cohérent  et  vérifie \par $\bullet$   ${\mathcal
D}^{k}({\mathcal F})=0$ si $k\notin[{\rm Prof}({\mathcal F}), {\rm
Dim}({\mathcal F})]$,\par $\bullet$   ${\mathcal D}^{k}({\mathcal F})\simeq
{\mathcal E}xt^{-k}_{{\mathcal O}_{X}}({\mathcal F}, {\mathcal D}^{\bullet}_{X})$
et,\par $\bullet$  pour tout $U\in {\mathfrak U}$, $\Gamma(U, {\mathcal
D}^{k}({\mathcal F}))$ est isomorphe au dual fort de ${\rm H}^{k}_{c}(U,
{\mathcal F})$. \vspace{2mm}

\noindent
Il est intéressant de remarquer que pour ${\mathcal F}={\mathcal O}_{X}$, on obtient ${\mathcal D}^{m}({\mathcal O}_{X})=\omega^{m}_{X}$ le faisceau canonique de Grothendieck et que pour  ${\mathcal F}={\Omega}^{m-k}_{X}$,  ${\mathcal D}^{m}({\Omega}^{m-k}_{X})={\mathcal H}om({\Omega}^{m-k}_{X}, \omega^{m}_{X})$ qui est le faisceau des $k$-formes méromorphes régulières, comme nous allons le voir.  Signalons qu'en considérant le point de vue "homologique"  Golovin (\cite{Go}) généralise ce
résultat au cas où ${\mathcal F}$ est un faisceau
analytique sur $X$ sans hypothèses de cohérence. Il les appellent {\emph{faisceaux d'homologie}}
associé à ${\mathcal F}$ et les note  ${\mathcal H}_{k}({\mathcal
F})$. \vspace{2mm}

\noindent 
\subsubsection{\bf{Formes
méromorphes régulières et propriété de la trace}.}\vspace{1mm}

\noindent 

\noindent 
Dans le cadre algébrique, {\bf E.  Kunz} montre
dans une série d'articles (\cite{Ku1}, \cite{Ku2}, \cite{Ku2})  que toute
variété algébrique projective $X$ irréductible et de
dimension $n$ porte un faisceau dualisant ${\mathcal
K}_{X}=\tilde{\omega}^{n}_{X}$ dont la construction se fait par
recollement relativement aux projections finies. Plus
précisemment, pour tout morphisme fini et dominant $\pi:
X\rightarrow {\Bbb P}_{n}$, il définit un faisceau ${\mathcal O}_{X}$-
cohérent $\omega_{\pi}$ donné par l'isomorphisme
$$ \pi_{*}\omega_{\pi} =  {\mathcal H}om_{{\mathcal O}_{{\Bbb
P}_{n}}}(\pi_{*}{\mathcal O}_{X}, \Omega^{n}_{{\Bbb P}_{n}})$$ dont il
montre l'indépendance vis-à-vis de la projection finie choisie
en utilisant les traces de différentielles. Ce faisceau ${\mathcal
O}_{X}$-cohérent $\tilde{\omega}^{n}_{X}$ est de profondeur au
moins deux dans $X$, coincide avec le faisceau des formes usuelles
aux points lisses de $X$ et est muni d'un morphisme canonique
 $ {\mathcal C}_{X}:\Omega^{n}_{X}\rightarrow \tilde{\omega}^{n}_{X}$
et d'une {\emph {trace}} $\int_{X}: {\rm H}^{n}(X,
\tilde{\omega}^{n}_{X})\rightarrow {\bf k} $.\vspace{1mm}

\noindent Ce
faisceau est appelé le faisceau des {\emph{formes méromorphes
régulières}}.  Kunz montre que le faisceau
$\tilde{\omega}^{n}_{X}$, s'identifie à un sous faisceau du
faisceau des $n$-formes méromorphes entièrement  caractérisé
par la fameuse {\emph {propriété de la trace}} disant que  tout
morphisme fini, surjectif et séparable $f:U\rightarrow U'$ d'un
ouvert $U$ de $X$ sur une variété normale $U'$, induit le
diagramme commutatif
$$\xymatrix{f_{*}f^{*}\Omega^{n}_{U'}\ar[r]^{f^{*}}\ar[d]_{{\mathcal T}^{0}_{f}
\otimes Id}&f_{*}\Omega^{n}_{U}\ar[r]^{f_{*}{\mathcal C}_{U}}&f_{*}\tilde{\omega}^{n}_{U}\ar[d]^{{\mathcal T}^{n}_{f}}\\
 \Omega^{n}_{U'}\ar[rr]_{Id}&& \Omega^{n}_{U'}}$$
Par un calcul explicite, Kunz et Waldi (cf \cite{KW}) étendent cette
construction au cas des degrés intermédiaires en définissant le faisceau
des $r$- formes  méromorphes régulières $\tilde{\omega}^{r}_{X}$ pour lequel ils  
 mettent  en évidence l'isomorphisme canonique $\tilde{\omega}^{r}_{X}\simeq{\mathcal H}om(\Omega^{n-r}_{X}, \tilde{\omega}^{n}_{X})$.\vspace{1mm}

 \noindent
On peut remarquer que pour tout morphisme $\pi: X\rightarrow Y$ fini
et surjectif de variétés compactes de dimension $n$ tel que $X$
soit muni d'un faisceau dualisant ${\mathcal F}_{X}$, le faisceau ${\mathcal
H}om_{{\mathcal O}_{Y}}(\pi_{*}{\mathcal O}_{X}, {\mathcal F}_{X})$ est aussi
dualisant. Signalons, d'une part,  que Kunz généralise sa construction au cas
d'une variété algébrique de Cohen - Macaulay et que, d'autre part,
Lipman (\cite{L}) montre que $\tilde{\omega}^{n}_{X}$ est encore
dualisant si le corps de base est seulement supposé parfait.
\subsection{Formes méromorphes régulières et faisceau ${\omega}^{k}_{Z}$.}\vspace{1mm}

\noindent
Ces formes émanent naturellement d'une théorie de la dualité. A ce titre, le faisceau de Grothendieck est l'exemple type de faisceau de formes régulières en degré maximal.\vspace{1mm}

\noindent
On les retrouve dans \cite{E} sans qu'elles ne soient pas explicitement mises  en
évidence. En effet, utilisant le complexe résiduel ${\rm
K}^{\bullet}_{X}$  ( c'est-à-dire un complexe de ${\mathcal
O}_{X}$-modules injectifs d'amplitude $[-n, 0]$ et dont l'image dans
la catégorie dérivée est le {\emph{ complexe dualisant}}) d'un
schéma $X$ de dimension $n$ et de type fini sur un corps de
caractéristique nulle ${\bf k}$, il montre que  le bicomplexe
${\mathcal K}^{\bullet, \star}_{X}= {\mathcal H}om(\Omega^{\bullet}_{X},
{\rm K}^{\star}_{X})$ peut être muni d'une structure de complexe
différentiel de $(\Omega^{\bullet}_{X}, d)$-modules à droite
doté de deux différentielles $\delta$ et $d^{'}_{X}$, la
première étant naturellement induite par celle de ${\rm
K}^{\bullet}_{X}$, la seconde est  construite pour les besoins de la théorie  au moyen de la différentielle extérieure usuelle du complexe de de Rham (qui
n'est pas ${\mathcal O}_{X}$-linéaire!). Ces objets et ces considérations permettent de voir que  les
faisceaux cohérents $\omega^{\bullet}_{X}= {\mathcal
H}^{\bullet,0}({\mathcal K}^{\bullet, \star}_{X}[-n,-n])$ satisfont,
pour tout $r\leq n$, $\omega^{r}_{X} \simeq {\mathcal H}om_{{\mathcal O}_{X}}(\Omega^{n-r}_{X},
\omega^{n}_{X})$.
 Ces faisceaux que l'on peut appeler faisceaux des {\emph{formes
différentielles régulières}} sont tous de profondeur au moins
deux dans $X$ puisque le ${\mathcal O}_{X}$-module cohérent
$\omega^{n}_{X}:={\mathcal H}^{0}({\rm K}^{\bullet}_{X}[-n])$, qui est
le faisceau de Grothendieck, l'est. Dans le cas réduit de
dimension pure, ils coincident avec les faisceaux de
Kunz.\vspace{2mm}

\noindent
Dans le cadre analytique complexe, le faisceau des formes méromorphes régulières est bien défini sur un espace complexe réduit de dimension pure (disons $m$).  La définition donnée dans \cite{B3} est implicitement donnée dans \cite{Ku1}. D'ailleurs, toutes les propriétés annoncées dans \cite{B3} sont explicitement décrites dans les références déjà citées; la nouveauté de \cite{B3} réside dans son  interprétation en termes de courant $\bar\partial$-fermé.  \vspace{1mm}

\noindent
 Notons $\Sigma$ le lieu singulier de $Z$ et $j $ l'inclusion naturelle de son complémentaire  dans $Z$.   Si ce dernier 
est localement plongé 
avec la codimension $r$ dans une variété 
lisse de Stein $W$,  la classe fondamentale ${\rm C}^{W}_{Z}$ de $Z$ dans $W$, induit, par cup-produit, des morphismes  ${ \Omega^{k}_{Z}\rightarrow {\mathcal E}xt^{r}({\mathcal O}_{Z}, 
\Omega^{k+r}_{W})}$. On pose, alors,   pour tout $k\geq 0$,  d'après \cite{B3}, $ \omega^{k}_{Z} := Ker\{j_{*}j^{*}\Omega^{k}_{Z}\rightarrow 
{\mathcal H}^{1}_{\Sigma}({\mathcal E}xt^{r}({\mathcal O}_{Z}, \Omega^{k+r}_{W}))\}$
 qui n'est que le noyau de  la composée des 
flèches canoniques suivantes $\partial :j_{*}j^{*}\Omega^{k}_{Z}\rightarrow {\mathcal 
H}^{1}_{\Sigma}(\Omega^{k}_{Z})$ et 
 ${\mathcal H}^{1}_{\Sigma}(\Omega^{k}_{Z})\rightarrow 
{\mathcal H}^{1}_{\Sigma}({\mathcal E}xt^{r}({\mathcal O}_{Z}, \Omega^{k+r}_{W}))$
\vspace{1mm}

\noindent
  On dit qu'une section $\sigma$ du faisceau $j_{*}j^{*}\Omega^{k}_{Z}$  vérifie  {\emph {la propriété de la  trace}} si pour toute 
paramétrisation locale $\phi:Z\rightarrow U$
 qui est un revêtement ramifié de degré $r$ de lieu de ramification ${\rm R}(\phi)$ et de branches locales $(f_{l})_{l\in\{1,2,\cdots,r\}}$  et toute forme holomorphe $\alpha$ sur 
$Z$, la forme méromorphe 
$$\sum_{l=1}^{l= r}{f_{l}}^{*}(\alpha\wedge\sigma)$$ 
holomorphe sur $U\setminus{{\rm R}(\phi)\cap U}$, se prolonge analytiquement sur tout $U$.\vspace{1mm}

\noindent Cette propriété caractérise les sections du faisceau ${\omega}^{k}_{Z}$. Cela revient aussi à dire que $\sigma$ est une section de ce faisceau si et seulement si
$\sigma\wedge {\rm C}^{W}_{Z}$ définit génériquement se prolonge en une section  du faisceau 
${\mathcal E}xt^{p}_{{\mathcal O}_{W}}({\mathcal O}_{Z}, \Omega^{k+r}_{W})$.\vspace{1mm}
\subsubsection{\bf{Formes méromorphes régulières  et courants}.}\vspace{1mm}

\noindent \rm On
dispose, grâce à \cite{E1}, d'un résultat explicite comparant  les
résidus de Grothendieck et de Herrera. En effet, il y est montré
que, pour toute paramétristation locale d'un espace analytique
complexe de dimension pure $n$
${\xymatrix{X\ar@/_/[rr]_{f}\ar[r]^{\sigma}&Z:=U\times
B\ar[r]^{q}&U}}$ et toute système $(g_{1},\cdots,g_{p})$ de
fonctions holomorphes sur $Z$ et définissant une suite
régulière sur $X$ en un point générique $x$. Alors, si l'on
désigne par ${\rm Res}^{G}$ (resp. ${\rm Res}^{H}$) le résidu de
Grothendieck (resp. Herrera), on a, pour toute section $w$ de
$\Gamma(Z,\Omega^{n+p}_{Z})$,
$${\mathcal T}_{f}\left[{\rm Res}^{H}_{g_{j}}\left({w\over{g_{1}\cdots g_{p}}}
\right) \right] = (2i\pi)^{n}\sigma_{*}\left[X\right]\wedge{\rm
Res}^{G}\left[{w\over{g_{1}\cdots g_{p}}}\right]$$ ${\mathcal T}_{f}$
étant l'image directe au sens des courants. vspace{2mm}

\noindent L'approche de \cite{B3} est propre à la géométrie analytique complexe puisqu'elle décrit  le faisceau de Kunz  en terme de courants "holomorphes" sur $X$. Plus précisemment, les
considérations locales précédentes permettent d'écrire toute
$r$-forme méromorphe comme un quotient $\xi= { v\over{g}} $ où
$v$ est une $r$-forme holomorphe sur $Z$ et $g$ une fonction
holomorphe s'annulant sur le lieu singulier $\Sigma$ de $X$.
D'après \cite{HL}, $\xi$  définit un courant, 
appelé  {\bf { valeur principale}} et, souvent noté $\xi\wedge
[X]$, définit par
 $$\langle{\bf T}_{\xi},\phi\rangle:= lim_{\epsilon\rightarrow 0}\int_{X\cap\{|g|>\epsilon\}}\xi\wedge \phi$$
pour toute forme $\phi$ de type $(n-r,n)$,  $C^{\infty}$ et à
support compact dans $X$. \vspace{1mm}

\noindent 
Le lien avec les courants se fait, alors,  de la manière suivante:\vspace{1mm}

\noindent
en considérant sur $Z$ localement plongé dans $W$ lisse et de Stein,  une k-forme méromorphe $\xi= { v\over{g}} $ où $v$ est une k-forme holomorphe sur $W$ et $g$ une fonction holomorphe s'annulant sur le lieu singulier $\Sigma$ de $Z$, alors, d'après \cite{HL}, pour toute forme $C^{\infty}$ et à support compact dans $W$, la limite 
$$lim_{\epsilon\rightarrow 0}\int_{Z\cap\{|g|>\epsilon\}}\xi\wedge \phi$$
est le {\bf {courant valeur principale}}  associé à $\xi$ et notée $\xi\wedge [Z]$.  Le faisceau $\omega^{k}_{Z}$ est  défini comme étant le faisceau des courants valeurs principales  $\bar\partial$- fermé, modulo ${\mathcal O}_{X}$- torsion. \vspace{1mm}

\noindent
En le degré maximal, $\omega_{Z}^{m}$ s'identifie au  faisceau de Grothendieck, dualisant au sens 
de Kunz. \vspace{1mm}

\noindent On peut renvoyer le lecteur à \cite{Alek1} pour l'approche en terme de multi-résidus en relation avec le complexe logarithmique.
\subsection{Notations et définitions.} On note \vspace{1mm}

\noindent 
$\bullet$ ${\mathcal M}_{1}$ (resp.  ${\mathcal M}_{0}$) l'ensemble des morphismes d'espaces complexes réduits $\pi:X\rightarrow Y$ arbitraires (resp. à fibres de dimension constantes et surjectifs). \vspace{1mm}

\noindent
$\bullet$ Si $G^{X}_{\bullet}:X\rightarrow{\mathcal F}^{\bullet}_{X}$ est une construction assignant à un espace analytique complexe réduit donné un certain sous faisceau cohérent ${\mathcal F}^{\bullet}_{X}$ de $\omega^{\bullet}_{X}$ coincidant génériquement avec le faisceau des formes holomorphes usuelles $\Omega^{\bullet}_{X}$ et compatible avec les projections c-à-d:
tout morphisme de projection $q:X\times Y\rightarrow Y$ induit un morphisme naturel $q^{*}{\mathcal F}^{\bullet}_{Y}\rightarrow {\mathcal F}^{\bullet}_{X\times Y}$. De plus, si $p:X\times Y\rightarrow X$ est la projection canonique sur $X$ lisse de dimension $n$ et ${\mathcal F}^{\bullet}_{X\times Y}$ coincide génériquement avec le faisceau des $(n+k)$-formes usuelles, $q^{*}{\mathcal F}_{Y}\otimes p^{*}\Omega^{n}_{X}$ est un facteur direct de ${\mathcal F}_{X\times Y}$.

\vspace{1mm}

\noindent
On dira, alors,  que:\vspace{1mm}

\indent $\star$ $G^{X}_{\bullet}$ (ou, par abus, de notation, ${\mathcal F}^{\bullet}_{X}$) est ${\mathcal M}_{1}$ (resp. ${\mathcal M}_{0}$)-stable par image réciproque si tout $f:X\rightarrow Y$ de  ${\mathcal M}_{1}$ (resp. ${\mathcal M}_{0}$) induit un morphisme d'image réciproque  $f^{*}{\mathcal F}^{\bullet}_{S}\rightarrow {\mathcal F}^{\bullet}_{X}$ prolongeant celle des formes holomorphes usuelles, \vspace{1mm}

\indent
$\star$  $G^{X}_{\bullet}$ (ou, par abus, de notation, ${\mathcal F}^{\bullet}_{X}$) est ${\mathcal M}_{1}$ (resp. ${\mathcal M}_{0}$)-stable\footnote{Rappelons que tout morphisme $\pi$ de ${\mathcal M}_{0}$ admet des factorisations locales du type $\xymatrix{X\ar[r]_{f}\ar@/^1pc/[rr]^{\pi}&S\times U\ar[r]_{q}&S}$. } par image directe finie si pour pour tout morphisme fini et surjectif  $f:X\rightarrow Y$ de ${\mathcal M}_{1}$ (resp. ${\mathcal M}_{0}$), il existe un morphisme d'image directe $f_{*}{\mathcal F}^{\bullet}_{X}\rightarrow {\mathcal F}^{\bullet}_{Y}$ rendant commutatif le diagramme
$$\xymatrix{f_{*}{\mathcal F}^{\bullet}_{X}\ar@{^{(}->}[r]\ar[d]&f_{*}\omega^{\bullet}_{X}\ar[d]\\
{\mathcal F}^{\bullet}_{Y}\ar@{^{(}->}[r]&\omega^{\bullet}_{Y}}$$
et induit par l'image directe définie au sens des courants
$f_{*}\omega^{\bullet}\rightarrow \omega^{\bullet}_{Y}$.  
\defn{}{}\label{def1} Si $f:X\rightarrow Y$ est un morphisme fini et surjectif  d'espaces analytiques réduits de dimension pure $m$ et ${\mathcal F}^{k}_{X}$ (resp. ${\mathcal F}^{k}_{Y}$) est un sous faisceau cohérent de $\omega^{k}_{X}$ (resp. $\omega^{k}_{Y}$)  coincidant génériquement avec le faisceau des formes holomorphes usuelles sur $X$ (resp. $Y$), on dit que ${\mathcal F}^{k}_{X}$ possède la {\it propriété de la trace} relativement à $f$ ou est {it stable par trace} si
le morphisme trace générique $f_{*}{\mathcal F}^{k}_{X}\rightarrow {\mathcal F}^{k}_{Y}$ induit par l'image directe définie au sens des courants
$f_{*}\omega^{k}_{X}\rightarrow \omega^{k}_{Y}$ se prolonge globalement et on a le diagramme commutatif
$$\xymatrix{f_{*}{\mathcal F}^{k}_{X}\ar@{^{(}->}[r]\ar[d]&f_{*}\omega^{k}_{X}\ar[d]\\
{\mathcal F}^{k}_{Y}\ar@{^{(}->}[r]&\omega^{k}_{Y}}$$
dans lequel les morphismes horizontaux sont des morphismes trace.\vspace{2mm}

\noindent\rm
Pour un faisceau cohérent ${\mathcal F}$,  son  $\omega^{m}_{X}$-dual (ou $\omega_{X}$-dual) (resp. ${\mathcal L}^{m}_{X}$-dual)  est noté ${\mathcal D}_{\omega_{X}}({\mathcal F}):={\mathcal H}om_{{\mathcal O}_{X}}({\mathcal F}, \omega^{m}_{X})$ (resp. ${\mathcal D}_{{\mathcal L}_{X}}({\mathcal F}):={\mathcal H}om_{{\mathcal O}_{X}}({\mathcal F}, {\mathcal L}^{m}_{X})$).

\section{\color{blue}{Propriétés fondamentales et résultats basiques.}}

\subversionInfo
\subsection{Les propriétés.}
\approvals{Mohamed & yes}

\Prop{}{}\label{prop1} Soit $X$ un espace analytique complexe réduit et de dimension pure $m$. Alors, pour tout entier $k\geq 0$, on a\vspace{1mm}

\noindent 
{\bf(i)} $\omega^{k}_{X}$ est de profondeur au moins $2$;  en particulier, pour  tout sous espace $T$ de codimension 2 dans $X$, on a l'annulation des faisceaux de cohomologie ${\mathcal H}^{0}_{T}(\omega^{k}_{X})={\mathcal H}^{1}_{T}(\omega^{k}_{X})=0$ et $\omega^{k}_{X}\simeq{\mathcal H}om_{{\mathcal O}_{X}}(\Omega^{m-k}_{X}, 
\omega^{m}_{X})$.  Ces faisceaux sont stables par différentiation extérieure et sont $\omega$-réfléxif mais non ${\mathcal O}_{X}$-réfléxif en général\footnote{Une étude un peu plus détaillée des foncteurs de bidualités se trouve dans \cite{K6} . }. Si $X$ est normal, ils sont ${\mathcal O}_{X}$-réfléxifs et munis d'un cup produit interne.   \vspace{1mm}

\noindent
 {\bf(ii)} Tout morphisme $\pi:X\rightarrow S$ morphisme propre surjectif d'espaces analytiques complexes réduits de dimension pure $m$ induit, pour tout entier $k\leq m$,  un morphisme 
 $$\pi_{*}\omega^{k}_{X}\rightarrow \omega^{k}_{S}$$
 injectif  dans le cas d'une modification propre et surjectif dans le cas d 'un morphisme fini et ouvert.  De plus, on a l'isomorphisme $\pi_{*}\omega^{k}_{X}\simeq{\mathcal H}om_{{\mathcal O}_{S}}(\pi_{*}(\Omega^{m-k}_{X}), 
\omega^{m}_{S})$.\vspace{1mm}

\noindent
 {\bf{(iii)}} Soit ${\mathcal F}$ un sous faisceau cohérent de $j_{*}j^{*}\Omega^{k}_{X}$ dont le support rencontre le lieu singulier de $X$ en un fermé d'intérieur vide dans le support de ${\mathcal F}$. Si pour  {\bf une} paramétrisation locale $f: X\rightarrow V$, toute section de ${\mathcal F}$ vérifie la {\it{propriété de la trace}}, alors ${\mathcal F}$ est un sous faisceau de ${\omega}^{k}_{X}$.\vspace{1mm}

 \noindent
 {\bf(iv)} Si $Y$ est un sous espace analytique non entièrement contenu dans ${\rm Sing}(X)$ muni d'un plongement local $i$  dans $X$ et  tel que $Y\cap{\rm Sing}(X))$ soit de codimension au moins deux dans $Y$, alors, on a un morphisme de restriction
${\omega}^{k}_{X}\rightarrow i_{*}{\omega}^{k}_{Y}$ compatible avec la restriction des formes holomorphes usuelles.\vspace{1mm}

 \noindent
{\bf(v)}  Si $\pi:X\rightarrow S$ est un morphisme  ouvert d'espaces complexes réduits tel que le lieu singulier  ${\rm Sing}(S)$ (resp. $\pi^{-1}({\rm Sing}(S))$)  soit de codimension au moins deux dans $S$ (resp. $X$),  on a une image réciproque ${\omega}^{k}_{S}\rightarrow\pi_{*}{\omega}^{k}_{X}$ prolongeant de façon naturelle l'image réciproque des formes holomorphes. En particulier, c'est le cas pour tout morphisme équidimensionnel sur une base normale.\vspace{1mm}

\noindent \rm
\begin{proof} {\bf(i)} Les trois premières propriétés sont classiques (cf \cite{Ku1}, \cite{Ke1}, \cite{B3}). L'isomorphisme $\omega^{k}_{X}\simeq{\mathcal H}om_{{\mathcal O}_{X}}(\Omega^{m-k}_{X}, \omega^{m}_{X})$ montre clairement que $\omega^{k}_{X}$ est de profondeur au moins deux puisqu'il en est ainsi de $\omega^{m}_{X}$. La stabilité par différenciation  extérieure résulte de la construction même de  \cite{B3} ou du diagramme commutatif
$$\xymatrix{0\ar[r]&\omega^{k}_{X}\ar[r]\ar@{.>}
[d]&j_{*}j^{*}\Omega^{k}_{X}\ar[d]_{d}\ar[r]&{\mathcal H}^{1}_{{\rm Sing}(X)}(\Omega^{k}_{X})\ar[d]_{d}\ar[r]&0\\
0\ar[r]&\omega^{k+1}_{X}\ar[r]&j_{*}j^{*}\Omega^{k+1}_{X}\ar[r]&{\mathcal H}^{1}_{{\rm Sing}(X)}(\Omega^{k+1}_{X})\ar[r]&0}$$
La $\omega$-{\emph{réfléxivité}} est une conséquence immédiate de:\vspace{1mm}

\indent $\bullet$ l'accouplement de Yoneda donnant le morphisme\footnote{généralement non  injectif car le produit tensoriel est, par essence, porteur de torsion dans la majeure partie des cas.} 
$$\omega^{k}_{X}\otimes {\mathcal H}om(\omega^{k}_{X}, \omega^{m}_{X})\simeq {\mathcal H}om({\mathcal O}_{X}, \omega^{k}_{X})\otimes {\mathcal H}om(\omega^{k}_{X}, \omega^{m}_{X})\rightarrow  \omega^{m}_{X}$$
duquel résulte le  morphisme canonique 
$$\omega^{k}_{X}\rightarrow {\mathcal H}om({\mathcal H}om(\omega^{k}_{X}, \omega^{m}_{X}), \omega^{m}_{X})$$
injectif puisque génériquement bijectif et $\omega^{k}_{X}$ est sans torsion, \vspace{1mm}

\indent $\bullet$ de l'isomorphisme $\omega^{k}_{X}\simeq{\mathcal H}om(\Omega^{m-k}_{X}, \omega^{m}_{X})$ qui permet, grâce à la flèche naturelle 
$$\Omega^{m-k}_{X}\rightarrow {\mathcal H}om({\mathcal H}om(\Omega^{m-k}_{X}, \omega^{m}_{X}), \omega^{m}_{X})$$
de construire un morphisme 
$${\mathcal H}om({\mathcal H}om(\omega^{k}_{X}, \omega^{m}_{X}), \omega^{m}_{X})\rightarrow \omega^{k}_{X}$$
là aussi nécessairement injectif pour les raisons sus-mentionnées. D'où la conclusion.\vspace{1mm}

\noindent 
L'existence d'un  {\emph{cup produit interne }} comme la propriété  {\bf{(iv)}} ne sont qu'une conséquence de la normalité de $X$ et la profondeur suffisante de ces faisceaux pour autoriser les prolongements en dehors de la codimension $2$. \vspace{1mm}

\noindent
 En effet, ce cup produit est clairement induit par le cup produit
$$\Gamma(Reg(X),\omega^{k}_{X})\otimes\Gamma(Reg(X),\omega^{k'}_{X})\rightarrow\Gamma(Reg(X),\omega^{k+k'}_{X})$$
sachant que, par normalité de $X$, on a, pour tout entier $k$,  $j_{*}j^{*}\omega^{k}_{X}\simeq\omega^{k}_{X}$ et que, par conséquent
$$\omega^{k}_{X}\simeq {\mathcal H}om(\omega^{k'}_{X}, \omega^{k+k'}_{X})$$
puisque deux faisceaux cohérents de profondeur au moins deux coincidant sur la partie régulière de $X$ coincident globalement.
\vspace{1mm}

\noindent
 {\bf(ii)} Si $\pi:X\rightarrow S$ est un morphisme propre surjectif d'espaces analytiques complexes réduits de dimension pure $m$, il est génériquement fini et l'image directe au sens des courants nous donne, pour tout entier $k\leq m$, un morphisme naturel 
 $$\pi_{*}{\mathcal D}^{(k,0)}_{X}\rightarrow {\mathcal D}^{(k,0)}_{S}$$
 qui envoie tout courant $\bar\partial$-fermé et sans torsion sur un courant $\bar\partial$-fermé et sans torsion. D'où le morphisme 
 $$\pi_{*}\omega^{k}_{X}\rightarrow \omega^{k}_{S}$$
 Il est bien connu qu'il est injectif dans le cas d'une modification et surjectif dans le cas fini.\vspace{1mm}

\noindent L'isomorphisme 
$$\pi_{*}{\mathcal H}om({\mathcal F}, \omega^{m}_{X})\simeq {\mathcal H}om(f_{*}{\mathcal F}, \omega^{m}_{Y})$$
peut s'établir de plusieurs façons et l'une d'entre elles consiste à le voir comme une conséquence immédiate du \lemmaref{L2}.  Mais on peut le voir rapidement en appliquant la dualité de Andréotti-Kas \cite{AK}  sur $X$ et $Y$ pour avoir les isomorphismes topologiques
$${\rm I}\!{\rm H}om(X; {\mathcal F}, \omega^{m}_{X})\simeq \big({\rm H}^{m}_{c}(X, {\mathcal F})\big)^{'}$$
$${\rm I}\!{\rm H}om(Y; f_{*}{\mathcal F}, \omega^{m}_{Y})\simeq \big({\rm H}^{m}_{c}(Y, f_{*}{\mathcal F})\big)^{'}$$
la finitude de  $f$ permet de conclure.\vspace{1mm}

\noindent
{\bf(iii)} C'est la caractérisation de $\omega^k_{X}$ ou simplement dûe à l'isomorphisme
$$f_{*}{\mathcal H}om({\mathcal F}, \omega^{m}_{X})\simeq {\mathcal H}om({\mathcal H}om(f_{*}{\mathcal F}, \Omega^{m}_{U})$$
pour toute paramétrisation locale $f:X\rightarrow U$.\vspace{1mm}

\noindent
{\bf(iv)} La {\emph{restriction}} à un sous espace vérifiant les conditions d'incidence citées repose sur le principe du prolongement analytique ou de Hartogs. En effet,  si $X$ est un espace analytique  de dimension $m$ et $Y$ un sous ensemble analytique non entièrement contenu dans le lieu singulier, $\Sigma$,  de $X$. En supposant  que la codimension de $Y\cap \Sigma$ soit supérieure ou égale à 2 dans $Y$), on peut définir une restriction de ces courants en procédant génériquement. En effet,  on a le diagramme 
$$\xymatrix{\Gamma(X,\omega^{k}_{X})\ar[r]\ar[d]&\Gamma(X-\Sigma, \omega^{k}_{X})\simeq \Gamma(X-\Sigma, \Omega^{k}_{X})\ar[d]\\
\Gamma(Y,\omega^{k}_{Y})\eq[r]&\Gamma(Y-\Sigma\cap Y, \Omega^{k}_{Y})}$$ 
on utilise seulement la rrestriction usuelle au niveau des formes holomorphes $(X- X\cap \Sigma)$ à $ (Y-\Sigma\cap Y)$ et le prolongement en dehors de la codimension deux  des sections du faisceau $\omega^{k}_{Y}$.\vspace{1mm}

\noindent 
{\bf(v)} Comme dans {\bf(iv)}, l'existence de cette image réciproque repose esssentiellement sur la possibilité de prolonger les sections en dehors de la codimension deux. On peut remarquer que nos hypothèses imposent l'équidimensionnalité de $\pi$ et que, par conséquent, l'image réciproque d'un fermé d'intérieur vide $\Sigma$ dont le complémentaire est de codimension au moins deux dans $S$ est aussi un fermé d'intérieur vide (ouverture de $\pi$) de codimension au moins deux dans $X$ (équidimensionnalite de $\pi$). Alors, pour tout entier $k\leq dim(S)$, on a un morphisme de ${\mathcal O}_{X}$-modules cohérents $\pi^{*}(\omega^{k}_{S})\rightarrow \omega^{k}_{X}$ 
prolongeant l'image réciproque des formes holomorphes usuelles.\vspace{1mm}

\noindent
On  s'en convainc aisément en commençant par remarquer qu'il nous suffit de traiter le cas où $\Sigma$ est le lieu singulier de $S$. Alors, comme $\pi^{-1}(\Sigma)$ est de codimension au moins deux dans $X$ et que toute  section $\xi$ du faisceau $\omega^{k}_{S}$ sur un certain ouvert de Stein $U$ de $S$ admet une  restriction non triviale sur la partie régulière par absence de torsion, on la relève  en une forme holomorphe  sur l'ouvert $\pi^{-1}(U)\setminus \pi^{-1}(\Sigma\cap U)$. Le morphisme canonique $\Omega^{k}_{X}\rightarrow\omega^{k}_{X}$ permet de voir cette forme comme une section du faisceau $\omega^{k}_{X}$ sur l'ouvert $\pi^{-1}(U)\setminus \pi^{-1}(\Sigma\cap U)$.  Comme le faisceau $\omega^{k}_{X}$ est de profondeur au moins deux sur un espace complexe réduit (en particulier, ils vérifient le principe du prolongement de Hartog),  cette image réciproque se prolonge, alors, en une section globale de $\omega^{k}_{X}$ ce que l'on peut résumer par le diagramme commutatif 
$$\xymatrix{\Gamma(U, \omega^{k}_{S})\ar[r]\ar[d]&\Gamma(U\setminus U\cap \Sigma, \Omega^{k}_{S})\ar[d]\\
\Gamma(\pi^{-1}(U), \omega^{k}_{X})\simeq\Gamma(\pi^{-1}(U)\setminus \pi^{-1}(U\cap \Sigma), \omega^{k}_{X})&\Gamma(\pi^{-1}(U)\setminus \pi^{-1}(U\cap \Sigma), \Omega^{k}_{X})\ar[l]}$$
D'où,  le morphisme image réciproque 
$ \omega^{k}_{S}\rightarrow\pi_{*} \omega^{k}_{X}\,\blacksquare$
\end{proof}

\vspace{1mm}
\section{\color{blue}{Remarques et quelques exemples.}}

\subversionInfo
\subsection{\bf{Comportement par image réciproque et restriction}.} 
\approvals{Mohamed &yes}
On peut attirer l'attention du lecteur sur le fait qu'en général ce type de faisceau se comporte mal par cup produit, par restriction et par image réciproque.  Pour illustrer notre propos considérons l'exemple simple donné par la surface faiblement normale 
$$Z=\lbrace{(x,y,z)\in {\Bbb C}^{3}/ xy^{2} +z^{2}=0}\rbrace$$
communément appelée `` parapluie de Whitney '' dont le lieu singulier est la droite $\lbrace{y=z=0}\rbrace$. Le morphisme de normalisation est donné par l'application de ${\Bbb C}^{2}$ sur $ Z$ envoyant $(u,v)$ sur $(u^{2},v,uv)$. Les formes méromorphes $\displaystyle{{xdy\over{z}}}$ et $\displaystyle{{dz\over{y}}}$ définissent des sections du faisceau $\omega^{1}_{Z}$ sans être des sections du faisceau ${\mathcal L}^{1}_{Z}$ puisqu'elles n'admettent pas de prolongement holomorphe sur la normalisée ${\Bbb C}^{2}$ (ce qui enlève tout espoir d'une image réciproque par morphisme arbitraire!). 
\subsection{\bf{Comportement par cup produit}.} 
\approvals{Mohamed &yes}
  Le cup produit des formes  $\displaystyle{{xdy\over{z}}}$ et $\displaystyle{{dz\over{y}}}$ qui est  $\displaystyle{{xdy\wedge dz\over{zy}}}$ n'est pas une section de $\omega^{2}_{Z}$ car $\displaystyle{xy\over{z}}$ ne se prolonge pas  holomorphiquement  sur  $Z$.  En effet,  comme $\xi:=\displaystyle{dy\wedge dz\over{y^{2}}}$ est un générateur  du faisceau canonique $\omega^{2}_{Z}$, notre forme s'écrit $\displaystyle{xy\over{z}}\xi$. Mais $\displaystyle{{xy\over{z}}=-{z\over{y}}}$ qui est une fonction méromorphe localement bornée (son carré est holomorphe) ne se prolonge pas analytiquement sur $Z$ car si tel était le cas elle serait holomorphe sur tout sous espace de $Z$; ce qui est manifestement faux  sur la courbe $Z\cap\{(x,y,z)\in {\Bbb C}^{3}: x=y\}$ (la restriction est encore  méromorphe localement bornée et non prolongeable holomorphiquement sur cette courbe).\vspace{1mm}

\indent
Même le cup produit  avec des sections de  ${\mathcal L}^{\bullet}_{Z}$ ne stabilise pas $\omega^{\bullet}_{Z}$  alors qu'on aurait pu penser le contraire. \vspace{1mm}
 
 \noindent En effet, en considérant la forme $\displaystyle{{ydx\over{z}}}$ qui définit une section du faisceau ${\mathcal L}^{1}_{Z}$ puisqu'elle se relève holomorphiquement sur la normalisée, on a  $\displaystyle{{dz\over{y}}}\wedge{{ydx\over{z}}}={{dz\wedge dx\over{z}}}$ qui, pour les mêmes raisons que ci-dessus,  n'est pas une section du faisceau $\omega^{2}_{Z}$.\vspace{1mm}

 \noindent
 D'ailleurs, l'existence d'un cup produit $\omega^{k}_{Z}\otimes{\mathcal L}^{m-k}_{Z}\rightarrow\omega^{m}_{Z} $ ( resp. $\omega^{k}_{Z}\otimes{\omega}^{m-k}_{Z}\rightarrow\omega^{m}_{Z} $) correspond exactement à la donnée d'un morphisme
 $$\omega^{k}_{Z}\rightarrow {\mathcal H}om({\mathcal L}^{m-k}_{Z}, \omega^{m}_{Z})$$
 (resp. $\omega^{k}_{Z}\rightarrow {\mathcal H}om({\omega}^{m-k}_{Z}, \omega^{m}_{Z})$)
 nécessairement injectif puisque c'est un isomorphisme générique entre faisceaux sans torsion. Mais on a déjà un morphisme naturel et injectif
 $${\mathcal H}om({\mathcal L}^{m-k}_{Z}, \omega^{m}_{Z})\rightarrow\omega^{k}_{Z}$$ 
 (resp. ${\mathcal H}om({\omega}^{m-k}_{Z}, \omega^{m}_{Z})\rightarrow \omega^{k}_{Z}$) 
 et, donc, l'isomorphisme
  $$\omega^{k}_{Z}\simeq{\mathcal H}om({\mathcal L}^{m-k}_{Z}, \omega^{m}_{Z})$$
  (resp. $\omega^{k}_{Z}\simeq{\mathcal H}om({\omega}^{m-k}_{Z}, \omega^{m}_{Z})$)
  qui n'est réalisé que si $Z$ est normal. En effet, si $\nu:\hat{Z}\rightarrow Z$ est la normalisation, on a
  $$\nu_{*}({\mathcal O}_{\hat Z})\simeq {\mathcal Hom}({\mathcal L}^{m}_{Z}, \omega^{m}_{Z})$$
  et, donc, $$\omega^{0}_{Z}\simeq{\mathcal L}^{0}_{Z}$$
  qui est une condition imposant la normalité de $Z$.\vspace{1mm}
  
  \noindent
  De même, le second isomorphisme 
  $$\omega^{k}_{Z}\simeq{\mathcal H}om({\omega}^{m-k}_{Z}, \omega^{m}_{Z})$$
  donne, en particulier, $$\omega^{0}_{Z}\simeq{\mathcal H}om({\omega}^{m}_{Z}, \omega^{m}_{Z})$$
  Or, pour $Z$ Cohen Macaulay, on a l'isomorphisme
  $${\mathcal O}_{Z}\simeq {\mathcal H}om({\omega}^{m}_{Z}, \omega^{m}_{Z})$$
 qui se déduit de l'auto-dualité
 $${\mathcal O}_{Z}\simeq{\rm I}\!{\rm R}{\mathcal H}om({\mathcal D}^{\bullet}_{Z}, {\mathcal D}^{\bullet}_{Z})$$
 avec pour dualisant  ${\mathcal D}^{\bullet}_{Z}\simeq \omega^{m}_{Z}[m]$.
 D'où l'isomorphisme
 $$\omega^{0}_{Z}\simeq {\mathcal O}_{Z}$$
 traduisant exactement la normalité de $Z$!\vspace{1mm}

 \noindent
  On peut remarquer que dans cet exemple,  la stabilité par cup produit des formes vérifiant des équations de dépendance intégrale est satisfaite. En effet, les formes 
$\xi_{1}:={{zdx\over{x}}}$ et  $\xi_{2}:=2{{zdy\over{y}}}$ vérifient les relations
$$\xi_{1}+\xi_{2}=2dz,\,\,\,\xi_{1}\otimes\xi_{2}=2ydx\otimes dy$$
et, comme on le vérifie facilement, les cup produit de ces formes par chacune des formes ${{xdy\over{z}}}$ ou ${{dz\over{y}}}$ donnent des sections du faisceau  $\omega^{2}_{Z}$. Mais cela semble propre à cette situation particulière  où la normalisée est lisse.\vspace{2mm}

\noindent

\approvals{Mohamed & yes}

\phantomsection\addcontentsline{toc}{part}{\color{blue}{Formes méromorphes et modifications.}}
%
%
\svnid{$Id: S07-proof-setup.tex 269 2020-01-20 11:28:53Z kebekus $}
\section{\color{blue}{le faisceau ${\mathcal L}^{\bullet}_{X}$.}}
Pour $X$ espace analytique complexe réduit et $\pi:\widetilde{X}\rightarrow X$ une désingularisation et $j$ un entier,  on pose
$${\mathcal L}^{j}_{X}:=\pi_{*}\Omega^{j}_{\widetilde{X}}$$
qui est le faisceau des $j$-formes méromorphes sur $X$ dont l'image réciproque se prolonge en une $j$-forme holomorphe sur la désingularisée (cette définition ne dépend pas du choix de la résolution). Il a été introduit en degré maximal dans \cite{GR70}. 
Les principales propriétés de ces faisceaux sont leur stabilité par cup-produit, par différentiation extérieure, par image réciproque arbitraire et par intégration sur les fibres. Une étude détaillée fait l'objet de \cite{K0} et \cite{K4} dans lesquels on montre, entre autres, l'existence de morphismes  $${\mathcal T}^{j}_{\pi,{\mathcal L}}: {\rm I}\!{\rm R}^{n}\pi_{!}{\mathcal L}^{n+j}_{X}\rightarrow{\mathcal L}^{j}_{S}$$ 
( en remplaçant $\pi_{!}$ par $\pi_{*}$ dans le cas propre) 
vérifiant un certain nombre de propriétés dont celle d'être compatible à tous changements de base. Ces morphismes induisent des morphismes d'intégration de classes de cohomologie dans le contexte suivant: 
si $Z$ et $S$  sont deux espaces  analytiques complexes de dimension pure avec $Z$ dénombrable 
à l'infini et $S$ réduit, $(X_{s})_{s\in S}$ une famille analytique de $n$- cycles de $Z$ 
paramétrée par $S$ dont le graphe $X$ est muni du morphisme $\pi: X\rightarrow S$ (induit par la projection canonique sur $S$), $\Phi$ et  $\Psi$ deux familles de supports avec $\Phi$ paracompactifiante, $\Psi$ contenant, pour tout point  $s$ de $S$, les supports $\vert{X_s}\vert$ et $\Phi \cap \Psi$ contenue dans la famille des compacts de $Z$ de sorte que $X\cap (S\times \Phi)$ soit contenu dans la famille des ferm\'es $S$-propre, alors on a des morphismes d'intégration 
$${\rm H}^{n}_{\Phi}(Z, {\mathcal L}^{n+q}_{Z})\longrightarrow 
{\rm H}^{0}(S, {\mathcal L}^{q}_{S})$$
que l'on peut généraliser en considérant des familles de supports $\Phi$, $\Psi$ choisies comme ci-dessus et $\theta$ une famille de supports sur $S$ en supposant que les familles  $X\cap (S\times (\Phi\cap \Psi))$ et $\theta$ soient $\pi$-adaptées c'est-à-dire que $\pi(X\cap (S\times (\Phi\cap \Psi)))$ soit contenue dans $\theta$ et que, pour tout fermé $F$ de $\theta$ la restriction de $\pi$ à $\pi^{-1}(F)\cap X\cap (S\times (\Phi\cap \Psi)) $ soit propre sur $F$.\vspace{1mm}

\noindent
On aura, alors, des morphismes d'intégration
$${\rm H}^{n+p}_{\Phi}(Z, {\mathcal L}^{n+q}_{Z})\longrightarrow 
{\rm H}^{p}_{\Theta}(S, {\mathcal L}^{q}_{S})$$
dotés des mêmes propriétés fonctorielles.\vspace{1mm}

\noindent
Ces faisceaux possèdent d'autres vertus se dévoilant en présence   de singularités particulières(cf \cite{K6}). 
\section{\color{blue}{Les faisceaux $\widehat{\Omega}^{\bullet}_{X}$ et $\overline{\Omega}^{\bullet}_{X}$. }}
\subversionInfo
\subsection{Image réciproque de faisceaux cohérents  modulo torsion.}
\approvals{Mohamed& yes}
Il est bien connu que l'image réciproque analytique d'un faisceau cohérent sans torsion est généralement de torsion (cf par exemple \cite{RS}). Le problème se posait de savoir si pour tout faisceau cohérent ${\mathcal F}$ sur un espace complexe réduit $X$, s'il était possible de trouver un espace complexe réduit $X_{{\mathcal F}}$ et une modification propre $\phi_{\mathcal F}:X_{{\mathcal F}}\rightarrow X$ vérifiant les deux conditions:\vspace{1mm}

\indent
$(\star)$ $\displaystyle{\phi^{*}{\mathcal F}/{\mathcal T}_{\phi}}$
 soit localement libre (${\mathcal T}_{\mathcal F}$étant le sous faisceau de torsion) et\vspace{1mm}
 
 \indent $(\star\star)$ pour tout espace complexe réduit  $Y$ et toute modification propre $\psi:Y\rightarrow X$ vérifiant $(\star)$, il existe un unique morphisme $\theta:Y\rightarrow X_{{\mathcal F}}$ tel que $\psi=\phi_{\mathcal F}\circ \theta$.\vspace{2mm}
 
 \noindent 
 
 L'existence d'une telle donnée a été établie dans  \cite{Ros} et ses propriétés étudiées dans \cite{Riem1}. L'opération consistant à prendre l'image réciproque analytique d'un faisceau sans torsion et à quotienter par la torsion de cette préimage se trouve déjà dans \cite{GR70} et envisagée plus tard dans \cite{RS} ou, dans un contexte un peu différent, dans \cite{You}.  
\vspace{1mm}

\noindent
\subsection{{Désingularisation.}}
\vspace{1mm}

\noindent
 Pour tout espace analytique complexe réduit $X$ de dimension $m$ et toute désingularisation $\phi:\widetilde{X}\rightarrow X$, les faisceaux cohérents  ${\mathcal L}^{j}_{X}:=\phi_{*}\Omega^{j}_{\widetilde{X}}$ dont les sections sont des formes méromorphes vérifiant la propriété de la trace possèdent de bonnes propriétés fonctorielles (cf \cite{K0}, \cite{K4}) parmi lesquelles celles d'être stables par cup produit, par différentiation extérieure, par image réciproque analytique (avec un défaut mineur de compatibilité avec la composition!) et par image directe supérieure (intégration sur les fibres) c'est-à-dire que pour tout morphisme $\pi:X\rightarrow S$ propre surjectif génériquement à fibres de dimension pure $n$, on a un morphisme d'intégration
 $\displaystyle{{\rm I}\!{\rm R}^{n}\pi_{*}{\mathcal L}^{n+j}_{X}\rightarrow {\mathcal L}^{j}_{S}}$ doté d'un certain nombre de propriétés fonctorielles et dont la formation est compatible avec certaines opérations naturelles (pour plus de détails cf \cite{K4}).\vspace{1mm}

 \noindent 
 On désignera dans toute la suite  par $\widetilde{\Omega}^{q}_{X}$ le faisceau des $q$-formes holomorphes modulo torsion. Pour  une désingularisation $\pi:\widetilde{X}\rightarrow X$, on pose
$$\overline{\Omega}^{q}_{X}:=\pi_{*}(\pi^{*}(\widetilde{\Omega}^{q}_{X})/{\mathcal T}_{\pi})$$
qui est clairement un sous faisceau de ${\mathcal L}^{q}_{X}$.
\subsection{{Modification de Nash}}.\vspace{1mm}
 
 \noindent 
 Soit ${\mathcal F}$ est un faisceau cohérent sur $X$, Alors, tout point $x$ de $X$ admet un voisinage ouvert $U$ sur lequel on a la suite exacte
 $${\mathcal O}^{p}_{U}\rightarrow{\mathcal O}^{q}_{U}\rightarrow {\mathcal F}\vert_{U}\rightarrow 0$$ 
 permettant de voir l'algèbre symétrique ${\rm Sym}_{{\mathcal O}_{U}}{\mathcal F}\vert_{U}$ comme un quotient de ${\mathcal O}_{U}[Z_1,\cdots,Z_q]$ par un idéal engendré par un nombre fini d'éléments linéaire en les $Z_j$ et décrivant un sous espace linéaire  de $U\times {\Bbb C}^{q}$. Ces données se recollent (cf\cite{Fi}) pour constituer ce que l'on appelle l'espace linéaire associé à ${\mathcal F}$ noté ${\rm L}({\mathcal F})$ ou ${\rm Specan}_{X}{\rm Sym}_{{\mathcal O}_{X}}({\mathcal F})$. \vspace{1mm}

 \noindent
 Pour un espace complexe $X$ donné, l'espace tangent de Zariski est,  par définition, l'espace linéaire ${\rm L}({\Omega}^{1}_{X})$ noté ${\rm  T}_{X}$,  muni de la projection ${\rm 
 T}_{X}\rightarrow X$  linéaire dans les fibres ${\rm  T}_{X}(x)$ qui s'identifient à l'espace tangent de Zariski usuel $\displaystyle{{\rm Hom}_{{\Bbb C}}(\frak{m}/\frak{m}, {\Bbb C})}$.\vspace{1mm}

 \noindent
Si ${\rm Grass}_{d}(\Omega^{1}_{X})$ est la grassmanienne des sous espaces de dimension $d$ de ${\rm  T}_{X}(x)$, on dispose d'une projection sur $X$ et d'une propriété universelle caractérisant les quotients localement libres disant que pour $f:Y\rightarrow X$ donné, $f^{*}\Omega^{1}_{X}$ admet un quotient localement libre si et seulement si $f$ se factorise au travers de  ${\rm Grass}_{d}(\Omega^{1}_{X})$. 
 L'adhérence, dans ${\rm Grass}_{d}(\Omega^{1}_{X})$,  de l'image de la section naturelle $s:{\rm Reg}(X)\rightarrow {\rm Grass}_{d}(\Omega^{1}_{X})$  est un sous espace analytique fermé réduit ${\rm N}(X)$ muni d'un morphisme propre et biméromorphe $\nu$  sur $X$. Cette modification est  appelée modification de Nash de $X$. On notera  $\widehat{\Omega}^{1}_{X}:=\nu_{*}(\nu^{*}(\Omega^{1}_{X}){\mathcal T})$ l'image réciproque analytique modulo torsion. \vspace{1mm}

\noindent 
 Pour une modification propre $\nu:\overline{X}\rightarrow X$ avec  $\overline{X}$ normal et telle que $\nu^{*}(\widetilde{\Omega}^{q}_{X})/{\mathcal T}^{q}_{\nu}$ soit localement libre, on pose 
$$\widehat{\Omega}^{q}_{X}:={\nu}_{*}({\nu}^{*}(\widetilde{\Omega}^{q}_{X})/{\mathcal T}_{\nu})$$

Dans la terminologie de \cite{B5} où est faite une étude détaillée de ce faisceau $\alpha^{q}_{X}$, ({\bf {Prop 2.1.1}} p.53; {\bf{Def.2.1.2}}, p.54; {\bf{Def.3.0.1}} p.58), $\nu$ est appelée {\emph{normalisation}} du faisceau $\Omega^{q}_{X}$. Elle est minimale dans le sens où tout morphisme $\phi:Y\rightarrow X$ avec $Y$ normal et ${\phi}^{*}(\widetilde{\Omega}^{q}_{X})/{\mathcal T}_{\phi}$ est localement libre se factorise au travers de $\nu$. Toute modification propre de $X$ est dominée par une désingularisation de $\overline{X}$. Pour l' indépendance de la construction  vis à vis de la modification choisie, on renvoie le lecteur  au {\bf{Cor.3.0.3}} p.59. La propriété de  la clôture intégrale est prouvée dans la  {\bf{prop 2.2.4}} p.57; {\bf{prop 5.2.1}} p.72) et sa stabilité  par image réciproque avec compatibilité avec la composition dans  (cf \cite{B6}, {\bf{thm 1.0.1}}).\vspace{1mm}

 \noindent
\subsection{{Propriétés de ces faisceaux.}}\vspace{1mm}

\noindent
Dans \cite{Sp}, on trouve la description suivante: \vspace{1mm}
 
 \noindent si $X$ est un  espace analytique complexe normal et  ${\mathcal F}$, un faisceau cohérent sans torsion  de fibré linéaire ${\mathcal T}$ de rang $r$ localement trivial en dehors d'un sous espace fermé $\Sigma$. Si ${\mathcal M}_{X}$ désigne le faisceau des fonctions méromorphes sur $X$, on pose ${\mathcal M}_{X}({\mathcal F}):={\mathcal M}_{X}\otimes_{{\mathcal O}_{X}} {\mathcal F}$ et\vspace{1mm}
 
 \indent $\bullet$ $\overline{\mathcal F}:=\{\xi\in {\mathcal M}_{X}({\mathcal F}): \xi{\,\rm{soit\,entier}\,sur}\,{\rm S}^{1}({\mathcal F})\} $
 \vspace{1mm}
 
 \indent $\bullet$ $\widetilde{\mathcal F}$ le faisceau associé au préfaisceau
 $$\displaystyle{\widetilde{\mathcal F}(U):=\{\xi\in {\mathcal F}(U\setminus \Sigma)/\, \xi:{\rm T}^{1}(X)\vert_{U\setminus \Sigma}\rightarrow {\Bbb C}\,\,{{\rm localement}}\,{\rm {born\acute ee}}\}}$$ \vspace{1mm}
 
 \indent
 $\bullet$ si $\nu:\overline{X}\rightarrow X$ est l'éclatement de Nash normalisé associé à ${\mathcal T}$, ${\mathcal G}$ le fibré localement trivial prolongeant $\nu^{*}({\mathcal T})\vert_{\overline{X}\setminus \Sigma}$ à $\overline X$ et ${\mathcal G}^{\nu}$ le faisceau dual. Alors, on a (cf \cite{Sp})
 \Prop{}{}\label{prop2'} \vspace{1mm}

 \noindent 
{\bf(i)} $$\overline{\mathcal F}\simeq \widetilde{\mathcal F}\simeq\nu_{*}({\mathcal G}^{\nu})$$
{\bf(ii)} pour $\displaystyle{{\mathcal F}:=\Omega^{1}_{X}}$, toute forme méromorphe sur $X$ qui, considérée comme fonction méromorphe sur ${\rm T}^{1}(X)$, se prolonge holomorphiquement sur globalement, définit une forme holomorphe sur $X$.\vspace{1mm}

\noindent
{\bf(iii)} si $X$ et ${\rm T}^{1}(X)$ sont normaux, $\Omega^{1}_{X}=\overline{\Omega}^{1}_{X}$. En particulier  vrai, pour une hypersurface normale de ${\Bbb C}^{m+1}$ dont le lieu singulier est de codimension supérieure ou égale à $3$ (et donc pour toute hypersurface  à singularité isolée de dimension plus grande que $3$).\rm

 \Prop{}{}\label{prop2} Soit $X$ un espace complexe réduit et $q$ un entier naturel. Alors, le faisceau $\widehat{\Omega}^{q}_{X}$ possède les propriétés suivantes:\vspace{1mm}

\noindent
 {\bf(i)} Sa construction est indépendante de la modification normalisante choisie ,\vspace{1mm}

 \noindent
 {\bf(ii)} il est sans ${\mathcal O}_{X}$-torsion, coincide avec $\widetilde{\Omega}^{q}_{X}$ sur la partie régulière de $X$ avec 
 $$\widetilde{\Omega}^{q}_{X}\subset\widehat{\Omega}^{q}_{X}\subset {\mathcal L}^{q}_{X}\subset \omega^{q}_{X} $$
 {\bf(iii)} il est stable par image réciproque dans le sens suivant:
 si $f:Z\rightarrow X$ est un morphisme d'espaces complexes réduits tel que $\displaystyle{f^{-1}({X})\not\subset {\rm Sing}(Z)}$, on a un morphisme ${\mathcal O}_{Z}$-linéaire $\displaystyle{f^{*}\widehat{\Omega}^{q}_{X}\rightarrow \widehat{\Omega}^{q}_{Z}}$ prolongeant naturellement celle des formes holomorphes et compatible avec la composition sous certaines conditions d'incidences.\vspace{1mm}

 \noindent
 {\bf(iv)} Il s'identifie à la clôture intégrale dans $\omega^{q}_{X}$
 du faisceau $\widetilde{\Omega}^{q}_{X}$.\vspace{1mm}
 
 \noindent {\bf v)} Ils sont stable par cup produit et s'écrivent localement comme somme finie de produit de fonction localement bornée génériquement continues par des formes holomorphes (cf \cite{B5}, {\bf{Thm 3.0.2}, p.58}).\rm 
 \vspace{1mm}

  \noindent
\Prop{}{}\label{prop3} Soit $X$ un espace complexe réduit et $q$ un entier naturel. Alors, le faisceau $\overline{\Omega}^{q}_{X}$ possède les propriétés {\bf(i)}, {\bf(ii)} et {\bf(iii)} énumérées dans la proposition précédente.\rm\vspace{1mm}

\noindent
\begin{proof} {\bf(i)} Le faisceau $\overline{\Omega}^{q}_{X}$ est  bien défini et ne dépend pas de la résolution choisie (comme pour le faisceau $\pi_{*}\Omega^{q}_{\widetilde X}$ dont l'étude est menée dans \cite{K4}). Cela  reposent essentiellement sur le fait que deux résolutions peuvent toujours être dominées par une troisième et qu'un morphisme de faisceaux cohérents envoie toujours la torsion de l'un sur celle de l'autre et induira naturellement un morphisme au niveau des quotients modulo torsion de ces  faisceaux. Plus précisément, étant données deux résolutions $\pi_{j}:X_{j}\rightarrow X$, $j=1,2$, il existe une résolution $\pi_{3}:X_{3}\rightarrow X$ et un diagramme commutatif
$$\xymatrix{&X_{3}\ar[dd]^{\psi}\ar[ld]_{\phi_1}\ar[rd]^{\phi_2}&\\
X_{1}\ar[rd]_{\pi_{1}}&&X_{2}\ar[ld]^{\pi_{2}}\\
&X&}$$
Il est très facile de voir que (cf \cite{B4}, {\bf{Lem.1.0.4}}, p.52)
$${\phi_{1}}^{*}({\pi_1}^{*}(\widetilde{\Omega}^{q}_{X})/{\mathcal T}_{\pi_{1}})/{\mathcal T}_{\phi_{1}} \simeq {\psi}^{*}(\widetilde{\Omega}^{q}_{X})/{\mathcal T}_{\psi}
\simeq {\phi_{2}}^{*}({\pi_2}^{*}(\widetilde{\Omega}^{q}_{X})/{\mathcal T}_{\pi_{2}})/{\mathcal T}_{\phi_{2}}$$
Comme $X_j$ est lisse (donc normal) et que ${\pi_j}^{*}(\widetilde{\Omega}^{q}_{X})/{\mathcal T}_{\pi_{j}}$ est sans torsion (donc de lieu singulier de  codimension au moins deux dans $X_j$), on obtient, en appliquant $\psi_{*}$ (\cite{RS}, {\bf{Th 4.1}} p.253), les isomorphismes
$${\pi_{1}}_{*}({\pi_1}^{*}(\widetilde{\Omega}^{q}_{X})/{\mathcal T}_{\pi_{1}})\simeq {\psi}^{*}(\widetilde{\Omega}^{q}_{X})/{\mathcal T}_{\psi}
\simeq {\pi_{2}}_{*}({\pi_2}^{*}(\widetilde{\Omega}^{q}_{X})/{\mathcal T}_{\pi_{2}})$$
et donc l'indépendance vis-à-vis de la résolution choisie.\vspace{1mm}

\noindent 
On a utilisé le fait que, pour une modification propre $\phi:Z\rightarrow Y$ avec  $Z$ lisse et $Y$ normal, les morphismes canoniques et injectifs
 $${\mathcal F}\rightarrow\phi_{*}\phi^{*}{\mathcal F},\,\, {\mathcal F}\rightarrow\phi_{*}({\phi}^{*}({\mathcal F})/{\mathcal T}_{\phi})\,\,\,$$
 sont des isomorphismes pour tout faisceau cohérent ${\mathcal F}$ sans torsion sur $Y$ dont on se convainc rapidement en remarquant que, pour le premier, il suffit de  voir qu'il est réalisé sur les faisceaux localement libres puisque $\phi_{*}{\mathcal O}_{Z}={\mathcal O}_{Y}$ par normalité de $Y$ et de passer au cas général en utilisant un début de résolution  localement libre à deux termes pour ${\mathcal F}$. Pour le second, il suffit de se rappeler que ${\mathcal F}$ étant sans torsion, peut être vu comme un sous faisceau d'un faisceau localement libre; ce qui permet de l'installer dans une suite exacte courte du type
 $$\xymatrix{0\ar[r]&{\mathcal F}\ar[r]&{\mathcal L}\ar[r]&{\mathcal T}\ar[r]&0}$$
 où ${\mathcal T}$ est un faisceau de torsion supporté par le lieu singulier de ${\mathcal L}$ qui est de codimension au moins deux puisque ${\mathcal F}$ est sans torsion et $Y$ normal.\vspace{1mm}
 \noindent
 On obtient la suite exacte 
 $$\xymatrix{0\ar[r]&{\mathcal T}or^{{\mathcal O}_{Y}}_{1}({\mathcal O}_{Z}, {\mathcal T})\ar[r]&\phi^{*}{\mathcal F}\ar[r]&\phi^{*}{\mathcal L}\ar[r]&\phi^{*}{\mathcal T}\ar[r]&0}$$
 Il est clair que ${\mathcal T}or^{{\mathcal O}_{Y}}_{1}({\mathcal O}_{Z}, {\mathcal T})$ s'identifie, par construction, au sous faisceau de torsion ${\mathcal T}_{\phi}$. Par conséquent, on a une suite exacte courte
 $$\xymatrix{0\ar[r]&\phi^{*}{\mathcal F}/{\mathcal T}_{\phi}\ar[r]&\phi^{*}{\mathcal L}\ar[r]&\phi^{*}{\mathcal T}\ar[r]&0}$$
qui donne, le diagramme commutatif à lignes exactes
$$\xymatrix{0\ar[r]&{\mathcal F}\ar@{^{(}->}[d]\ar[r]&\phi^{*}{\mathcal L}\eq[d]\ar[r]&\phi^{*}{\mathcal T}\ar[r]\eq[d]&0&\\
0\ar[r]&\phi_{*}((\phi^{*}{\mathcal F})/{\mathcal T}_{\phi})\ar[r]&\phi_{*}\phi^{*}{\mathcal L}\ar[r]&\phi_{*}\phi^{*}{\mathcal T}\ar[r]&{\rm I}\!{\rm R}^{1}\phi_{*}((\phi^{*}{\mathcal F})/{\mathcal T}_{\phi})\ar[r]&0}$$
 et donc l'assertion voulue avec, de plus, ${\rm I}\!{\rm R}^{1}\phi_{*}((\phi^{*}{\mathcal F})/{\mathcal T}_{\phi})=0$.\vspace{1mm}

 \noindent
{\bf(ii)} Ce faisceau est un sous faisceau  du faisceau ${\mathcal L}^{q}_{X}$ puisque, pour toute désingularisation  $\pi:\widetilde{X}\rightarrow X$,  $\pi^{*}(\widetilde{\Omega}^{q}_{X})/{\mathcal T}_{\pi}$ est un sous faisceau de $\Omega^{q}_{\widetilde{X}}$. Il est donc  sans torsion et coincident génériquement avec le faisceau des formes holomorphes usuelles. On a les inclusions 
 $$\widetilde{\Omega}^{q}_{X}\subset\overline{\Omega}^{q}_{X}\subset {\mathcal L}^{q}_{X}\subset\omega^{q}_{X}$$
   {\bf(iii)} La stabilité par image réciproque et sa compatibilité à la composition des morphismes sous les conditions d'incidences imposées dans \cite{B6} se prouve exactement comme pour $\widehat{\Omega}^{q}_{X}$\end{proof}

\noindent
Alors, 
\Prop{}{}\label{prop4}  On a un isomorphisme naturel $$\overline{\Omega}^{q}_{X}\simeq\widehat{\Omega}^{q}_{X}$$
\rm
\begin{proof} \vspace{1mm}

\noindent Si $\pi:\widetilde{X}\rightarrow X$ est une résolution donnée et $\nu$ une modification associée au faisceau $\widetilde{\Omega}^{q}_{X}$, on peut construire (en considérant, par exemple, la résolution du réduit du produit fibré $\widetilde{X}\times_{X} \overline{X}$) une résolution de $\overline{X}$ (cf \cite{B4}) $\phi:Z\rightarrow X$ se factorisant au travers de $\pi$ en le diagramme  commutatif
$$\xymatrix{&{Z}\ar[ld]_{\theta}\ar[rd]^{\phi}&\\
\widetilde{X}\ar[rr]_{\pi}&&X}$$
avec 
$$\theta^{*}(\pi^{*}(\widetilde{\Omega}^{q}_{X})/{\mathcal T}_{\pi}))/{\mathcal T}_{\theta}\simeq \phi^{*}(\widetilde{\Omega}^{q}_{X})/{\mathcal T}_{\phi}$$
et  $\phi^{*}(\widetilde{\Omega}^{q}_{X})/{\mathcal T}_{\phi}$ localement libre.\vspace{1mm}

\noindent
En appliquant le foncteur $\phi_{*}$ et sachant que $\theta_{*}\circ\theta^{*}={\rm Id}$ puisque $\widetilde{X}$ est normal car lisse, on obtient l'isomorphisme

$$\overline{\Omega}^{q}_{X}\simeq\widehat{\Omega}^{q}_{X}\,\,\,\,\blacksquare$$
On peut relever que $\overline{\Omega}^{0}_{X}={\mathcal L}^{0}_{X}$.
\end{proof}

\noindent 

\noindent

\begin{rem}  {\bf(i)} {\it{Les faisceaux $\overline{\Omega}^{\bullet}_{X}$ ne sont pas stables par différentiation extérieure}}:\rm\vspace{1mm}

\indent
$\bullet$ Sur $\displaystyle{X:=\{(x,y)\in {\Bbb C}^{2}: x^{2} -y^{3}=0\}}$ dont le tangeant de Zariski est donné par
$\displaystyle{{\rm T}^{1}(X):=\{(x,y,\alpha, \beta)\in X\times {\Bbb C}^{2}: x^{2} -y^{3}=0, 2x\alpha -3y^{2}\beta=0\}}$, la fonction méromorphe localement bornée $f(x,y)=\frac{x}{y}$ vérifiant une équation de dépendance intégrale puisque $f^{2}=y$ et définit, donc, une section du faisceau $\overline{\Omega}^{0}_{X}$ et a pour différentielle extérieure $\displaystyle{df= \frac{dx}{y}-\frac{xdy}{y^2}=\frac{dx}{3y}}$ qui, sur le tangeant de Zariski, s'écrit:
$\displaystyle{\xi:=\frac{\alpha}{3y}=\frac{3y\beta}{2x}}$ dont on vérifie facilement qu'elle n'est pas localement bornée puisque, par exemple sa restriction  sur la courbe décrite par $x=t^6$, $y=t^4$, $\alpha=3t^3$ et $\beta=2t$  donne un terme en $\frac{1}{t}$.\vspace{1mm}

\indent $\bullet$  Sur la surface faiblement normale $\displaystyle{X:=\{(x,y,z)\in {\Bbb C}^{3}: yx^{2}=z^{2}\}}$ de normalisation $(u,v)\rightarrow (u, v^2, uv)$, la fonction méromorphe localement bornée $\displaystyle{\frac{z}{x}}$ a pour  différentielle $\displaystyle{\frac{x}{2z}dy}$ qui est une section de ${\mathcal L}^{1}_{X}$ mais pas de $\overline{\Omega}^{1}_{X}$ dont on s'en convainc aisément en nous ramenant  au cas précédent en prenant $x=y$ ou en regardant sa restriction à la courbe $x=t$, $y=t^4$ et $z=t^3$.\vspace{1mm}

\noindent 
Remarquons, au passage, que la relation différentielle
$$ \frac{2yx}{z}dx + \frac{x^{2}}{z}dy= 2dz$$
qui s'écrit aussi 
$$ \frac{2z}{x}dx + \frac{z}{y}dy= 2dz$$
montre que ces deux formes méromorphes dont le produit tensoriel $\displaystyle{\frac{2z}{x}dx \otimes\frac{z}{y}dy=2xdx\otimes dy}$ est holomorphe définissent des sections du faisceau $\overline{\Omega}^{1}_{X}$. On relève aussi  que si le coefficient $\displaystyle{\frac{z}{x}}$ est localement borné, il n'en est pas de même de $\displaystyle{\frac{z}{y}}$ (considérer la restriction  à la courbe $x=t$, $y=t^4$ et $z=t^3$ qui donne un terme $\displaystyle{\frac{1}{t}}$!). Cependant la forme méromorphe $\displaystyle{\frac{z}{y}dy}$ s'écrit comme somme de produit de formes holomorphes par des fonctions génériquement continus et localement bornés!
\vspace{1mm}

\indent
$\bullet$ Pour $\displaystyle{X:=\{(x,y,z)\in {\Bbb C}^{3}: yx^{3}=z^{6}\}}$ de normalisation $(u,v)\rightarrow (u^2, v^6, uv)$, la relation différentielle
$$ \frac{3yx^{2}}{z^5}dx + \frac{x^{3}}{z^5}dy= 6dz$$
permet de voir que les formes $\displaystyle{\frac{yx^{2}}{z^k}dx}$ et $\displaystyle{\frac{x^{3}}{z^k}dy}$
sont des sections du faisceau  ${\mathcal L}^{1}_{X}$ pour tout entier $k\leq 5$, qu'elles sont sections du faisceau $\overline{\Omega}^{1}_{X}$ pour $k\leq 4$. Pour $k\leq 3$, elles vérifient la même équation de dépendance intégrale de degré $2$ alors que pour $k=4$, elles vérifient deux équations de dépendance intégrale de degré $3$ différentes. En effet,  la puissance tensorielle troisième de $\displaystyle{\frac{yx^{2}}{z^4}dx=\frac{z^{2}}{x}dx}$ $\displaystyle{\xi^{\otimes 3}_{1}}$ est $\displaystyle{ydx^{\otimes 3}}$,  d'où la première équation. Il est, alors, facile d'expliciter l'équation tensorielle de degré $3$ vérifiée par $\displaystyle{\frac{x^{3}}{z^4}dy=\frac{z^{2}}{y}dy}$.  \vspace{1mm}

\noindent Notons que parmi les formes $\displaystyle{\frac{ydx}{z^k}}$ et   $\displaystyle{\frac{xdy}{z^k}}$ qui sont des sections de $\omega^{1}_{X}$, seules  $\displaystyle{\frac{ydx}{z}}$,  $\displaystyle{\frac{xdy}{z}}$ et $\displaystyle{\frac{xdy}{z^2}}$ définissent des sections  du faisceau ${\mathcal L}^{1}_{X}$. \vspace{1mm}

\noindent 
Enfin, là encore  le coefficient $\displaystyle{g:={\frac{z^{2}}{x}}}$ est une fonction méromorphe localement bornée (de cube holomorphe) et définit, donc, une section du faisceau  $\overline{\Omega}^{0}_{X}:=\overline{\mathcal O}_{X}$ alors que ce n'est pas le cas pour $\displaystyle{\frac{z^{2}}{y}}$ (sa restriction à la courbe paramétrée par $(x=t, y=t^3, z=t)$ donne $\displaystyle{\frac{1}{t}}$). Cependant la forme méromorphe  $\displaystyle{\frac{z^{2}dy}{y}}$ s'écrit comme somme de produit de formes holomorphes par des fonctions génériquement continues et localement bornées.  \vspace{1mm}

\noindent
Pour les surfaces normales du type  $\displaystyle{X:=\{(x,y,z)\in {\Bbb C}^{3}: yx=z^{k}\}}$ avec $k\geq 3$ ou l'hypersurface normale de ${\Bbb C}^{4}$ donnée par $\displaystyle{X:=\{(x,y,z,t)\in {\Bbb C}^{3}: yx=zt\}}$, on renvoie le lecteur à \cite{B5} et précisément pour la dernière hypersurface au {\bf{lemma 6.4.2} p.80}. Remarquons que ces exemples présentent tous des singularités rationnelles et dans ce cas les faisceaux $\omega^{\bullet}_{X}$ et ${\mathcal L}^{\bullet}_{X}$ coincident.\vspace{1mm}

\noindent 
{\bf(ii)} {\it{Les faisceaux $\overline{\Omega}^{\bullet}_{X}$ ne sont pas stables par image directe de modifications propres}}\rm:\vspace{1mm}

\noindent 
Plus précisémment si $\phi:Y\rightarrow X$ est une modification propre, on a pas de morphisme naturel d'image directe $\phi_{*}\overline{\Omega}^{\bullet}_{Y}\rightarrow \overline{\Omega}^{\bullet}_{X}$ alors que de telles flèches existent pour les faisceaux $\omega^{\bullet}_{Y}$ et ${\mathcal L}^{\bullet}_{Y}$ (avec bijectivité pour ce dernier!). En effet, il est très facile de voir que 
l'image réciproque définit un morphisme injectif 
$\displaystyle{\overline{\Omega}^{\bullet}_{X}\rightarrow \phi_{*}\overline{\Omega}^{\bullet}_{Y}}$ (cf \cite{B5}, {\bf{Lem 5.2.3}, p.72}) génériquement bijective entre faisceaux sans torsion. Mais, alors, cela  impose l'isomorphisme $\overline{\Omega}^{\bullet}_{X}\simeq\phi_{*}\overline{\Omega}^{\bullet}_{Y}$; ce qui est généralement faux car cela entrainerait, en particulier,  que cela soit vrai pour toute résolution et donc que $\overline{\Omega}^{\bullet}_{X}\simeq {\mathcal L}^{\bullet}_{X}$ !
\end{rem} 
\vspace{1mm}

\noindent
\section{\color{blue}{Entre formes holomorphes et formes méromorphes régulières.}}

Les inclusions $$\widetilde{\Omega}^{q}_{X}\subset \overline{\Omega}^{q}_{X}=\widehat{\Omega}^{q}_{X}\subset {\mathcal L}^{q}_{X}\subset \omega^{q}_{X}$$ nous donnent les inclusions
$$\xymatrix{&{\mathcal H}om_{{\mathcal O}_{X}}({\omega}^{m-j}_{X},
{\mathcal L}^{m}_{X})\ar@{^{(}->}[r]\ar@{^{(}->}[d]&{\mathcal H}om_{{\mathcal O}_{Z}}({\omega}^{m-j}_{X},{\omega}^{m}_{X})\ar@{^{(}->}[d]\\
{\mathcal L}^{j}_{X}\ar@{^{(}->}[r]&{\mathcal H}om_{{\mathcal O}_{X}}({\mathcal L}^{m-j}_{X},
{\mathcal L}^{m}_{X})\ar@{^{(}->}[r]\ar@{^{(}->}[d]&{\mathcal H}om_{{\mathcal O}_{X}}({\mathcal L}^{m-j}_{X}{\omega}^{m}_{X})\ar@{^{(}->}[d]\\
&{\mathcal H}om_{{\mathcal O}_{X}}(\overline{\Omega}^{m-j}_{X}, {\mathcal L}^{m}_{X})\ar@{^{(}->}[r]\ar@{^{(}->}[d]&{\mathcal H}om_{{\mathcal O}_{X}}(\overline{\Omega}^{m-j}_{X},\omega^{m}_{X})\ar@{^{(}->}[d]\\
&{\mathcal H}om_{{\mathcal O}_{Z}}(\widetilde{\Omega}^{m-j}_{X}, {\mathcal L}^{m}_{X})\ar@{^{(}->}[r]&{\mathcal H}om_{{\mathcal O}_{X}}(\widetilde{\Omega}^{m-j}_{X},\omega^{m}_{X}):=\omega^{j}_{X}}$$

\noindent puisque  tous ces faisceaux coincident génériquement et sont sans torsion. Tous ceux  de la seconde colonne sont de profondeurs au moins deux dans $X$ et coincident donc en dehors de la codimension deux (en particulier, si $X$ est normal, ils coincident tous!). Ces inclusions sont généralement strictes comme on peut le voir sur certains exemples du {\bf \S 5}. \vspace{1mm}

\noindent 

Si  $Y$  est un espace complexe réduit de dimension $m$ et ${\mathcal F}$ un faisceau cohérent sur $Y$,  rappelons  que le  $\omega^{m}_{Y}$ (resp. ${\mathcal L}^{m}_{Y}$)-dual de ${\mathcal F}$ est noté ${\mathcal D}_{\omega_{Y}}({\mathcal F}):={\mathcal H}om_{{\mathcal O}_{Y}}({\mathcal F}, \omega^{m}_{Y})$ (resp. ${\mathcal D}_{{\mathcal L}_{Y}}({\mathcal F}):={\mathcal H}om_{{\mathcal O}_{Y}}({\mathcal F}, \omega^{m}_{Y})$ ). 

\Prop{}{}\label{prop5} Soit $X$ un espace analytique réduit.  On a:\vspace{2mm}

\noindent 
{\bf(i)} $${\mathcal D}_{{\mathcal L}_{X}}(\widetilde{\Omega}^{m-j}_{X})\simeq {\mathcal D}_{{\mathcal L}_{X}}(\overline{\Omega}^{m-j}_{X})$$
{\bf(ii)} si $\phi:X'\rightarrow X$ est le morphisme de normalisation,
$$\phi_{*}{\mathcal D}_{{\mathcal L}_{X'}}({\mathcal L}^{m-j}_{X'})\simeq {\mathcal D}_{{\mathcal L}_{X}}({\mathcal L}^{m-j}_{X})$$
$$\phi_{*}{\mathcal D}_{{\omega}_{X'}}({\mathcal L}^{m-j}_{X'})\simeq {\mathcal D}_{{\omega}_{X}}({\mathcal L}^{m-j}_{X})$$
En particulier, si $X'$ est lisse, $$\phi_{*}{\mathcal D}_{{\mathcal L}_{X'}}({\mathcal L}^{m-j}_{X'})\simeq {\mathcal D}_{{\mathcal L}_{X}}({\mathcal L}^{m-j}_{X})\simeq{\mathcal D}_{{\omega}_{X}}({\mathcal L}^{m-j}_{X})$$ 
et un morphisme injectif (génériquement bijectif)
$$\phi_{*}{\mathcal D}_{{\mathcal L}_{X'}}(\overline{\Omega} ^{m-j}_{X'})\rightarrow {\mathcal D}_{{\mathcal L}_{X}}(\overline{\Omega} ^{m-j}_{X})$$
\noindent{\bf(iii)} 
 les faisceaux  ${\mathcal D}_{{\omega}_{X}}({\mathcal L}^{m-j}_{X})$, ${\mathcal D}_{{\omega}_{X}}(\phi_{*}\overline{\Omega}^{m-j}_{X'})$ et $ {\mathcal D}_{{\omega}_{X}}(\phi_{*}\widetilde{\Omega}^{m-j}_{X'})$ sont isomorphes et s'injectent dans le faisceau ${\mathcal D}_{{\omega}_{X}}(\overline{\Omega}^{m-j}_{X})$.\vspace{1mm}

\noindent 
{\bf(iv)} En posant $\overline{\mathcal L}^{j}_{X}:={\mathcal D}_{{\mathcal L}_{X}}({\mathcal L}^{m-j}_{X})$, $\widehat{\mathcal L}^{j}_{X}:={\mathcal D}_{{\mathcal L}_{X}}(\widetilde{\Omega}^{m-j}_{X})$ et $\widetilde{\mathcal L}^{j}_{X}:={\mathcal D}_{{\omega}_{X}}({\mathcal L}^{m-j}_{X})$, on a les cup-produits 
$$\widehat{\mathcal L}^{j}_{X}\otimes \widetilde{\Omega}^{k}_{X}\rightarrow \widehat{\mathcal L}^{j+k}_{X},\,\,\overline{\mathcal L}^{j}_{X}\otimes{\mathcal L}^{k}_{X}\rightarrow\overline{\mathcal L}^{j+k}_{X},\,\,\widetilde{\mathcal L}^{j}_{X}\otimes{\mathcal L}^{k}_{X}\rightarrow\widetilde{\mathcal L}^{j+k}_{X} $$
\rm
\begin{proof} {\bf(i)} L'isomorphisme $${\mathcal H}om_{{\mathcal O}_{X}}(\overline{\Omega}^{m-j}_{X}, {\mathcal L}^{m}_{X})\simeq {\mathcal H}om_{{\mathcal O}_{X}}(\widetilde{\Omega}^{m-j}_{X}, {\mathcal L}^{m}_{X})$$
s'explique très simplement eu égard à la définition et à la formule d'adjonction. En effet, si $\pi:\widetilde{X}\rightarrow X$ est une désingularisation, on a, puisque ${\mathcal T}_{\pi}$ est de torsion,
$${\mathcal H}om(\pi^{*}({\widetilde{\Omega}}^{m-j}_{X})/{\mathcal T}_{\pi}, {\Omega}^{m}_{\widetilde{X}})\simeq {\mathcal H}om(\pi^{*}({\widetilde{\Omega}}^{m-j}_{X}), {\Omega}^{m}_{\widetilde{X}})$$
qui donne, par adjonction, l'isomorphisme voulue.\vspace{1mm}

\noindent 
{\bf(ii)} De façon générale, si $\phi:Y\rightarrow T$ est une modification finie d'espaces complexes réduits de dimension $m$, on a, pour tout faisceau cohérent ${\mathcal F}$ et   tout faisceau cohérent sans torsion ${\mathcal G}$ sur $Y$, une suite exacte courte
$$0\rightarrow {\mathcal N}\rightarrow {\phi}^{*}{\phi}_{*}{\mathcal F}\rightarrow {\mathcal F}\rightarrow 0$$
et, donc, la suite exacte
$$0\rightarrow{\mathcal H}om({\mathcal F}, {\mathcal G})\rightarrow{\mathcal H}om({\phi}^{*}{\phi}_{*}{\mathcal F}, {\mathcal G})\rightarrow{\mathcal H}om({\mathcal N}, {\mathcal G})$$
Comme ${\mathcal N}$ est à support dans un fermé d'intérieur vide (car $\phi$ est une modification) et que  ${\mathcal G}$ est sans torsion, on obtient l'isomorphisme
$${\mathcal H}om({\mathcal F}, {\mathcal G})\simeq {\mathcal H}om({\phi}^{*}{\phi}_{*}{\mathcal F}, {\mathcal G})$$
auquel on applique le foncteur ${\phi}_{*}$ et la formule d'adjonction pour en déduire l'isomorphisme
$$\phi_{*}{\mathcal H}om({\mathcal F}, {\mathcal G})\simeq {\mathcal H}om({\phi}_{*}{\mathcal F}, \phi_{*}{\mathcal G})$$
Par ailleurs, $\phi$ étant fini, on a aussi l'isomorphisme (de dualité)
$$\phi_{*}{\mathcal H}om({\mathcal F}, \omega^{m}_{Y})\simeq {\mathcal H}om({\phi}_{*}{\mathcal F}, \omega^{m}_{T})$$
On en déduit aisément les isomorphismes voulus sachant que $\phi_{*}{\mathcal L}^{i}_{X'}={\mathcal L}^{i}_{X}$.\vspace{1mm}

\noindent
Du morphisme canonique et injectif
$$\phi^{*}(\widetilde{{\Omega}}^{m-j}_{X})/{\mathcal T}_{\phi}\rightarrow \widetilde{{\Omega}}^{m-j}_{X'}$$
résulte le morphisme 
$${\mathcal H}om({\widetilde{\Omega}}^{m-j}_{X'}, {\mathcal L}^{m}_{{X'}})\rightarrow {\mathcal H}om(\phi^{*}({\widetilde{\Omega}}^{m-j}_{X})/{\mathcal T}_{\phi}, {\mathcal L}^{m}_{{X'}})$$
nécessairement injectif puisque bijectif entre faisceaux sans torsion. Comme dans {\bf(i)}, on a
$${\mathcal H}om(\phi^{*}({\widetilde{\Omega}}^{m-j}_{X})/{\mathcal T}_{\phi}, {\mathcal L}^{m}_{{X'}})\simeq {\mathcal H}om(\phi^{*}({\widetilde{\Omega}}^{m-j}_{X}), {\mathcal L}^{m}_{{X'}})$$
il nous suffit d'appliquer le foncteur exact à gauche  $\phi_{*}$ pour l'injection désirée.\vspace{2mm}

\noindent 
{\bf(iii)} On peut toujours supposer donné un diagramme commutatif de modifications propres
$$\xymatrix{&\tilde{X}\ar[ld]_{\theta}\ar[rd]^{\pi}&\\
X'\ar[rr]_{\phi}&&X}$$
dans lequel $\phi$ est la normalisation de $X$ et $\tilde X$ une variété lisse.\vspace{1mm}

\noindent
Comme $X'$ est normal, on a 
$${\mathcal H}om_{{\mathcal O}_{X'}}({\mathcal L}^{m-j}_{X'}, {\omega}^{m}_{X'})\simeq{\mathcal H}om_{{\mathcal O}_{X'}}(\overline{\Omega}^{m-j}_{X'},\omega^{m}_{X'})\simeq {\mathcal H}om_{{\mathcal O}_{X'}}(\widetilde{\Omega}^{m-j}_{X'},\omega^{m}_{X'})$$
puisque ceux sont des faisceaux de profondeur au moins deux coincidant  en dehors du lieu singulier de $X'$ qui est de codimension au moins deux.\vspace{1mm}

\noindent
On en déduit les isomorphismes
$${\mathcal H}om_{{\mathcal O}_{X}}({\mathcal L}^{m-j}_{X}, {\omega}^{m}_{X})\simeq{\mathcal H}om_{{\mathcal O}_{X}}(\phi_{*}\overline{\Omega}^{m-j}_{X'},\omega^{m}_{X})\simeq {\mathcal H}om_{{\mathcal O}_{X}}(\phi_{*}\widetilde{\Omega}^{m-j}_{X'},\omega^{m}_{X})$$
Comme $\phi$ est une modification,  $\phi_{*}\overline{\Omega}^{m-j}_{X'}$ est un sous faisceau de $\overline{\Omega}^{m-j}_{X}$. On s'en convainc aisément en appliquant  $\phi_{*}$ aux  injections
$$\overline{\Omega}^{m-j}_{X'}\subset {\mathcal L}^{m-j}_{X'}\subset \omega^{m-j}_{X'}$$
donnant
$$\xymatrix{\phi_{*}\overline{\Omega}^{m-j}_{X'}\ar@{^{(}->}[r]& \phi_{*}{\mathcal L}^{m-j}_{X'}\eq[d]\ar@{^{(}->}[r]&\phi_{*}\omega^{m-j}_{X'}\ar@{^{(}->}[d]\\
\overline{\Omega}^{m-j}_{X}\ar@{^{(}->}[r]\ar@{^{(}->}[u]& {\mathcal L}^{m-j}_{X}\ar@{^{(}->}[r]&\omega^{m-j}_{X}}$$
et il est facile de voir  que l'image  $\phi_{*}\overline{\Omega}^{m-j}_{X'}$ contient  $\overline{\Omega}^{m-j}_{X}$ (cf \cite{B5}, {\bf{Lem.5.2.3}}, p. 72).\vspace{1mm}

\noindent 
En effet, l'image réciproque $\phi^{*}\widetilde{\Omega}^{k}_{X}\rightarrow \widetilde{\Omega}^{k}_{X'}$ passe au quotient et induit le morphisme injectif ${\mathcal A}:=(\phi^{*}(\widetilde{\Omega}^{k}_{X}))/{{\mathcal T}_{\phi}}\rightarrow \widetilde{\Omega}^{k}_{X'}$ auquel on applique le foncteur $\theta^{*}$ et quotientant par les sous faisceaux de  torsion respectifs, on obtient le morphisme injectif
$$\theta^{*}({\mathcal A})/{{\mathcal T}_{\theta}}\rightarrow \theta^{*}(\widetilde{\Omega}^{k}_{X'})/{{\mathcal T}'_{\theta}}$$
(dont le conoyau est supporté par le diviseur exceptionnel) qui induit, à son tour, les injections 
$$\theta_{*}(\theta^{*}({\mathcal A})/{{\mathcal T}_{\theta}})\subset\overline{\Omega}^{k}_{X'}$$
et
$$\phi_{*}(\theta_{*}(\theta^{*}({\mathcal A})/{{\mathcal T}_{\theta}}))\subset\phi_{*}(\overline{\Omega}^{k}_{X'})$$
Mais, par construction, on a $$\overline{\Omega}^{k}_{X}:=\pi_{*}(\pi^{*}(\widetilde{\Omega}^{k}_{X})/{\mathcal T}_{\pi}) \simeq\phi_{*}(\theta_{*}(\theta^{*}({\mathcal A})/{{\mathcal T}_{\theta}}))$$
et la conclusion qui s'en suit avec le diagramme
$$\xymatrix{{\mathcal D}_{\omega_{X}}({\mathcal L}^{m-j}_{X})\eq[r]&
{\mathcal D}_{\omega_{X}}(\phi_{*}\overline{\Omega}^{m-j}_{X'})\ar@{^{(}->}[d]\eq[r]& {\mathcal D}_{\omega_{X}}(\phi_{*}\widetilde{\Omega}^{m-j}_{X'})\ar@{^{(}->}[d]\ar@{^{(}->}[r]&\omega^{j}_{X}\\
&{\mathcal D}_{\omega_{X}}(\overline{\Omega}^{m-j}_{X})\ar@{^{(}->}[r]& {\mathcal D}_{\omega_{X}}(\widetilde{\Omega}^{m-j}_{X})&}$$

\vspace{1mm}

\noindent 
{\bf(iv)} En ce qui concerne le cup produit, il suffit d'écrire
$${\mathcal H}om(\widetilde{\Omega}^{k'}_{X}, {\mathcal H}om(\widetilde{\Omega}^{k}_{X}, {\mathcal L}^{m}_{X})\simeq {\mathcal H}om(\widetilde{\Omega}^{k'}_{X}\otimes \widetilde{\Omega}^{k}_{X}, {\mathcal L}^{m}_{X})$$ 
et d'utiliser le cup produit classique $\widetilde{\Omega}^{k'}_{X}\otimes \widetilde{\Omega}^{k}_{X}\rightarrow \widetilde{\Omega}^{k+k'}_{X}$ pour en déduire le morphisme
$${\mathcal H}om(\widetilde{\Omega}^{k+k'}_{X}, {\mathcal L}^{m}_{X})\rightarrow {\mathcal H}om(\widetilde{\Omega}^{k'}_{X}, {\mathcal H}om(\widetilde{\Omega}^{k}_{X}, {\mathcal L}^{m}_{X}) $$
et, par suite, le cup produit recherché. On raisonne de la même façon pour les deux autres en utilisant, cette fois, le cup produit
$${\mathcal L}^{k'}_{X}\otimes {\mathcal L}^{k}_{X}\rightarrow {\mathcal L}^{k+k'}_{X}$$
On obtient, alors, les flèches 
$${\mathcal D}_{{\mathcal L}_{X}}(\widetilde{\Omega}^{k+k'}_{X})\otimes \Omega^{k'}_{X}\rightarrow {\mathcal D}_{{\mathcal L}_{X}}(\widetilde{\Omega}^{k}_{X})$$
$${\mathcal D}_{{\mathcal L}_{X}}({\mathcal L}^{k+k'}_{X})\otimes{\mathcal L}^{k'}_{X}\rightarrow {\mathcal D}_{{\mathcal L}_{X}}({\mathcal L}^{k}_{X}),\,\,\, {\mathcal D}_{{\omega}_{X}}({\mathcal L}^{k+k'}_{X})\otimes{\mathcal L}^{k'}_{X}\rightarrow {\mathcal D}_{{\omega}_{X}}({\mathcal L}^{k}_{X})$$
et, par suite, en tenat compte des notations, celles de l'énoncé$\,\blacksquare$
\end{proof} \vspace{2mm}

\noindent 
La proposition qui suit concerne des faisceaux particuliers que l'on aura à taiter dans la suite;
\Prop{}{}\label{pro2}  Soit $X$ un espace analytique réduit et  ${\mathcal F}_{X}$ un faisceau cohérent sans torsion dont la construction est compatible aux projections\footnote{toute projection $p:X\times U\rightarrow X$ induit un morphisme naturel $p^{*}{\mathcal F}_{X}\rightarrow {\mathcal F}_{X\times U}$ faisant de $p^{*}{\mathcal F}_{X}\otimes \Omega^{n}_{U}$ un facteur direct de ${\mathcal F}_{X\times U}$  . } et muni des inclusions  $\widetilde{\Omega}^{\bullet}_{X}\subset {\mathcal F}_{X}\subset \omega^{\bullet}_{X}$ bijectives sur la partie régulière ${\rm Reg}(X)$ de $X$. Alors,  si ${\mathcal F}_{X}$ est stable par image réciproque équidimensionnelle ou ouverte (resp. image directe finie), ${\mathcal D}_{\omega_{X}}({\mathcal F}_{X})$ est stable par image directe ou par trace (resp. image réciproque).\rm\vspace{1mm}

\noindent
\begin{proof} {\bf(i)} Supposons la construction de ${\mathcal F}_{X}$  stable par image réciproque équidimensionnelle. Cela signifie que tout morphisme $f:X\rightarrow Y$ surjectif et à fibres de dimension constante $n$, induit un morphisme  ${\mathcal O}_{X}$-linéaire $f^{*}{\mathcal F}_{Y}\rightarrow {\mathcal F}_{X}$ prolongeant l'image réciproque des formes holomorphes usuelles de sorte à avoir la commutativité  des diagrammes
$$\xymatrix{f^{*}\widetilde{\Omega}^{\bullet}_{Y}\ar[r]\ar[d]&f^{*}{\mathcal F}_{Y}\ar[d]&{\rm ou}&\widetilde{\Omega}^{\bullet}_{Y}\ar@{^{(}->}[r]\ar[d]&{\mathcal F}_{Y}\ar[d]\\
\widetilde{\Omega}^{\bullet}_{X}\ar[r]&{\mathcal F}_{X}&&f_{*}\widetilde{\Omega}^{\bullet}_{X}\ar@{^{(}->}[r]&f_{*}{\mathcal F}_{X}}$$
Alors, pour $f$ fini et surjectif, on a le diagramme commutatif
$$\xymatrix{f_{*}{\mathcal D}_{\omega_{X}}({\mathcal F}_{X})\eq[r]&{\mathcal H}om({f_{*}}{\mathcal F}_{X}, \omega^{m}_{Y})\ar[d]_{\Psi}\ar@{^{(}->}[r]&{\mathcal H}om({f_{*}}\widetilde{\Omega}^{\bullet}_{X}, \omega^{m}_{Y})\ar[d]\eq[r]&f_{*}\omega^{m-\bullet}_{X}\ar[d]^{{\rm Tr}_{f}}\\
{\mathcal D}_{\omega_{Y}}({\mathcal F}_{Y})\eq[r]&{\mathcal H}om({\mathcal F}_{Y}, \omega^{m}_{Y})\ar@{^{(}->}[r]&{\mathcal H}om(\widetilde{\Omega}^{\bullet}_{Y}, \omega^{m}_{Y})\eq[r]
&\omega^{m-\bullet}_{Y}}$$
qui montre que $\Psi$ est un morphisme de type trace puisque naturellement induit par la trace usuelle (ou l'image directe au sens des courants)  $f_{*}\omega^{m-\bullet}_{X}\rightarrow\omega^{m-\bullet}_{Y}$ sachant qu'il coincide génériquement avec ce dernier et que les faisceaux sont sans torsion.\vspace{1mm}

\noindent
 Il en sera, donc, ainsi si  ${\mathcal F}_{X}$ est l'un des faisceaux $$\widetilde{\Omega}^{\bullet}_{X}\subset \overline{\Omega}^{\bullet}_{X}\subset\widehat{\Omega}^{\bullet}_{X}\subset {\mathcal L}^{\bullet}_{X}$$
  {\bf(ii)} Supposons, à présent, que  ${\mathcal F}_{X}$ soit stable par  image directe finie ce qui signifie, pour nous, qu'en tout point $x$ de $X$ et toute paramétrisation locale   $f:X\rightarrow T$, le morphisme trace habituel $f_{*}\omega^{\bullet}_{X}\rightarrow\omega^{\bullet}_{T} $ induit un morphisme trace $f_{*}{\mathcal F}_{X}\rightarrow {\mathcal F}_{T}$. \vspace{1mm}

  \noindent
Si $\pi:X\rightarrow Y$ est un morphisme à fibres de dimension constante, on peut le factoriser localement (comme de coutume) en un morphisme fini et une projection $\displaystyle{\xymatrix{X\ar@/_1pc/[rr]_{\pi}\ar[r]_{f}&Y\times U\ar[r]_{q}&Y}}$. \vspace{1mm}

\noindent Comme on le verra, notre intérêt porte sur les faisceaux ${\mathcal F}_{X}$ coincidant génériquement avec le faisceau des $(n+k)$-formes holomorphes usuelles. \vspace{1mm}

\noindent 
On veut montrer que toute section $\xi$ de ${\mathcal D}_{\omega_{Y}}({\mathcal F}_{Y})$ donne par image réciproque une section $\pi^{*}(\xi)$ de ${\mathcal D}_{\omega_{X}}({\mathcal F}_{X})$.\vspace{1mm}

\noindent
Vu la factorisation de ces morphismes, il nous suffit de nous en convaincre pour une projection et un morphisme fini   (ouvert et surjectif). 
Dans le cas de la projection, la nature particulière de ces faisceaux pour lesquels on a des morphismes de projection $${\mathcal F}_{Y\times U}\rightarrow q^{*}{\mathcal F}_{Y}\otimes p^{*}\Omega^{n}_{U}$$
permet de mettre en évidence une image réciproque
$$q^{*}{\mathcal D}_{\omega_{Y}}({\mathcal F}_{Y})\rightarrow{\mathcal D}_{\omega_{Y}}({\mathcal F}_{Y})$$  construite par le biais du diagramme commutatif
$$\xymatrix{ q^{*}{\mathcal H}om({\mathcal F}_{Y}, \omega^{r}_{Y})\ar[d]\eq[r]&{\mathcal H}om(q^{*}{\mathcal F}_{Y}, q^{*}\omega^{r}_{Y})\eq[r]&{\mathcal H}om(q^{*}{\mathcal F}_{Y}\otimes p^{*}\Omega^{n}_{U}, q^{*}\omega^{r}_{Y}\otimes p^{*}\Omega^{n}_{U})\eq[d]\\
{\mathcal H}om({\mathcal F}_{Y\times U}, \omega^{m}_{Y\times U})&&{\mathcal H}om(q^{*}{\mathcal F}_{Y}\otimes p^{*}\Omega^{n}_{U}, \omega^{m}_{Y\times U})\ar[ll]}$$
On peut, alors, se concentrer sur le cas d'un morphisme fini.
Considérons, alors, une section  $\xi$ de  ${\mathcal H}om({\mathcal F}_{Y}, \omega^{m}_{Y\times U})$. Si $\alpha$ est une section arbitraire  de ${\mathcal F}_{X}$, la forme $f^{*}(\xi)\wedge \alpha$ vérifie
$${\rm Tr}_{f}(f^{*}(\xi)\wedge \alpha)={\rm deg}(f) \xi\wedge {\rm Tr}_{f}(\alpha)$$
Comme, par hypothèse, ${\rm Tr}_{f}(\alpha)$ définit une section du faisceau ${\mathcal F}_{Y\times U}$, on a, par définition, $\xi\wedge {\rm Tr}_{f}(\alpha)\in \omega^{m}_{Y\times U}$. Mais cela équivaut  au fait que $f^{*}(\xi)\wedge \alpha$  est une section du faisceau  $\omega^{m}_{X}$, c'est-à-dire $f^{*}(\xi)\in{\mathcal H}om({\mathcal F}_{X}, \omega^{m}_{X})$. D'où  une image réciproque
$${\mathcal D}_{\omega_{Y\times U}}({\mathcal F}_{Y\times U})\rightarrow f_{*}{\mathcal D}_{\omega_{X}}({\mathcal F}_{X})$$
et, par suite, le morphisme image réciproque recherché
$$\xymatrix{\pi^{*}{\mathcal D}_{\omega_{Y}}({\mathcal F}_{Y})=f^{*}(q^{*}{\mathcal D}_{\omega_{Y}}({\mathcal F}_{Y}))\ar[rr]\ar[rd]&&f^{*}({\mathcal D}_{\omega_{Y\times U}}({\mathcal F}_{Y\times U}))\ar[ld]\\
&{\mathcal D}_{\omega_{X}}({\mathcal F}_{X})&}$$
On peut remarquer, au passage, que, par hypothèse, le faisceau ${\mathcal F}_{Y\times U}$ est un facteur direct du faisceau $f_{*}{\mathcal F}_{X}$ (en raison de la trace et du pull back) et que, par conséquent, ${\mathcal D}_{\omega_{Y\times U}}({\mathcal F}_{Y\times U})$ est aussi facteur direct de  $f_{*}{\mathcal D}_{\omega_{X}}({\mathcal F}_{X})$. On a aussi un diagramme commutatif 
$$\xymatrix{ q^{*}{\mathcal D}_{\omega_{Y}}({\mathcal F}_{Y})\ar[dd]\eq[r]&{\mathcal H}om(q^{*}{\mathcal F}_{Y}, q^{*}\omega^{r}_{Y})\ar[rd]&\\
&{\mathcal D}_{\omega_{Y\times U}}({\mathcal F}_{Y\times U})\eq[r]&{\mathcal H}om({\mathcal F}_{Y\times U}, \omega^{m}_{Y\times U})\ar[ld]\\
f_{*}{\mathcal D}_{\omega_{X}}({\mathcal F}_{X})\eq[r]&{\mathcal H}om(f_{*}{\mathcal F}_{X}, \omega^{m}_{Y\times U})&}$$

On peut remarquer, au passage, qu'étant donné que  le faisceau ${\mathcal L}^{\bullet}_{Z}$ est stable par image réciproque et directe pour tout espace complexe réduit $Z$, son $\omega_{Z}$-dual l'est  aussi\,$\blacksquare$
\end{proof}
\vspace{1mm}

\noindent
\begin{rem} En général, l'inclusion $${\mathcal H}om_{{\mathcal O}_{X}}(\overline{\Omega}^{m-j}_{X},\omega^{m}_{X})\subset \omega^{j}_{X}$$
est stricte. En effet, même dans le cas très spécial où la normalisation est une désingularisation $\pi:\tilde{X}\rightarrow X$, on a
$${\mathcal H}om_{{\mathcal O}_{X}}(\overline{\Omega}^{m-j}_{X},\omega^{m}_{X})\simeq \pi_{*}{\mathcal H}om_{{\mathcal O}_{\tilde X}}(\pi^{*}(\widetilde{\Omega}^{m-j}_{X})/{{\mathcal T}_{\pi}}, {\Omega}^{m}_{\tilde{X}})\simeq\pi_{*}{\mathcal H}om_{{\mathcal O}_{\tilde X}}(\pi^{*}(\widetilde{\Omega}^{m-j}_{X}), {\Omega}^{m}_{\tilde{X}})$$
dont on déduit, grâce à la formule d'adjonction:
$${\mathcal H}om_{{\mathcal O}_{X}}(\overline{\Omega}^{m-j}_{X},\omega^{m}_{X})\simeq{\mathcal H}om_{{\mathcal O}_{X}}(\widetilde{\Omega}^{m-j}_{X}, {\mathcal L}^{m}_{X})$$ 
Alors, si l'inclusion n'était pas stricte, on aurait l'isomorphisme
$${\mathcal H}om_{{\mathcal O}_{X}}(\widetilde{\Omega}^{m-j}_{X}, {\mathcal L}^{m}_{X})\simeq {\mathcal H}om_{{\mathcal O}_{X}}(\widetilde{\Omega}^{m-j}_{X}, {\omega}^{m}_{X})$$
qui, pour $j=m$,  donne l'égalité ${\mathcal L}^{m}_{X}=\omega^{m}_{X}$ qui, non seulement, impose la normalité de l'espace (que l'on a pas ici !) mais aussi d'être en présence de singularités très spéciales (par exemple si $X$ est de Cohen Macaulay, elles seront rationnelles!).\end{rem}\vspace{1mm}

\noindent
\section{\color{blue}{Remarques générales et exemples.}}
\subsection{Les inclusions 
$\widetilde{\Omega}^{\bullet}_{X}\subsetneq\overline{\Omega}^{\bullet}_{X}\subsetneq {\mathcal
L}^{\bullet}_{X}\subsetneq \omega^{\bullet}_{X}$:}\vspace{1mm}

\par \noindent
Ces inclusions généralement strictes induisent les inclusions 
$${\mathcal L}^{j}_{X}\subsetneq {\mathcal H}om_{{\mathcal O}_{X}}({\mathcal L}^{m-j}_{X},
{\mathcal L}^{m}_{X})\subsetneq {\mathcal H}om_{{\mathcal O}_{X}}(\Omega^{m-j}_{X},
{\mathcal L}^{m}_{X})\subsetneq{\mathcal H}om_{{\mathcal O}_{X}}\overline{(\Omega}^{m-j}_{X},
{\mathcal L}^{m}_{X})\subsetneq \omega^{j}_{X}$$
encore  strictes comme on peut le voir sur certains exemples simples.\vspace{1mm}

\noindent
\subsubsection{Espaces  faiblement normaux.} On commence par le cas classique du cusp pour lequel les premières inclusions sont strictes mais pas les secondes.\vspace{1mm}

\noindent 
\centerline{{\bf(i)} $X:=\{(x,y)/ x^{2} - y^{3} =0\}$}\vspace{1mm}

\noindent muni de son morphisme de normalisation donn\'e par l'application $\nu:{\Bbb C}\rightarrow{\Bbb C}^{2}$ qui \`a $t$ associe $\nu(t):=(t^{3}, t^{2})$.\vspace{1mm}

\noindent
On a les inclusions strictes
$$\overline{\Omega}^{j}_{X}\subsetneq{\mathcal L}^{j}_{X}\subsetneq\omega^{j}_{X},\,\,j=0,1$$
puisque, par exemple, $\displaystyle{y\over{x}}$ définit une section de  $\omega^{0}_{X}$ qui n'est pas dans ${\mathcal L}^{0}_{X}$ ( car $\nu^{*}({y\over{x}})={1\over{t}}$ non holomorphe!) qui, s'identifie, d'ailleurs à $\overline{\Omega}^{0}_{X}$. De même, on vérifie aisément que $\displaystyle{dy\over{x}}$ définit une section de $\omega^{1}_{X}$ mais non de ${\mathcal L}^{1}_{X}$. La forme méromorphe $\displaystyle{dx\over{y}}$, quant à elle, définit une section de ${\mathcal L}^{1}_{X}$ mais non de $\overline{\Omega}^{1}_{X}$. D'autre part, un calcul simple montre que  
$${\mathcal H}om_{{\mathcal O}_{X}}(\Omega^{1}_{X},
{\mathcal L}^{1}_{X}) ={\mathcal H}om_{{\mathcal O}_{X}}(\overline{\Omega}^{1}_{X}
{\mathcal L}^{1}_{X})={\mathcal H}om_{{\mathcal O}_{X}}(\Omega^{1}_{X},
{\omega}^{1}_{X})=\omega^{0}_{X}$$
il suffit pour cela de vérifier que  la section $\displaystyle{y\over{x}}$ de $\omega^{0}_{X} $  d\'efinit bien une section de ${\mathcal H}om_{{\mathcal O}_{X}}(\Omega^{1}_{X}, {\mathcal L}^{1}_{X})$. 
Par ailleurs, le morphisme de normalisation étant fini et la normalisée lisse, on a
$${\mathcal L}^{j}_{X}\simeq {\mathcal H}om_{{\mathcal O}_{X}}({\mathcal L}^{1-j}_{X},
{\mathcal L}^{1}_{X})\simeq {\mathcal H}om_{{\mathcal O}_{X}}({\mathcal L}^{1-j}_{X},
{\omega}^{1}_{X}),\,\forall\,j\in\{0,1\}$$
\vspace{1mm}

\indent
\centerline{{\bf(ii)} $X=\{(x,y,z)/ xy^{2} = z^{2}\}$}\vspace{3mm}

\noindent
La  normalisation est donnée  par le morphisme
$\nu: {\Bbb C}^{2}\rightarrow X$ qui \`a $(u,v)$ associe $(u^{2}, v,
vu)$. 
On a les inclusions strictes 
$$\overline{\Omega}^{j}_{X}\subsetneq{\mathcal L}^{j}_{X}\subsetneq\omega^{j}_{X},\,\,j=1,2$$ 
que l'on constate en remarquant  que la forme $\displaystyle{\frac{dx\wedge dy}{z}}$ est une section de $\omega^{2}_{X}$ mais pas de ${\mathcal L}^{2}_{X}$, la forme $\displaystyle{dz\over{y}}$ définit  une section de $\omega^{1}_{X}$ mais pas de ${\mathcal L}^{1}_{X}$ (car sa relevée est non holomorphe sur la normalisée!) et n'est pas non plus une section du faisceau 
 ${\mathcal H}om_{{\mathcal O}_{X}}(\Omega^{1}_{X},{\mathcal L}^{2}_{X})$  puisque le cup produit par $dx$ donne  une section du faisceau $\omega^{2}_{X}$ qui n'est pas dans 
 ${\mathcal L}^{2}_{X}$. De plus, la forme $\displaystyle{\frac{y}{z}dx}$ est une section de ${\mathcal L}^{1}_{X}$ mais pas de $\overline{\Omega}^{1}_{X}$.\vspace{1mm}

 \noindent
On montre facilement que $\overline{\Omega}^{0}_{X}={\mathcal L}^{0}_{X}=\omega^{0}_{X}$ et, donc
$${\mathcal L}^{0}_{X}={\mathcal H}om_{{\mathcal O}_{X}}({\mathcal L}^{2}_{X},
{\mathcal L}^{2}_{X})\simeq {\mathcal H}om_{{\mathcal O}_{X}}(\Omega^{2}_{X},
{\mathcal L}^{2}_{X})\simeq {\mathcal H}om_{{\mathcal O}_{X}}(\overline{\Omega}^{2}_{X},
{\omega}^{2}_{X})=\omega^{0}_{X}$$
Comme ci-dessus, la normalisée étant lisse, on a 
$${\mathcal L}^{i}_{X}= {\mathcal H}om_{{\mathcal O}_{X}}({\mathcal L}^{2-i}_{X},
{\mathcal L}^{2}_{X})\simeq{\mathcal H}om_{{\mathcal O}_{X}}({\mathcal L}^{2-j}_{X},
{\omega}^{2}_{X})$$
\vspace{1mm}
 
 \indent
\centerline{{\bf(iii)} $X=\{(x,y,z)/ xy^{2} = z^{3}\}$}\vspace{3mm}

\noindent
Le morphisme de normalisation est donné par le morphisme fini  $\nu: {\Bbb C}^{2}\rightarrow X$ qui \`a $(u,v)$ associe $(u^{3}, v^{3}, uv^{2})$. On a 
$${\mathcal L}^{0}_{X}\subsetneq\omega^{0}_{X}={\mathcal H}om_{{\mathcal O}_{X}}(\Omega^{2}_{X},
{\mathcal L}^{2}_{X})$$
comme on s'en convainc aisément en vérifiant que  la section $\displaystyle{z\over{y}}$ de $\omega^{0}_{X}$ ( et qui n'est pas dans $ {\mathcal L}^{0}_{X}$) definit bien une section du faisceau ${\mathcal H}om_{{\mathcal O}_{X}}(\Omega^{2}_{X}, {\mathcal L}^{2}_{X})$ contenant strictement le faisceau   ${\mathcal L}^{0}_{X}$ dont un générateur (sur ${\mathcal O}_{X}$) est $\displaystyle{\frac{xy}{z}}$. Pour les mêmes raisons que ci-dessus, on a 
$${\mathcal L}^{i}_{X}= {\mathcal H}om_{{\mathcal O}_{X}}({\mathcal L}^{2-i}_{X},
{\mathcal L}^{2}_{X})$$
et les inclusions strictes
$${\mathcal L}^{1}_{X}\subsetneq {\mathcal H}om_{{\mathcal O}_{X}}(\Omega^{1}_{X}, {\mathcal L}^{2}_{X})\subsetneq \omega^{1}_{X}$$
comme on peut le constater en prenant, par exemple, les formes m\'eromorphes:\vspace{1mm}

\indent $\bullet$ $\displaystyle{dz\over{y}}$ et $\displaystyle{\frac{zdy}{y^{2}}}$  définissent  des sections du faisceau $\omega^{1}_{X}$ mais pas du  faisceau ${\mathcal L}^{1}_{X}$. Ceux ne sont pas non plus des sections du faisceau  ${\mathcal H}om_{{\mathcal O}_{X}}(\Omega^{1}_{X}, {\mathcal L}^{2}_{X})$ car leur cup produit par la forme $dx$ donne, par image réciproque sur la normalisée (au coefficient près), la forme   $\displaystyle{{u^{3}\over{v^{2}}}du\wedge dv}$ qui n'est manifestement pas holomorphe sur ${\Bbb C}^{2}$.\vspace{1mm}

\indent
$\bullet$ $\displaystyle{\frac{ydx}{z^{2}}}$ est une section de $\omega^{1}_{X}$ mais non de  ${\mathcal L}^{1}_{X}$ (puisque d'image r\'eciproque $\displaystyle{\frac{3du}{v}}$) mais définit une section de ${\mathcal H}om_{{\mathcal O}_{X}}(\Omega^{1}_{X}, {\mathcal L}^{2}_{X})$ puisque les images réciproques $\nu^{*}(\xi\wedge dx)$,   $\nu^{*}(\xi\wedge dy)$ et
$\nu^{*}(\xi\wedge dz)$ sont holomorphes sur ${\Bbb C}^{2}$.\vspace{1mm}

\vspace{1mm}

\noindent 
La relation différentielle
$$y^{2}dx +2xydy=3z^{2}dz$$
nous donne
$$\frac{zdx}{x} + \frac{zdy}{y}=3dz$$
avec des formes méromorphes $\frac{zdx}{x}$ et $\frac{zdy}{y}$ sections du faisceau ${\mathcal L}^{1}_{Z}$ puisque leurs images réciproques par la normalisation  sont holomorphes. Mais elles ne définissent pas des sections de $\overline{\Omega}^{1}_{X}$.

Bien que leur somme soit holomorphe, leur produit tensoriel a pour coefficient la fonction méromorphe $\frac {z^{2}}{xy}$ qui ne se prolonge pas holomorphiquement sur $Z$ car son image réciproque dans la normalisée est donnée par la forme méromorphe $\frac{v}{u}$; on peut donc garantir qu'elles ne peuvent pas vérifier  d'équations de dépendance intégrale de degré $2$ et, d'ailleurs, d'aucun degré vu la relation différentielle minimale et l'équation de définition de $X$. 
Relation minimale qui montre que ni l'une ni l'autre ne peut s'écrire comme somme de produit de formes holomorphes à coeficients génériquement continus et localement bornés!
\vspace{1mm}

\noindent
 \centerline{{\bf(iv)} $X=\{(x,y,z,t)\in {\Bbb C}^{4}/ xy^{2}-tz^{3}=0\}$}\vspace{3mm}

 \noindent $X$ est une hypersurface de ${\Bbb C}^{4}$ non normale puisque de lieu singulier (réduit) ${\rm Sing}(X):=\{(x,y,z,t)\in X: y=z=0\}$ de codimension $1$ dans $X$.  \vspace{1mm}
 
 \noindent
 On peut voir, par exemple, que la forme $\displaystyle{\frac{xdy-ydx}{z^{2}}}$ est une section  de $\omega^{1}_{X}$ et non de ${\mathcal L}^{1}_{X}$; les calculs sont facilités par le fait que l'éclatement le long du lieu singulier donne une carte lisse et une à singularités canoniques.  En regardant les préimages des cup-produit de cette forme par les éléments différentiels d'ordre deux, on vérifie sans difficulté qu'elle définit une section du faisceau ${\mathcal H}om_{{\mathcal O}_{X}}(\Omega^{2}_{X}, {\mathcal L}^{3}_{X})$. Enfin, en examinant la préimage de son cup-produit par la forme $\displaystyle{\frac{z}{y}dt\wedge dz}$ qui est une section du faisceau ${\mathcal L}^{2}_{X}$, on constate qu'elle  n'est pas section du faisceau ${\mathcal H}om_{{\mathcal O}_{X}}({\mathcal L}^{2}_{X}, {\mathcal L}^{3}_{X})$. \vspace{1mm}

 \noindent
 Le lecteur pourra vérifier que  la forme
 $$\xi:=\frac{dx\wedge dy}{z} -\frac{ydx\wedge dz}{z^2} -\frac{xdy\wedge dz}{z^2}$$
  qui est une section de $\omega^{2}_{X}$ mais pas ${\mathcal L}^{2}_{X}$ est, en fait, une section de ${\mathcal H}om_{{\mathcal O}_{X}}({\mathcal L}^{1}_{X}, {\mathcal L}^{3}_{X})$. 

\subsubsection{Espaces normaux.}\vspace{1mm}

\noindent
Remarquons que si $X$ est normal et ${\mathcal F}$  un faisceau cohérent sur $X$  tel que
$\widetilde{\Omega}^{m-k}_{X}\subset {\mathcal F}\subset  \omega^{m-k}_{X}$, on a l'isomorphisme naturel 
$$\omega^{k}_{X}\simeq {\mathcal H}om({\mathcal F}, \omega^{m}_{X})$$
dû au fait que ${\mathcal F}$ coincide nécessairement avec le faisceau $\Omega^{m-k}_{X}$ sur la partie régulière, que $\omega^{m}_{X}$ est de profondeur au moins deux et que, $X$ étant normal, son lieu singulier est de codimension au moins deux.\vspace{1mm}

\noindent
\centerline{{\bf(i)}  $X=\{(x,y,z)/ x^{3} + y^{3} = z^{4}\}$}\vspace{3mm}

\noindent est une surface normale à singularité isolée en l'origine. 
Comme ${1\over{3}} +{1\over{3}} + {1\over{4}}<1$, la singularité n'est pas rationnelle et, donc, l'inclusion ${\mathcal L}^{2}_{X}\subset \omega^{2}_{X}$ est stricte (vérification immédiate en relevant le générateur $\displaystyle{\frac{dy\wedge dz}{x^{2}}}$ de $\omega^{2}_{X}$).\vspace{1mm}

\noindent
On vérifie facilement que la forme $\xi_{1}:=\displaystyle{\frac{xdy -ydx}{z^{3}}}$ définit une section du faisceau $\omega^{1}_{X}$ mais pas de ${\mathcal L}^{1}_{X}$ (en regardant les images réciproques dans les différentes cartes des éclatements).\vspace{1mm}

\noindent
Plus précisémment,  $\xi$ définit une section du faisceau  ${\mathcal H}om_{{\mathcal O}_{X}}(\Omega^{1}_{X}, {\mathcal L}^{2}_{X})\simeq {\mathcal H}om_{{\mathcal O}_{X}}({\mathcal L}^{1}_{X}, {\mathcal L}^{2}_{X}) $. Là encore, on s'en convainc aisément en prenant les images réciproques dans les différents éclatements des formes $\xi\wedge dx$, $\xi\wedge dy$
 et $\xi\wedge dz$. \vspace{1mm}

\noindent
Même conclusions pour  $X:=\{(x,y,z)\in {\Bbb C}^{3}/ x^{4} + y^{4} + z^{5}=0\}$ une surface normale de ${\Bbb C}^{3}$ qui n'est pas à singularité rationnelle et la forme méromorphe $\xi:= \displaystyle{xdy - ydx\over{z^{4}}}$. 
$\bullet$ Avec des calculs un peu plus fastidieux, on montre que, pour la surface normale \vspace{1mm}

\centerline{{\bf(ii)} $X=\{(x,y,z)/ yx^{3} + zy^{2} +z^{4}=0\}$ }\vspace{3mm}

\noindent qui n'est pas à singularité rationnelle puisque le relèvement de la forme $\omega_{0}:=\displaystyle{\frac{dy\wedge dz}{3yx^{2}}}$ ne se prolonge   pas  holomorphiquement sur la désingularisée, la forme méromorphe $\xi:=\displaystyle{\frac{3ydz -2zdy}{yx^{2}}}$ définit une section du faisceau $\omega^{1}_{X}$ mais pas de ${\mathcal L}^{1}_{X}$. De plus, on vérifie à la suite de  calculs ne présentant aucune  difficulté (mais un peu long!),  que $\xi$ est en fait une section du faisceau  ${\mathcal H}om_{{\mathcal O}_{X}}(\Omega^{1}_{X}, {\mathcal L}^{2}_{X})$ ( on se ramène à étudier les formes $\xi\wedge dx=-5x\omega_{0}$, $\xi\wedge dy=-9y\omega_{0}$ et $\xi\wedge dz=-6z\omega_{0}$ dans les différentes cartes des éclatements). \vspace{1mm}

\noindent 
On remarque au passage que le cup produit de $\xi$ par la forme $\displaystyle{\frac{zdy}{y}}$, section du faisceau ${\mathcal L}^{1}_{X}$, est une section du faisceau ${\mathcal L}^{2}_{X}$. Cela nous incite à penser que $\xi$ pourrait bien définir une section du faisceau ${\mathcal H}om_{{\mathcal O}_{X}}({\mathcal L}^{1}_{X}, {\mathcal L}^{2}_{X})$. Mais pour en être sûr, il faudrait pouvoir déterminer toutes les sections du faisceau ${\mathcal L}^{1}_{X}$; ce qui n'est pas chose facile!\vspace{1mm}

\noindent 
\subsection{{\bf{ L'isomorphisme d'auto-dualité.}}}
Les exemples précédents montrent que si  l'espace complexe  de Gorenstein (ou de Cohen Macaulay) $X$ possède un espace normalis\'e  lisse ou si c'est une $V$-vari\'et\'e (i.e quotient d'une vari\'et\'e lisse sous l'action d'un groupe d'automorphismes  discret sans points fixes), on a l'isomorphisme
$${\mathcal L}^{i}_{X}= {\mathcal H}om_{{\mathcal O}_{X}}({\mathcal L}^{m-i}_{X},
{\mathcal L}^{m}_{X})$$
Il est, alors, légitime  de se demander si c'est le cas, en général,  pour un espace complexe réduit $X$ de dimension $m$. Mais il n'en est rien comme le montre les exemples:\vspace{2mm}

\noindent 
\centerline{$X:=\big\{(x,y,z,t)\in {\Bbb C}^{4}: f(x,y,z,t)=x^4 +y^5 +z^5 +t^5=0\big\}$}\vspace{3mm}

\noindent
hypersurface normale à singularité isolée et non rationnelle puisque $3/5 +1/4<1$. Le faisceau dualisant $\omega^{3}_{X}$ est localement libre engendré par l'une des formes
$$\frac{dy\wedge dz\wedge dt}{4x^3}=-\frac{dx\wedge dz\wedge dt}{5y^4}=\frac{dx\wedge dy\wedge dt}{5z^4}=-\frac{dx\wedge dy\wedge dz}{5t^4}$$
que l'on désignera par $\omega_{1}$, $\omega_{2}$, $\omega_{3}$ et $\omega_{4}$ respectivement.\vspace{1mm}

\noindent
Soit
$$\xi:=\frac{tdx\wedge dy}{5z^4} -\frac{ydx\wedge dt}{5z^4} +\frac{xdy\wedge dt}{4z^4}$$
On vérifie facilement (par cup produit par $\frac{df}{f}$) que $\xi$ est une section du faisceau $\omega^{2}_{X}$. \vspace{1mm}

\noindent
En regardant, l'éclatement en $z$ dont la préimage stricte est lisse, on voit que l'image réciproque de cette forme n'est pas holomorphe puisqu'elle fait apparaitre, après simplification,  un coefficient en $\frac{1}{z}$ qui est non holomorphe sur cette préimage; ce qui signifie qu'elle ne définit pas une section du faisceau   ${\mathcal L}^{2}_{X}$ (d'ailleurs, $d\xi=\frac{-3\omega_{3}}{4}$ !). On vérifie, alors, par un calcul simple, que les formes
$$\xi\wedge dz=z\omega_{3},\,\xi\wedge dx=\frac{5x}{4}\omega_{3},\,\xi\wedge dy=y\omega_{3}, \,\xi\wedge dt=t\omega_{3}$$
définissent toutes des sections du faisceau ${\mathcal L}^{3}_{X}$ et que, par conséquent, $\xi$ définit une section non triviale du faisceau ${\mathcal H}om({\Omega}^{1}_{X}, {\mathcal L}^{3}_{X})$.\vspace{1mm}

\noindent Comme la singularité est de codimension $3$, on a 
$$\Omega^{1}_{X}={\mathcal L}^{1}_{X}=\omega^{1}_{X}$$
et donc
$${\mathcal H}om({\mathcal L}^{1}_{X}, {\mathcal L}^{3}_{X})\simeq{\mathcal H}om({\Omega}^{1}_{X}, {\mathcal L}^{3}_{X})$$
Cela montre que l'injection
$${\mathcal L}^{2}_{X}\into{\mathcal H}om({\mathcal L}^{1}_{X}, {\mathcal L}^{3}_{X})$$
est stricte. On en déduit aussitôt qu'il est impossible d'avoir l'isomorphisme
$${\mathcal L}^{i}_{X}\simeq {\mathcal H}om_{{\mathcal O}_{X}}({\mathcal L}^{m-i}_{X},
{\mathcal L}^{m}_{X})$$
pour $X$ réduit quelconque. En effet, si tel était le cas, on l'aurait aussi pour le normalisé!
Soit $\nu:\widehat{X}\rightarrow X$ le morphisme de normalisation. Alors,
$$\nu_{*}{\mathcal H}om_{{\mathcal O}_{\widehat{X}}}({\mathcal L}^{m-i}_{\widehat X},
{\mathcal L}^{m}_{\widehat X})\simeq {\mathcal H}om_{{\mathcal O}_{X}}({\mathcal L}^{m-i}_{X},
{\mathcal L}^{m}_{X})$$
implique que
$$\nu_{*}{\mathcal L}^{i}_{\widehat X}\simeq \nu_{*}{\mathcal H}om_{{\mathcal O}_{\widehat{X}}}({\mathcal L}^{m-i}_{\widehat X},
{\mathcal L}^{m}_{\widehat X})$$
et, par conséquent, $\nu$ étant fini, 
$${\mathcal L}^{i}_{\widehat X}\simeq {\mathcal H}om_{{\mathcal O}_{\widehat{X}}}({\mathcal L}^{m-i}_{\widehat X},
{\mathcal L}^{m}_{\widehat X})$$
\vspace{2mm}

\noindent 
Le lecteur pourra vérifier que, de façon analogue, pour l'hypersurface normale de ${\Bbb C}^{4}$ 
\vspace{1mm}

\centerline{ $X := \{(x,y,z,t)\in {\Bbb C}^{4}/f(x,y,z,t)=  x^{4} + y^{4} + z^{4} + t^{5}=0\}$}
\noindent qui n'est pas à singularité rationnelle puisque $3/4 + 1/5 =19/20<1$, la forme 
$$ \xi= {zdx\wedge dy\over{t^{4}}} - {ydx\wedge dz\over{t^{4}}} +{xdy\wedge dz\over{t^{4}}} $$
définit une section de ${\mathcal H}om({\Omega}^{1}_{X}, {\mathcal L}^{3}_{X})$ mais pas de ${\mathcal L}^{2}_{X}$. \vspace{1mm}

\indent 
$\bullet$ Il  est facile de voir que ce n'est pas une section de ${\mathcal L}^{2}_{X}$  puisque l'éclatement de la singularité dans la première carte  $y=\alpha.x$, $z=\beta.x$, $t=\gamma.x$, donne pour préimage stricte $\widetilde{X} := \{(\alpha,\beta,\gamma,x)\in {\Bbb C}^{4}/ 1 + \alpha^{4} + \beta^{4} + x\gamma^{5}=0\}$ et une image réciproque de $\xi$ égale à $\displaystyle{d\alpha\wedge d\beta\over{x\gamma^{4}}}$. Mais cette dernière  n'est pas holomorphe sur  $\widetilde{X}$ car la fonction ${\gamma\over{x}}$ ne s'y prolonge pas analytiquement .\vspace{1mm}

\noindent En revanche, elle définit une section du faisceau   ${\mathcal H}om({\Omega}^{1}_{X}, {\mathcal L}^{3}_{X})$ puisque  les formes méromorphes 
$\xi\wedge dx=-5x\omega_{0}$, 
$\xi\wedge dz = -5z\omega_{0}$, 
$\xi\wedge dy =-5y\omega_{0}$ et  
$\xi\wedge dt = -4t\omega_{0}$, $\omega_{0}$ étant un générateur de $\omega^{3}_{X}$, 
se prolongent analytiquement sur chacune des cartes de l'éclatement (et donc sur la désingularisée...).   En effet, comme les variables jouent des rôles symétriques, on peut se contenter de regarder ${\xi\wedge dx}$ qui, par exemple, donne dans la première carte, la forme ${{dx\wedge d\alpha\wedge d\beta\over{x\gamma^{4}}}}$ dont l'holomorphie ne suscite aucun doute sur  $\widetilde{\Gamma} := \{(\alpha,\beta,\gamma,t)\in {\Bbb C}^{4}/ 1 + \alpha^{4} + \beta^{4} + x\gamma^{5}=0\}$.  Les cartes en $z$ et en $y$ conduisent au même calcul et à la même conclusion.  Dans la dernière carte   $x=\alpha.t$, $y=\beta.t$, $z=\gamma.t$, la préimage stricte est $\widetilde{X} := \{(\alpha,\beta,\gamma,t)\in {\Bbb C}^{4}/ \alpha^{4} + \beta^{4} + \gamma^{4} + t=0\}$ et l'image réciproque de ${\xi\wedge dx}$  est 
${{\alpha^{2}\over{t}}{dt\wedge d\beta\wedge d\gamma}+ {\alpha\beta\over{t}}{d\alpha\wedge dt\wedge d\gamma}+
{\alpha\gamma\over{t}}{d\alpha\wedge d\beta\wedge dt}}$
qui est holomorphe sur $\widetilde{X}$.\vspace{2mm}

\noindent
 \subsection{{\bf{Non stabilité  des faisceaux ${\mathcal H}om({\Omega}^{m-q}_{X}, {\mathcal L}^{m}_{X})$ et ${\mathcal H}om({\mathcal L}^{m-q}_{X}, {\mathcal L}^{m}_{X})$ par image réciproque arbitraire et par différentiation extérieure.}}}\vspace{1mm}

\noindent S'ils étaient stable par image réciproque arbitraire, ils le seraient pour toute désingularisation et, donc, coincideraient nécessairement avec le faisceau ${\mathcal L}^{\bullet}_{X}$.\vspace{1mm}

\noindent Pour la différentiation extérieure, on peut reprendre  les exemples de surface  déjà cités:\vspace{1mm}

\indent $X:=\big\{(x,y,z)\in {\Bbb C}^{3}: xy^{2}=z^{3}\big\}$ et de la forme $\xi:=\displaystyle{\frac{ydx}{z^{2}}}$ qui est une section de  ${\mathcal H}om_{{\mathcal O}_{X}}(\Omega^{1}_{X}, {\mathcal L}^{2}_{X})$ dont la différentielle extérieure $d\xi=-7\displaystyle{\frac{dx\wedge dy}{3z^{2}}}$ n'est pas du tout une section du faisceau ${\mathcal L}^{2}_{X}$.\vspace{1mm}

\noindent de même pour $X:=\lbrace (x,y,z)\in{\Bbb C}^3: x^3 +y^3 +z^4=0\rbrace$ et la forme $\xi:=\displaystyle{\frac{ydx - xdy}{z^3}}$ qui est une section  du faisceau ${\mathcal H}om({\Omega}^{1}_{X}, {\mathcal L}^{2}_{X})$ dont la différentielle extérieure  $d\xi=\displaystyle{\frac{dx\wedge dy}{4z^3}}$ est le générateur de $\omega^{2}_{X}$ !\vspace{1mm}

\indent $X:=\big\{(x,y,z,t)\in {\Bbb C}^{4}: f(x,y,z,t)=x^4 +y^5 +z^5 +t^5=0\big\}$ et la forme
$$\xi:=\frac{tdx\wedge dy}{5z^4} -\frac{ydx\wedge dt}{5z^4} +\frac{xdy\wedge dt}{4z^4}$$
qui est une section du faisceau ${\mathcal H}om({\Omega}^{1}_{X}, {\mathcal L}^{3}_{Z})={\mathcal H}om({\mathcal L}^{1}_{X}, {\mathcal L}^{3}_{X})$ 
dont la différentielle extérieure est  $d\xi=-\displaystyle{\frac{3}{4}\omega_{3}}$ qui n'est pas une section du faisceau   ${\mathcal L}^{3}_{X}$ (car singularité non rationnelle!).

\vfill\eject
\phantomsection\addcontentsline{toc}{part}{Les preuves des \theoremref{T1} et  \theoremref{T'1}.}

%
\svnid{$Id: S01-intro.tex 272 2020-01-21 14:29:32Z kebekus $}
\section{\color{blue}{Preuve du \theoremref{T1}.}}
\subversionInfo

{\bf(i){Le morphisme ${\bf{\pi}}^{*}: \pi^{*}({\mathcal D}_{\omega_{S}}({\mathcal F}^{q}_{S}))\rightarrow{\mathcal D}_{\omega_{X}}({\mathcal F}^{q}_{X})$}.} \vspace{1mm}

\noindent 
On va montrer le résultat pour ${\mathcal F}^{\bullet}_{S}:=\omega^{\bullet}_{S}$ et son analogue sur $X$, le raisonnement est strictement le même pour le faisceau ${\mathcal L}^{\bullet}_{S}$. Pour alléger les notations, on pose $\displaystyle{\widetilde{\omega}^{\bullet}_{S}:={\mathcal D}_{\omega_{S}}(\omega^{r-\bullet}_{S})}:={\mathcal H}om(\omega^{r-\bullet}_{S}, \omega^{r}_{S}))$ et de même pour les analogues sur $X$.
 \vspace{2mm}

\indent
$\star$ Montrons l'existence d'un morphisme ${\mathcal O}_{X}$-linéaire d'image réciproque $\displaystyle{{\bf \pi}^{*}:\pi^{*}(\widetilde{\omega}^{q}_{S}) \rightarrow \widetilde{\omega}^{q}_{X}} $ prolongeant naturellement celle des formes holomorphes.\vspace{1mm}

\noindent
Comme le problème est de nature locale sur $X$ et $S$ et que le morphisme $\pi:X\rightarrow S$ est ouvert à fibres de dimension constante $n$, on peut se donner une paramétrisation locale au voisinage d'un point $x$ de $X$, que l'on peut écrire par abus de notation, $\xymatrix{ X\ar@/^1pc/[rr]^{\pi}\ar[r]_{f}& S\times U\ar[r]_{p}&S}$. On désigne par $p'\!:S\times U\rightarrow U$  la seconde  projection.  On a, alors,
$\displaystyle{p^{*}(\widetilde{\omega}^{q}_{S})\simeq{\mathcal H}om(p^{*}{\omega}^{r-q}_{S},  p^{*}{\omega}^{r}_{S})}$
et, en  simplifiant les notations, la décomposition  
$$\omega^{k}_{S\times U}=\bigoplus_{i+j=k} p^{*}(\omega^{i}_{S})\otimes_{{\mathcal O}_{S\times U}}p'^{*}(\Omega^{j}_{U})=\bigoplus_{i+j=k}\omega^{i}_{S}\widehat\otimes_{\Bbb C}\Omega^{j}_{U}$$
donnant, en particulier, l'isomorphisme
$$\omega^{r+n}_{S\times U}\simeq p^{*}(\omega^{r}_{S})\otimes_{{\mathcal O}_{S\times U}}p'^{*}(\Omega^{n}_{U})$$ 
D'autre part, comme ${\mathcal O}_{U}\simeq \Omega^{n}_{U}$, on a les isomorphismes naturels
$$p^{*}{\omega}^{q}_{S}\simeq{\mathcal H}om(p^{*}{\omega}^{r-q}_{S},  p^{*}{\omega}^{r}_{S})\simeq{\mathcal H}om(p^{*}({\omega}^{r-q}_{S})\otimes p'^{*}(\Omega^{n}_{U}) ,  p^{*}({\omega}^{r}_{S})\otimes p'^{*}(\Omega^{n}_{U}))$$
Il suffit alors d'utiliser le morphisme de projection 
$$\omega^{n+r-q}_{S\times U}\rightarrow p^{*}(\omega^{r-q}_{S})\otimes p'^{*}(\Omega^{n}_{U})$$ 
pour obtenir le morphisme naturel
$$p^{*}(\widetilde{\omega}^{q}_{S})\rightarrow{\mathcal H}om(\omega^{n+r-q}_{S\times U}, \omega^{r+n}_{S\times U})=\widetilde{\omega}^{q}_{S\times U} $$
On peut dès lors nous concentrer essentiellement sur le cas  d'un morphisme fini $\pi:X\rightarrow S$ d'espaces complexes de dimension $m$. \vspace{1mm}

\noindent Considérons, alors, une section $\xi$ du faisceau $\widetilde{\omega}^{q}_{S}$ et montrons que son image réciproque $\pi^{*}(\xi)$
qui est  une $q$-forme holomorphe sur $X\setminus\pi^{-1}({\rm Sing}(S))$ et méromorphe sur $X$ définit une section du faisceau  $\widetilde{\omega}^{q}_{X}$ c'est-à-dire qu'elle satisfait la propriété 
$$\pi^{*}(\xi)\wedge \alpha \in \omega^{m}_{X},\,\,\forall\,\alpha\in \omega^{m-q}_{X}$$
Mais, par  définition de $\omega^{m}_{X}$, cela revient à montrer que ce cup-produit vérifie {\emph{la propriété de la trace}}.  Il nous suffit donc de prouver que, pour une paramétrisation locale donnée  $h:X\rightarrow V$ (sur un ouvert de Stein d'un espace numérique), la trace de $ \pi^{*}(\xi)\wedge \alpha$ se prolonge analytiquement sur $V$. Comme on a le choix pour $h$ et que $\pi$ est finie,  on  la choisit de sorte à ce qu'elle se factorise au travers de $S$ en le diagramme $\xymatrix{ X\ar@/^1pc/[rr]^{h}\ar[r]_{\pi}& S\ar[r]_{\phi}&V}$ où $\phi$ est une paramétrisation locale de $S$ en  $\pi(x)$. 
On a, alors, grâce à la compatibilité de la trace avec la composition des morphismes finis, la relation 
$${\rm Tr}_{h}(\pi^{*}(\xi)\wedge \alpha)={\rm Tr}_{\phi}({\rm Tr}_{\pi}(\pi^{*}(\xi)\wedge \alpha))={\rm deg}(\pi).{\rm Tr}_{\phi}(\xi\wedge {\rm Tr}_{\pi}(\alpha))$$
 Comme $\xi$ est une section du faisceau $\widetilde{\omega}^{q}_{S}$,  le cup-produit de cette section méromorphe par toute $(m-q)$-forme méromorphe section du faisceau $\omega^{m-q}_{S}$ définit une section de faisceau $\omega^{m}_{S}$.  Comme  ${\rm Tr}_{\pi}(\alpha)$ définit une section du faisceau $\omega^{m-q}_{S}$ (puisqu'image directe d'un courant $\bar\partial$-fermé sans torsion), on a, donc, en particulier, que $\xi\wedge {\rm Tr}_{\pi}(\alpha))$ est une section de $\omega^{m}_{S}$ et, par conséquent, le prolongement analytique de la forme ${\rm Tr}_{\phi}(\xi\wedge {\rm Tr}_{\pi}(\alpha))$. D'où la conclusion.
 \vspace{1mm}

 \noindent
 On remarque qu'il suffisait de vérifier que ${\rm Tr}_{\pi}(\pi^{*}(\xi)\wedge \alpha))$ définit une section du faisceau $\omega^{m}_{S}$ et ce, pour toute section $\alpha$ du faisceau $\omega^{m-q}_{X}$ (ou tout autre faisceau stable  par trace).\vspace{2mm}

\noindent 
$\bullet$ {\bf{La compatibilité avec la différentielle extérieure}} au sens des courants résulte, en partie, de l'existence, pour tout entier $j\geq q$, de morphismes
$$\widetilde{\omega}^{q}_{X}\rightarrow {\mathcal H}om(\omega^{j-q}_{X}, \omega^{j}_{X})$$
dont la construction découle simplement de la définition de $\omega^{j}_{X}$ en écrivant
$${\mathcal H}om(\omega^{j-q}_{X}, \omega^{j}_{X})\simeq {\mathcal H}om(\omega^{j-q}_{X}, {\mathcal H}om(\Omega^{m-j}_{X}, \omega^{m}_{X}))\simeq {\mathcal H}om(\omega^{j-q}_{X}\otimes\Omega^{m-j}_{X}, \omega^{m}_{X}) $$
et utilisant le cup produit 
$$\omega^{j-q}_{X}\otimes\Omega^{m-j}_{X}\rightarrow \omega^{m-q}_{X}$$
nous donnant finalement le morphisme
$$\widetilde{\omega}^{q}_{X}:={\mathcal H}om(\omega^{m-q}_{X}, \omega^{m}_{X})\rightarrow{\mathcal H}om(\omega^{j-q}_{X}, \omega^{j}_{X}) $$
On en déduit, alors, pour $j=m-1$, le cup produit
$$\widetilde{\omega}^{q}_{X}\otimes{\omega}^{m-1-q}_{X}\rightarrow {\omega}^{m-1}_{X}$$
qui permet de garantir que, pour toute section  $\xi$ de $\widetilde{\omega}^{q}_{X}$ et toute section $\alpha$ de $\omega^{m-q-1}_{X}$, le cup produit $\xi\wedge \alpha$ définit bien une section du faisceau  ${\omega}^{m-1}_{X}$. Ainsi, la relation 
$$ d(\xi\wedge \alpha)=d\xi\wedge \alpha \pm\xi\wedge d\alpha$$
montre clairement que, pour tout $\alpha$, la forme $d\xi\wedge \alpha$ est une section du faisceau ${\omega}^{m}_{X}$; ce qui signifie exactement que  $d\xi$ définit une section de $\widetilde{\omega}^{q+1}_{X}$.\vspace{1mm}

\noindent 
$\bullet$ {\bf{La compatibilité avec la composition des morphismes}} de cette nature relève de la compatibilité de la trace avec la composée des morphismes finis ou du fait que les faisceaux soient sans torsion.  Pour être un peu plus clair, supposons donné  le diagramme $$\xymatrix{X\ar[rr]^{f}\ar[rd]_{h}&&Y\ar[ld]^{g}\\
&S&}$$
avec des  morphismes surjectifs à fibres de dimension constante. Alors, 
 il nous faut montrer que la composée des morphismes
 $$\xymatrix{f^{*}(g^{*}\widetilde{\omega}^{q}_{S})\ar[r]\ar[rd]&f^{*}\widetilde{\omega}^{q}_{S}\ar[d]\\
 &\widetilde{\omega}^{q}_{X}}$$
 coincide avec le morphisme $$h^{*}\widetilde{\omega}^{q}_{S}\rightarrow \widetilde{\omega}^{q}_{X}$$
 En d'autres termes
 $$h^{*}(\xi)=f^{*}(g^{*}\xi),\,\,\forall\,\xi\,{\rm section\,de}\,\widetilde{\omega}^{q}_{S}$$
 On se ramène, à l'aide de factorisation locale,  à la situation suivante: 
 $$\xymatrix{X\ar[rr]^{f'}\ar[rd]_{h'}&&Y\times V\ar[ld]^{g'}\\
&S\times U\times V&}$$
avec $f'$, $g'$ et $h'$ finis, ouverts  et surjectifs. \vspace{1mm}

\noindent La formule de projection, l'image directe au sens des courants et la formule de composition des traces nous donnent, pour tout section $\alpha$ du faisceau $\omega^{m-q}_{X}$,
$${\rm Tr}_{h'}(h'^{*}(\xi)\wedge \alpha)={\rm deg}(h').\xi\wedge{\rm Tr}_{h'}(\alpha)$$
$${\rm Tr}_{g'}({\rm Tr}_{f'}(f'^{*}(\xi)\wedge \alpha)) ={\rm deg}(f'). {\rm Tr}_{g'}(g'^{*}(\xi)\wedge {\rm Tr}_{f'}(\alpha))={\rm deg}(f').{\rm deg}(g').\xi\wedge\underbrace{{\rm Tr}_{g'}(\underbrace{{\rm Tr}_{f'}(\alpha)}_{\in \omega^{m-q}_{Y\times V}}}_{\in \omega^{m-q}_{S\times \times U\times V}})$$
avec $${\rm deg}(h')={\rm deg}(f').{\rm deg}(g')$$
qui, par définition de ces faisceaux, garantit la compatibilité des formules.\vspace{1mm}

\noindent On peut aussi utiliser l'absence de torsion dans ces faisceaux pour établir cette invariance. En effet, on peut mettre en évidence des ouverts denses $U$, $V$ et $W$ de $X$, $Y$ et $S$ respectivement sur lesquels les morphismes $f$, $g$ et $h$ soient lisses et pour lesquels  $h^{*}(\xi\vert_{W})$ et $f^{*}(g^{*}(\xi\vert_{W}))$ coincident puisque ce n'est que l'image réciproque des formes holomorphes. Comme $\widetilde{\omega}^{q}_{X}$ est sans torsion, cette égalité générique est satisfaite globalement.\vspace{2mm}

\indent 
Le même raisonnement vaut pour le faisceau $\widetilde{\mathcal L}^{q}_{X}:={\mathcal H}om({\mathcal L}^{m-q}_{X}, \omega^{m}_{X})$ sachant que ${\mathcal L}^{\bullet}_{X}$ est muni de morphismes trace relativement aux projections finies.  De même, on met en évidence un morphisme naturel
$$\widetilde{\mathcal L}^{q}_{X}\otimes {\mathcal L}^{m-1-q}_{X}\rightarrow \omega^{m-1}_{X}$$
permettant de vérifier que si $\xi\in \widetilde{\mathcal L}^{q}_{X}$ alors $d\xi \in \widetilde{\mathcal L}^{q+1}_{X}$.\vspace{1mm}

\noindent La compatibilité avec la composition des morphismes se traite de façon identique.\vspace{1mm}

\noindent A noter, au passage, que les sous  faisceaux ${\mathcal H}om({\mathcal L}^{m-q}_{X}, {\mathcal L}^{m}_{X})$ et  ${\mathcal H}om(\widetilde{\Omega}^{m-q}_{X}, {\mathcal L}^{m}_{X})$ de $\omega^{q}_{X}$ ne sont pas stables par différentiation extérieure en général (cf exemples).
\vspace{2mm}

\noindent
{\bf(ii) Les faisceaux
${\mathcal D}_{{\mathcal L}_{S}}({\mathcal L}_{S}^{r-\bullet}):={\mathcal H}om({\mathcal L}_{S}^{r-\bullet}, {\mathcal L}_{S}^{r})$.\vspace{2mm}

\indent 
$\star$ Le  pull back $\pi^{*}{\mathcal H}om({\mathcal L}^{r-q}_{S}, {\mathcal L}^{r}_{S})\rightarrow{\mathcal H}om({\mathcal L}^{m-q}_{X}, {\mathcal L}^{m}_{X})$.}\vspace{2mm}

\noindent 
Il nous arrivera de noter $\overline{\mathcal L}^{\bullet}_{Y}:= {\mathcal H}om({\mathcal L}^{s-\bullet}_{Y}, {\mathcal L}^{s}_{Y})$ pour un espace complexe réduit $Y$ de dimension pure $s$. Dans toute la suite, on conservera les notations de {\bf(i)}) \vspace{1mm}

\noindent Pour réduire le problème au cas d'un morphisme fini comme dans {\bf(i)}, il faut s'assurer de l'existence d'une image réciproque pour une projection locale $\displaystyle{p:S\times U\rightarrow S}$. Pour cela, l'interprétation en termes de courants permet de voir que l'on a encore des décompositions analogues en remplaçant $\omega^{\bullet}_{S\times U}$ (resp. $\omega^{\bullet}_{S}$) par ${\mathcal L}^{\bullet}_{S\times U}$ (resp. ${\mathcal L}^{\bullet}_{S}$) et, donc, reprendre les mêmes arguments adaptés à ces faisceaux.\vspace{1mm}

\noindent 
On peut noter que l'existence de ce pull-back peut-être établie en considérant une désingularisation $\phi:\widetilde{S}\rightarrow S$ de $S$ et le diagramme de changement de base
$$\xymatrix{\widetilde{S}\times U\ar[r]^{\tilde\phi}\ar[d]_{\tilde{p}}&S\times U\ar[d]^{p}\\
\widetilde{S}\ar[r]_{\phi}&S}$$

Si $\xi$ est une section du faisceau $\overline{\mathcal L}^{k}_{S}$, on a, par définition,
$$\xi\wedge \alpha \in {\mathcal L}^{r}_{S},\,\,\forall\,\alpha\in {\mathcal L}^{r-k}_{S}$$
Alors,  la relation
$$\tilde{p}^{*}(\phi^{*}(\xi\wedge \alpha))=
\tilde{\phi}^{*}({p}^{*}(\xi\wedge \alpha)),$$
montre que la forme ${p}^{*}(\xi\wedge \alpha)$ définit une section du faisceau ${\mathcal L}^{r}_{S\times U}$.\vspace{1mm}

\noindent D'autre part, pour tester le cup produit ${p}^{*}(\xi) \wedge \beta$ pour une section arbitraire $\beta$ du faisceau  ${\mathcal L}^{m-k}_{S\times U}$ (avec $m:=n+r$), il nous suffit, pour des raisons de dimension, de l'évaluer sur celles qui s'écrivent
$\beta:=p^{*}(\alpha) \wedge p'^{*}(\alpha')$ avec $\alpha$ (resp. $\alpha'$) section de ${\mathcal L}^{r-k}_{S}$ (resp. de $\Omega^{n}_{U})$. On a, alors,
$${p}^{*}(\xi) \wedge \beta={p}^{*}(\xi) \wedge p^{*}(\alpha) \wedge p'^{*}(\alpha')$$
impliquant que ${p}^{*}(\xi)$ définit naturellement une section du faisceau  $\overline{{\mathcal L}}^{k}_{S}$.
\vspace{2mm}

\noindent

On suppose donc $\pi$ fini, surjectif et ouvert (avec $m=r$). Comme ${\mathcal H}om({\mathcal L}^{m-q}_{S}, {\mathcal L}^{m}_{S})$ est un sous faisceau cohérent 
de ${\mathcal H}om({\mathcal L}^{m-q}_{S}, {\omega}^{m}_{S})$ qui est stable par image réciproque d'après {\bf(i)}, toute section $\xi$ du premier faisceau a pour image réciproque $\pi^{*}(\xi)$ une section du faisceau  ${\mathcal H}om({\mathcal L}^{m-q}_{X}, {\omega}^{m}_{X})$ vérifiant
$$\pi^{*}(\xi)\wedge \alpha \in \omega^{m}_{X},\,\,\forall\,\alpha\in {\mathcal L}^{m-q}_{X}$$
On veut montrer que $\pi^{*}(\xi)\wedge \alpha$ est, en fait, une section du faisceau ${\mathcal L}^{m}_{X}$ c'est-à-dire que pour une résolution des singularités de $X$, elle se relève en une forme holomorphe. Cela revient à montrer l'existence de la flèche en pointillé dans le diagramme
$$\xymatrix{{\mathcal H}om({\mathcal L}^{m-q}_{S}, {\mathcal L}^{m}_{S})\ar@{^{(}->}[r] \ar@{.>}[d]&{\mathcal H}om({\mathcal L}^{m-q}_{S}, {\omega}^{m}_{S})\ar[d]^{f^{*}}\\
{\mathcal H}om({\mathcal L}^{m-q}_{X}, {\mathcal L}^{m}_{X})\ar@{^{(}->}[r]&{\mathcal H}om({\mathcal L}^{m-q}_{X}, {\omega}^{m}_{X})}$$
Comme les faisceaux ${\mathcal H}om({\mathcal L}^{m-q}_{S}, {\mathcal L}^{m}_{S})$ (resp. ${\mathcal H}om({\mathcal L}^{m-q}_{X}, {\mathcal L}^{m}_{X})$) sont stables par modifications finies, on peut supposer $X$ et $S$ normaux. En effet, si $\nu:\hat{S}\rightarrow S$ est la normalisation de $S$, on a (déjà vu) l'isomorphisme
$$\nu_{*}{\mathcal H}om({\mathcal L}^{m-q}_{\hat{S}}, {\mathcal L}^{m}_{\hat{S}})\simeq {\mathcal H}om({\mathcal L}^{m-q}_{S}, {\mathcal L}^{m}_{S})$$
de même sur $X$. Il suffit, alors, de travailler sur le diagramme commutatif
$$\xymatrix{\widehat{X}\ar@/_2pc/[dd]_{\phi}\ar[rd]^{\psi}\ar[d]_{\theta}&\\
X\times_{S} \hat{S}\ar[d]_{\hat{\pi}}\ar[r]_{\hat{\nu}}&X\ar[d]^{\pi}\\
\widehat{S}\ar[r]_{\nu}&S}$$
où $\theta$ est la normalisation de l'espace réduit sous jacent au produit fibré. 
La formule de projection et la finitude des morphismes permet de voir que si l'assertion est vraie pour $(\widehat{X}, \widehat{S},\phi)$ elle l'est pour $({X}, {S}, \pi)$.\vspace{1mm}

\noindent
Soit $\phi:\widetilde{X}\rightarrow X$ une désingularisation de $X$. Notons $\psi:=\pi\circ \phi$. Alors, comme  
$$ \phi^{*}(\pi^{*}(\xi)\wedge \alpha)=\psi^{*}(\xi)\wedge \phi^{*}\alpha$$
on a, en vertu de la formule de projection, 
$$\psi_{*}( \phi^{*}(\pi^{*}(\xi)\wedge \alpha))=\xi\wedge \psi_{*}\phi^{*}\alpha=\xi\wedge \pi_{*}\alpha$$
car $\phi_{*}\circ \phi^{*}={\rm Id}$ puisque $X$ est normal. Comme $\pi_{*}\alpha$ est donnée par la trace ou l'image directe au sens des courants, elle définit une section du faisceau ${\mathcal L}^{m-q}_{S}$. Mais $\xi$ étant, par hypothèse,  une section du faisceau ${\mathcal H}om({\mathcal L}^{m-q}_{S}, {\mathcal L}^{m}_{S})$, il s'en suit que $\xi\wedge \pi_{*}\alpha$ définit une section du faisceau ${\mathcal L}^{m}_{S}$. Par conséquent, la forme  $\phi^{*}(\pi^{*}(\xi)\wedge \alpha)$ est  holomorphe (se prolonge holomorphiquement) sur $\widetilde{X}$ imposant à $\pi^{*}(\xi)\wedge \alpha$ d'être une section du faisceau ${\mathcal L}^{m}_{X}$. D'où la conclusion.
\vspace{1mm}

\indent
$\bullet$ {\bf{La compatibilité avec la composition des morphismes}} de cette nature résulte de  l'absence de torsion dans ces faisceaux. En effet, on peut mettre en évidence des ouverts denses $U$, $V$ et $W$ de $X$, $Y$ et $S$ respectivement sur lesquels les morphismes $f$, $g$ et $h$ soient lisses et pour lesquels  $h^{*}(\xi\vert_{W})$ et $f^{*}(g^{*}(\xi\vert_{W}))$ coincident puisque ce n'est que l'image réciproque des formes holomorphes. Comme $\overline{\mathcal L}^{q}_{X}$ est sans torsion, cette égalité générique est satisfaite globalement$\,\blacksquare$
\vspace{1mm}

\noindent

\subsection{\color{blue}{Preuve du \corollaryref{Cor0'}}} \vspace{1mm}

\noindent
Soit $\pi:X\rightarrow S$ un morphisme ouvert, à fibres de dimension constante entre les espaces complexes réduits $X$ et $S$  de dimension pure $m$ et $r$ respectivement. Montrons que l'existence  d'un morphisme de faisceaux 
$\displaystyle{{\bf{\pi}}^{*}\!\!:{\omega}^{r}_{S}\rightarrow \pi_{*}{\omega}^{r}_{X}}$  prolongeant l'image réciproque des formes holomorphes usuelles avec la condition ${\mathcal O}_{S}\simeq {\mathcal H}om(\omega^{r}_{S}, \omega^{r}_{S})$ impose à $\pi$ d'être analytiquement  géométriquement plat. On en déduira, qu'en particulier, tout morphisme ouvert à fibres de dimension constante $n$ sur une base de Cohen Macaulay et ayant cette propriété du pull-back est analytiquement géométriquement plat.\vspace{1mm}

\noindent 
Pour s'assurer de  la platitude géométrique analytique, on utilise la caractérisation donnée dans le (\cite{K1}, {\bf{thm 1}}, p.40) en  montrant qu'à toute factorisation locale
$$\xymatrix{X\ar[rr]^{f}\ar[rd]_{\pi}&&S\times U\ar[ld]^{p}\\
&S&}$$
dans la quelle $f$ est une paramétrisation locale (en un point $x_0$ de $X$) et $p$ la projection canonique,  est associé un unique morphisme de type trace 
$f_{*}\Omega^{n}_{X/S}\rightarrow \Omega^{n}_{S\times U/S}$
prolongeant le morphisme trace usuel
$f_{*}f^{*}(\Omega^{n}_{S\times U/S}\rightarrow \Omega^{n}_{S\times U/S}$ défini génériquement pour $\pi$ et où $n:=m-r$.\vspace{1mm}

\noindent
On suppose  $S$ non normal car, sinon,  $\pi$ est  analytiquement géométriquement plat. Commençons par \vspace{1mm}

\noindent 
On note $r$ (resp. $n+r$) la dimension de $S$ (resp. $X$). Alors, on déduit du morphisme d'image réciproque ${\bf \pi}^{*}:\omega^{r}_{S}\rightarrow \pi_{*}\omega^{r}_{X}$ 
la factorisation naturelle $$\xymatrix{\Omega^{n}_{X/S}\otimes\pi^{*}(\omega^{r}_{S})\ar[r]\ar[rd]&\Omega^{n}_{X/S}\otimes\omega^{r}_{X}\ar[d]\\
&\omega^{n+r}_{X}}$$
qui, grâce à la formule de projection pour le foncteur $f_*$ et au morphisme trace $f_{*}\omega^{n+r}_{X}\rightarrow \omega^{n+r}_{S\times U}$ , donne le diagramme commutatif
 $$\xymatrix{f_{*}(\Omega^{n}_{X/S})\otimes p^{*}(\omega^{r}_{S})\ar[r]\ar[rd]&f_{*}\omega^{n+r}_{X}\ar[d]\\
&\omega^{n+r}_{S\times U}}$$
et, donc, le morphisme
$$f_{*}(\Omega^{n}_{X/S})\rightarrow {\mathcal H}om(p^{*}(\omega^{r}_{S}), \omega^{n+r}_{S\times U})$$
Mais il est facile de voir, qu'en degré maximal, on a  la décomposition  $$\omega^{n+r}_{S\times U}\simeq p^{*}(\omega^{r}_{S})\otimes_{{\mathcal O}_{S\times U}}
p'^{*}(\Omega^{n}_{U})$$
avec $p':S\times U\rightarrow U$ la projection usuelle. Or cela s'écrit aussi
$$\omega^{n+r}_{S\times U}\simeq p^{*}(\omega^{r}_{S})\otimes_{{\mathcal O}_{S\times U}}
\Omega^{n}_{S\times U/S}$$
Ainsi, 
$${\mathcal H}om(p^{*}(\omega^{r}_{S}), \omega^{n+r}_{S\times U})\simeq {\mathcal H}om(p^{*}(\omega^{r}_{S}), p^{*}(\omega^{r}_{S}))\otimes_{{\mathcal O}_{S\times U}}
\Omega^{n}_{S\times U/S}) $$
Alors, sous l'hypothèse ${\mathcal O}_{S}\simeq {\mathcal H}om(\omega^{r}_{S}, \omega^{r}_{S})$, on en déduit le morphisme
$$f_{*}(\Omega^{n}_{X/S})\rightarrow\Omega^{n}_{S\times U/S}$$
qui, par construction, est un morphisme de type trace. Vu le choix arbitraire de la factorisation, on en déduit la platitude géométrique analytique de $\pi$.\vspace{1mm}

\noindent
$\bullet$ Si $S$ est  de Cohen Macaulay, cette hypothèse est encore satisfaite. En effet,  si ${\mathcal D}^{\bullet}_{S}$ est le complexe dualisant de $S$, on a, de façon générale,
$${\mathcal O}_{S}\simeq {\rm I}\!{\rm R}{\mathcal H}om({\mathcal D}^{\bullet}_{S}, {\mathcal D}^{\bullet}_{S})$$
Comme $S$ est de Cohen Macaulay, ${\mathcal D}^{\bullet}_{S}\simeq \omega^{r}_{S}[r]$ et, par suite,
$${\mathcal O}_{S}\simeq {\rm I}\!{\rm R}{\mathcal H}om(\omega^{r}_{S}, \omega^{r}_{S})\simeq {\mathcal H}om(\omega^{r}_{S}, \omega^{r}_{S})$$

\noindent Remarquons, au passage, que si $S$ est de Gorenstein, on a l'isomorphisme
$$\Omega^{n}_{S\times U/S}\simeq {\mathcal H}om(p^{*}(\omega^{r}_{S}), \omega^{n+r}_{S\times U})$$
puisque $\times U$ étant aussi de Gorenstein, on a, pour complexes dualisants,
$${\mathcal D}^{\bullet}_{S}:=\omega^{r}_{S}[r]\simeq {\mathcal O}_{S}[r],\,\,{\mathcal D}^{\bullet}_{S\times U}:=\omega^{r+n}_{S\times U}[n+r]\simeq {\mathcal O}_{S\times U}[n+r]$$
Comme $p$ étant plat, on a 
 $${\mathcal H}om(p^{*}(\omega^{r}_{S}), \omega^{r}_{S\times U})={\rm I}\!{\rm R}{\mathcal H}om({\mathcal D}^{\bullet}_{S}[-r], {\mathcal D}^{\bullet}_{S\times U}[-n-r])=p^{!}({\mathcal O}_{S})[-n]\simeq \Omega^{n}_{S\times U/S} $$
 le dernier isomorphisme n'est rien d'autre que l'isomorphisme de Verdier \cite{V} que l'on construira d'une autre façon dans \cite{K8} et \cite{K9}$\,\blacksquare$
\vspace{1mm}

\noindent
 \section{\color{blue}{Preuve du \theoremref{T'1}.}}
\subversionInfo

On va faire la démonstration pour ${\mathcal F}^{\bullet}_{X}:=\widetilde{\Omega}^{\bullet}_{X}$ (resp. ${\mathcal F}^{\bullet}_{S}:=\widetilde{\Omega}^{\bullet}_{S}$ ). Dans ce cas, selon les notations standards $\widetilde{\mathcal F}^{\bullet}_{X}:=\omega^{\bullet}_{X}$  (resp. $\widetilde{\mathcal F}^{\bullet}_{S}:=\omega^{\bullet}_{S}$ )\vspace{1mm}

\noindent
{\bf(i) Cas propre.}\vspace{1mm}

\noindent
Soit $\pi:X\rightarrow S$ un morphisme propre $n$-géométriquement plat d'espaces complexes réduits de dimension pure $m$ et $r$ respectivement. On veut montrer que, pour tout entier $q$,  il existe des morphismes ${\mathcal O}_{S}$-linéaires de faisceaux cohérents
$${\mathcal T}^{\pi}_{q}:{\rm I}\!{\rm R}^{n}\pi_{*}\omega^{n+q}_{X}\rightarrow\omega^{q}_{S} $$
compatible aux restrictions ouvertes sur $X$ et aux changements de base plats.\vspace{1mm}

\noindent
Pour cela, on  va utiliser la dualité analytique pour un morphisme propre (\cite{RRV71}) pour $q=r$.  Le cas général s'en déduira sans difficultés.\vspace{1mm}

\noindent
 Dans la suite, on désigne par ${\mathcal D}^{\bullet}_{X}$ (resp. ${\mathcal D}^{\bullet}_{S}$) le complexe dualisant de $X$ (resp. $S$) dont la cohomologie ou l'homologie est concentrée sur $[-m, 0]$ (resp.  $[-r, 0]$). On commence par remarquer qu'un tel morphisme, dans la catégorie dérivée, détermine (et est déterminé par)  l'unique morphisme $\xymatrix{{\mathfrak T}_{r}:{\rm I}\!{\rm R}\pi_{*}\omega^{m}_{X}[n]\ar[r]&\omega^{r}_{S}} $. \vspace{1mm}

\noindent Alors, comme  par définition (ou construction),  $\omega^{m}_{X}={\mathcal H}^{-m}({\mathcal D}^{\bullet}_{X})$ (resp.$\omega^{r}_{S}={\mathcal H}^{-r}({\mathcal D}^{\bullet}_{S})$) , on a un morphisme canonique  $\xymatrix{\omega^{m}_{X}\ar[r]&{\mathcal D}^{\bullet}_{X}[-m]}$ et, donc, un morphisme naturel
 $\xymatrix{{\rm I}\!{\rm R}\pi_{*}\omega^{m}_{X}[n]\ar[r]&{\rm I}\!{\rm R}\pi_{*}{\mathcal D}^{\bullet}_{X}[-r]}$ que l'on compose avec la trace $\xymatrix{{\rm I}\!{\rm R}\pi_{*}{\mathcal D}^{\bullet}_{X}={\rm I}\!{\rm R}\pi_{*}\pi^{!}{\mathcal D}^{\bullet}_{S}\ar[r]&{\mathcal D}^{\bullet}_{S}}$ pour avoir la flèche
$\xymatrix{{\mathfrak T}^{r}_{\pi}:{\rm I}\!{\rm R}\pi_{*}\omega^{m}_{X}[n]\ar[r]& {\mathcal D}^{\bullet}_{S}[-r]}$ dont  la cohomologie de degré $0$,
 nous donne finalement le morphisme
$${\mathcal T}^{r}_{\pi}:{\rm I}\!{\rm R}^{n}\pi_{*}\omega^{m}_{X}\rightarrow \omega^{r}_{S}$$
On en déduit aisément les morphismes  ${\mathcal T}^{q}_{\pi}$. En effet, partant de l'isomorphisme $\omega^{n+q}_{X}\simeq{\mathcal H}om(\Omega^{r-q}_{X}, \omega^{m}_{X})$, utilisant l'image réciproque naturelle  $\xymatrix{\Omega^{r-q}_{S}\ar[r]&{\rm I}\!{\rm R}\pi_{*} \Omega^{r-q}_{X}}$ et la dualité relative pour le morphisme $\pi$,  on construit le diagramme 
$$\xymatrix{{\rm I}\!{\rm R}\pi_{*}\omega^{m}_{X}[n]\ar[d]\ar[r]&{\rm I}\!{\rm R}\pi_{*}{\rm I}\!{\rm R}{\mathcal H}om(\Omega^{r-q}_{X}, \omega^{m}_{X}[m])[-r]\ar[r]&{\rm I}\!{\rm R}\pi_{*}{\rm I}\!{\rm R}{\mathcal H}om(\Omega^{r-q}_{X}, {\mathcal D}^{\bullet}_{X}[-r])\eq[d]\\
{\rm I}\!{\rm R}{\mathcal H}om(\Omega^{r-q}_{S}, {\mathcal D}^{\bullet}_{S}[-r])&&{\rm I}\!{\rm R}{\mathcal H}om({\rm I}\!{\rm R}\pi_{*}\Omega^{r-q}_{X}, {\mathcal D}^{\bullet}_{S}[-r])\ar[ll]}$$
Alors, il suffit simplement de remarquer que l'on a 
$${\mathcal H}^{0}({\rm I}\!{\rm R}{\mathcal H}om(\Omega^{r-q}_{S}, {\mathcal D}^{\bullet}_{S}[-r]))\simeq\omega^{q}_{S} $$
dont on s'en convainc aisément en prenant un plongement local $i$ de $S$ dans une variété de Stein $V$ de dimension $N$ et utilisant la dualité analytique pour un plongement
$$i_{*}{\rm I}\!{\rm R}{\mathcal H}om(\Omega^{r-q}_{S}, {\mathcal D}^{\bullet}_{S}))\simeq {\rm I}\!{\rm R}{\mathcal H}om(i_{*}\Omega^{r-q}_{S}, \Omega^{N}_{V}[N])$$
qui donne en prenant la cohomologie d'ordre $-r$
$${\mathcal E}xt^{-r}(\Omega^{r-q}_{S}, {\mathcal D}^{\bullet}_{S})\simeq {\mathcal E}xt^{N-r}(i_{*}\Omega^{r-q}_{S}, \Omega^{N}_{V})\simeq{\mathcal H}om(i_{*}\Omega^{r-q}_{S}, {\mathcal E}xt^{N-r}( i_{*}{\mathcal O}_{S}, \Omega^{N}_{V})\simeq i_{*}\omega^{q}_{S} $$
sachant que $\displaystyle{\omega^{q}_{S}\simeq{\mathcal H}om(\Omega^{r-q}_{S}, \omega^{r}_{S})}$, $i_{*}\omega^{r}_{S}={\mathcal Ext}^{N-r}(i_{*}{\mathcal O}_{S}, {\Omega}^{N}_{Y})$ et que ${\mathcal Ext}^{j}(i_{*}{\mathcal O}_{S}, {\Omega}^{N}_{Y})=0$ pour tout $j<N-r$.\vspace{1mm}

\noindent
Ainsi, en prenant la cohomologie d'ordre $0$ de la première flèche verticale, on obtient le morphisme
$${\mathcal T}^{q}_{\pi}:{\rm I}\!{\rm R}^{n}\pi_{*}\omega^{n+q}_{X}\rightarrow \omega^{q}_{S}$$
\vspace{1mm}

\noindent
Pour les faisceaux $\overline{\Omega}^{r-q}_{X}$, $\widehat{\Omega}^{r-q}_{X}$ ou ${\mathcal L}^{r-q}_{X}$, on a besoin des morphismes images réciproques
$$ \widehat{\Omega}^{r-q}_{S}\rightarrow\pi_{*}\widehat{\Omega}^{r-q}_{X},\,\,\overline{\Omega}^{r-q}_{S}\rightarrow\pi_{*}\overline{\Omega}^{r-q}_{X},\,\,{\mathcal L}^{r-q}_{S}\rightarrow\pi_{*}{\mathcal L}^{r-q}_{X}$$
dont on trouvera la construction dans \cite{B5} pour le premier et \cite{K4} pour les derniers.\vspace{1mm}

\noindent
On met ainsi en évidence les morphismes
$${\rm I}\!{\rm R}^{n}\pi_{*}{\mathcal D}_{X}(\Omega^{r-q}_{X})\rightarrow {\mathcal D}_{S}(\Omega^{r-q}_{S}),\,\,{\rm I}\!{\rm R}^{n}\pi_{*}{\mathcal D}_{X}(\overline{\Omega}^{r-q}_{X})\rightarrow {\mathcal D}_{S}(\overline{\Omega}^{r-q}_{S})$$
$${\rm I}\!{\rm R}^{n}\pi_{*}{\mathcal D}_{X}(\widehat{\Omega}^{r-q}_{X})\rightarrow {\mathcal D}_{S}(\widehat{\Omega}^{r-q}_{S}),\,\,{\rm I}\!{\rm R}^{n}\pi_{*}{\mathcal D}_{X}({\mathcal L}^{r-q}_{X})\rightarrow {\mathcal D}_{S}({\mathcal L}^{r-q}_{S})$$
$\bullet$ {\bf la compatibilité avec le changement de base:}
Tous ces morphismes sont compatibles aux restrictions ouvertes sur la base et la source mais seul le morphisme $${\rm I}\!{\rm R}^{n}\pi_{*}{\mathcal D}_{X}({\mathcal L}^{r-q}_{X})\rightarrow {\mathcal D}_{S}({\mathcal L}^{r-q}_{S})$$
est compatible aux changements de base à morphismes dans ${\mathcal M}_{0}$ (i.e ouvert et à fibres de dimension constante $n$). \vspace{1mm}

\noindent On pose $\widetilde{\mathcal L}^{\bullet}_{X}:={\mathcal D}_{X}({\mathcal L}^{r-\bullet}_{X})$ (resp. $\widetilde{\mathcal L}^{\bullet}_{S}:={\mathcal D}_{S}({\mathcal L}^{r-\bullet}_{S})$). Soit
$$\xymatrix{\widehat{X}\ar[d]_{\hat{\pi}}\ar[r]^{\Theta}&X\ar[d]^{\pi}\\
\widehat{S}\ar[r]_{\nu}&S}$$
est un diagramme de changement de base avec $\nu$ un élément de ${\mathcal M}_{0}$ et, donc, $\Theta$ aussi puisque constance de la dimension des fibres et ouverture sont des notions stables par changement de base arbitraire. \vspace{1mm}

\noindent Les morphismes précédemment définis et la première partie permettent de mettre en évidence le diagramme commutatif
$$\xymatrix{\nu^{*}{\rm I}\!{\rm R}^{n}\pi_{*}\widetilde{\mathcal L}^{n+q}_{X}\ar[d]_{\nu^{*}({\mathcal T}^{q}_{\pi})}\ar[r]&{\rm I}\!{\rm R}^{n}{\hat\pi}_{*}\widetilde{\mathcal L}^{n+q}_{\widehat{X}}\ar[d]^{{\mathcal T}^{q}_{\pi}}\\
\nu^{*}\widetilde{\mathcal L}^{q}_{{S}}\ar[r]&\widetilde{\mathcal L}^{q}_{\widehat{S}}}$$
dans lequel seule la composée supérieure mérite d'être expliquée. Pour ce faire, on remarque que l'image réciproque décrite dans la partie {{\bf I}} ci-dessus, $\displaystyle{\widetilde{\mathcal L}^{n+q}_{X}\rightarrow{\Theta}_{*}\widetilde{\mathcal L}^{n+q}_{\widehat{X}}}$ induit le morphisme  $\xymatrix{{\rm I}\!{\rm R}^{n}\pi_{*}\widetilde{\mathcal L}^{n+q}_{X}\ar[r]&{\rm I}\!{\rm R}^{n}\pi_{*}{\Theta}_{*}\widetilde{\mathcal L}^{n+q}_{\widehat{X}}}$. D'autre part, il existe un morphisme canonique $\xymatrix{{\rm I}\!{\rm R}^{n}\pi_{*}{\Theta}_{*}\widetilde{\mathcal L}^{n+q}_{\widehat{X}}\ar[r]&{\rm I}\!{\rm R}^{n}(\pi\circ{\Theta})_{*}\widetilde{\mathcal L}^{n+q}_{\widehat{X}}}$ qui est un morphisme latéral de la suite spectrale
$${\rm I}\!{\rm R}^{i}\pi_{*}{\rm I}\!{\rm R}^{j}{\Theta}_{*}\widetilde{\mathcal L}^{n+q}_{\widehat{X}}
\Longrightarrow{\rm I}\!{\rm R}^{i+j}(\pi\circ{\Theta})_{*}\widetilde{\mathcal L}^{n+q}_{\widehat{X}}$$
La deuxième suite spectrale associée à la composée $\nu\circ\hat{\pi}$ et de même aboutissement fournit un morphisme latéral
$$\xymatrix{{\rm I}\!{\rm R}^{n}(\pi\circ{\Theta})_{*}\widetilde{\mathcal L}^{n+q}_{\widehat{X}}\ar[r]&\nu_{*}{\rm I}\!{\rm R}^{n}{\hat\pi}_{*}\widetilde{\mathcal L}^{n+q}_{\widehat{X}}}$$
A noter que si $\nu$ est plat, la construction du diagramme est simplifiée puisque l'on a un isomorphisme de changement de base $\displaystyle{\nu^{*}{\rm I}\!{\rm R}^{n}{\pi}_{*}\widetilde{\mathcal L}^{n+q}_{X}\simeq{\rm I}\!{\rm R}^{n}{\hat\pi}_{*}(\Theta^{*}\widetilde{\mathcal L}^{n+q}_{X})}$.\vspace{1mm}

\noindent
$\bullet${\bf{La compatibilité avec la composition des morphismes:}}\vspace{1mm}

\noindent
Etant donné un diagramme commutatif d'espaces analytiques complexes
$$\xymatrix{X_{2}\ar[rr]^{\Psi}\ar[rd]_{\pi_{2}}&&X_{1}\ar[ld]^{\pi_{1}}\\
&S&}$$
 avec $\pi_{1}$ (resp. $\pi_{2}$)  propre universellement équidimensionnels de dimension relative $n_{1}$ (resp. $n_{2}$)  et $\Psi$ propre de dimension relative bornée par l'entier   $d:=n_{2}-n_{1}$. \vspace{1mm}
 
 \noindent Le morphisme $\Psi$ étant propre, on a, d'après ce qui précède, un morphisme canonique 
 $${\rm I}\!{\rm
 R}^{d}{\Psi}_{*}\omega^{n_{2}+q}_{X_{2}}\rightarrow \omega^{n_{1}+q}_{X_{1}}$$
 auquel on applique le foncteur ${\rm I}\!{\rm
 R}^{n_{1}}{\pi_{1}}_{*}$ pour avoir 
 $${\rm I}\!{\rm
 R}^{n_{1}}{\pi_{1}}_{*}\big({\rm I}\!{\rm
 R}^{d}{\Psi}_{*}\omega^{n_{2}+q}_{X_{2}})\rightarrow {\rm I}\!{\rm
 R}^{n_{1}}{\pi_{1}}_{*}\omega^{n_{1}+q}_{X_{1}}$$
 qui, grâce à l'isomorphisme de foncteurs ${\rm I}\!{\rm
 R}^{n_{1}}{\pi_{1}}_{*}\circ {\rm I}\!{\rm
 R}^{d}{\Psi}_{*}\simeq {\rm I}\!{\rm
 R}^{n_{2}}{\pi_{2}}_{*}$ dû à l'annulation des images directe supérieures au-delà de la dimension des fibres, nous conduit au diagramme commutatif de faisceaux analytiques cohérents
$$\xymatrix{{\rm I}\!{\rm
 R}^{n_{2}}{\pi_{2}}_{*}\omega^{n_{2}+q}_{X_{2}}\ar[rr]\ar[rd]_{{\mathfrak T}^{q}_{\pi_{2}}}&&
 {\rm I}\!{\rm
 R}^{n_{1}}{\pi_{1}}_{*}\omega^{n_{1}+q}_{X_{1}}\ar[ld]^{{\mathfrak T}^{q}_{\pi_{1}}}\\
&{\omega}^{q}_{S}&}$$
\indent{\bf(ii) Cas général .}\vspace{1mm}

\noindent
On va montrer qu'à tout  morphisme  d'espaces complexes réduits de dimension pure $m$ et $r$ respectivement, $\pi:X\rightarrow S$ dont les fibres sont de dimension $n$, est associé, en tout degré $q$,  un unique  morphisme  ${\mathcal O}_{S}$-linéaires de faisceaux analytiques (en fait quasi-cohérents au sens de la géométrie analytique, comme nous le verrons plus tard)
$${\mathcal T}^{\pi}_{q}:{\rm I}\!{\rm R}^{n}\pi_{!}\omega^{n+q}_{X}\rightarrow\omega^{q}_{S} $$
compatible aux restrictions ouvertes sur $X$ et aux changements de base plats. En conservant les notations précédentes, il nous suffit d'utiliser le morphisme trace de la dualité relative de Ramis-Ruget (\cite{RR74}) donné par ${\rm I}\!{\rm R}\pi_{!}{\mathcal D}^{\bullet}_{X}\rightarrow{\mathcal D}^{\bullet}_{S}$ et le morphisme canonique de \cite{RRV71} 
$${\rm I}\!{\rm R}\pi_{!}{\rm I}\!{\rm R}{\mathcal H}om({\mathcal A}^{\bullet}, {\rm I}\!{\rm R}{\mathcal H}om({\rm I}\!{\rm L}\pi^{*}{\rm I}\!{\rm R}{\mathcal H}om({\mathcal B}^{\bullet}, {\mathcal D}^{\bullet}_{S}), {\mathcal D}^{\bullet}_{X}) )\rightarrow{\rm I}\!{\rm R}{\mathcal H}om({\rm I}\!{\rm R}\pi_{*}{\mathcal A}^{\bullet}, {\mathcal B}^{\bullet})$$
défini pour tout morphisme $\pi:X\rightarrow S$ d'espaces analytiques complexes et tout complexe ${\mathcal A}^{\bullet}$ (resp. ${\mathcal B}^{\bullet}$) de ${\mathcal O}_{Z}$-modules à cohomologie cohérente ( à cohomologie cohérente et bornée à gauche).  
Il en découle, alors,  le diagramme commutatif 
$$\xymatrix{{\rm I}\!{\rm R}\pi_{!}\omega^{n+q}_{X}[n]\ar[r]\ar[rrd]&{\rm I}\!{\rm R}\pi_{!}{\rm I}\!{\rm R}{\mathcal H}om({\Omega}^{r-q}_{X}, {\mathcal D}^{\bullet}_{X}[-m+n])\ar[r]&{\rm I}\!{\rm R}{\mathcal H}om({\rm I}\!{\rm R}\pi_{*}{\Omega}^{r-q}_{X}, {\mathcal D}^{\bullet}_{S}[-r] )\ar[d]\\
&&{\rm I}\!{\rm R}{\mathcal H}om({\Omega}^{r-q}_{S}, {\mathcal D}^{\bullet}_{S}[-r] )}$$ 
En tant que morphisme dans la catégorie dérivée, la flèche diagonale est entièrement déterminée par la flèche obtenue en prenant la cohomologie de degré $0$  en raison des annulations des faisceaux de cohomologie ${\rm I}\!{\rm R}^{j}\pi_{!}\omega^{n+q}_{X}$ sont nuls pour tout $j\geq n$ et ${\mathcal H}^{j}({\rm I}\!{\rm R}{\mathcal H}om({\Omega}^{r-q}_{S}, {\mathcal D}^{\bullet}_{S}[-r] ))$ pour tout entier $j<0$. Ce dernier point est immédiat d'après la dualité de Serre puisque, pour tout ouvert de Stein $U$ dans $S$, la dualité (topologique) absolue donne
$${\rm Ext}^{-r-j}(U; {\Omega}^{r-q}_{S}, {\mathcal D}^{\bullet}_{S})\simeq\big({\rm H}^{r+j}_{c}(U,  {\Omega}^{r-q}_{S})\big)^{'}=0$$
De plus, on a 
$$\omega^{q}_{S}\simeq{\mathcal Ext}^{-r}({\Omega}^{r-q}_{S}, {\mathcal D}^{\bullet}_{S})$$
qui est la cohomologie de degré $0$ du complexe $\displaystyle{{\rm I}\!{\rm R}{\mathcal H}om({\Omega}^{r-q}_{S}, {\mathcal D}^{\bullet}_{S}[-r] )}$. D'où finalement  les morphismes
$${\mathcal T}^{q}_{\pi}:{\rm I}\!{\rm R}^{n}\pi_{!}\omega^{n+q}_{X}\rightarrow\omega^{q}_{S}$$
Les propriétés fonctorielles se prouvent de la même façon que dans le cas propre. La compatibilité avec la composition des morphismes équidimensionnels, résulte de suite spectrales standards et des annulations des images directes supérieures en les degrés supérieures à la dimension des fibres donnant
$${\rm I}\!{\rm R}^{n_{1}}{\pi_{1}}_{!}{\rm I}\!{\rm R}^{d_{1}}{\psi}_{!}\simeq{\rm I}\!{\rm R}^{n_{2}}{\pi_{2}}_{!}$$
De la même façon, on met en évidence les morphismes 
$${\rm I}\!{\rm R}^{n}\pi_{!}{\mathcal D}_{X}(\Omega^{r-q}_{X})\rightarrow {\mathcal D}_{S}(\Omega^{r-q}_{S}),\,\,{\rm I}\!{\rm R}^{n}\pi_{!}{\mathcal D}_{X}(\overline{\Omega}^{r-q}_{X})\rightarrow {\mathcal D}_{S}(\overline{\Omega}^{r-q}_{S})$$
$${\rm I}\!{\rm R}^{n}\pi_{!}{\mathcal D}_{X}(\widehat{\Omega}^{r-q}_{X})\rightarrow {\mathcal D}_{S}(\widehat{\Omega}^{r-q}_{S}),\,\,{\rm I}\!{\rm R}^{n}\pi_{!}{\mathcal D}_{X}({\mathcal L}^{r-q}_{X})\rightarrow {\mathcal D}_{S}({\mathcal L}^{r-q}_{S})$$
tous compatibles aux restrictions ouvertes sur $X$ et $S$. On établit la ${\mathcal M}_{0}$-compatibilité aux changements de base pour la dernière flèche et la compatibilité à la composition des morphismes dans le sens précisé ci-dessus comme dans le cas propre$\,\blacksquare$
\section{\color{blue}{Preuve du \theoremref{T2}.}}\vspace{1mm}
Soit $\pi:X\rightarrow S$ un morphisme géométriquement plat à fibres de dimension pure $n$. Alors, on montre qu' il existe un morphisme de ${\mathcal O}_{S}$-module (cohérents si le morphisme est propre)
$${\mathcal T}^{j}_{\pi}: {\rm I}\!{\rm R}^{n}\pi_{!}\overline{\Omega}^{n+j}_{X}\rightarrow\overline{\Omega}^{j}_{S}$$ 
vérifiant les propriétés annoncées dans le \theoremref{T2}.\vspace{1mm}

\noindent  Signalons que du {\bf {théorème 1.0.4}, p.17} de \cite{B6} et sous ses conditions, on peut déduire (avec les techniques habituelles d'intégration sur les cycles cf \cite{K1}, \cite{K2} ou \cite{K4}) l'existence d'un morphisme
$${\rm H}^{n}_{\Phi}(X, \overline{\Omega}^{n}_{X})\rightarrow \Gamma(S, \overline{\Omega}^{0}_{S})$$
en considérant une famille analytique de cycles paramétrée par un espace complexe réduit dont les supports sont contenus dans une certaine famille de supports $\Psi$ rencontrant la famille paracompactifiante $\Phi$ en des fermés $S$-propres (donc compact pour chaque $s$). Comme nous l'avons expliqué dans \cite{K0} ou \cite{K4}, cela induit un morphisme
$${\mathcal T}^{0}_{\pi}: {\rm I}\!{\rm R}^{n}\pi_{!}\overline{\Omega}^{n}_{X}\rightarrow\overline{\Omega}^{0}_{S}$$ 
doté des propriétés fonctorielles attendues.\vspace{1mm}

\noindent 
Mais l'existence d'une telle flèche impose,  relativement à toute factorisation locale $f:X\rightarrow S\times U$, l'existence d'un morphisme de type trace $f_{*}\overline{\Omega}^{n}_{X}\rightarrow\overline{\Omega}^{n}_{S\times U}$. Réciproquement, ces traces locales permettent la reconstruction du morphisme global
 ${\mathcal T}^{0}_{\pi}$.\vspace{1mm}
 
 \noindent
 Dans ce qui suit, nous allons mettre en évidence ces morphismes trace de façon complètement indépendantes pour aboutir à la construction des morphismes globaux ${\mathcal T}^{j}_{\pi}$.\vspace{2mm}
 
 \noindent
 Pour cela, nous avons besoin des résultats fondamentaux suivants:
\Prop{}{} \label{prop5}Soit $f:X\rightarrow Y$ un morphisme ouvert et fini d'espaces complexes réduits de dimension pure $m$. Alors, il existe un morphisme de type trace ${\mathcal T}^{m}: f_{*}(\overline{\Omega}^{m}_{X})\rightarrow \overline{\Omega}^{m}_{Y}$ naturellement induit par le morphisme trace $\frak{T}^{m}: f_{*}{\omega}^{m}_{X}\rightarrow {\omega}^{m}_{Y}$.\rm

\begin{proof} Considérons une résolution spéciale $\nu:{Y_{1}}\rightarrow Y$ et le diagramme de changement de base associé
$$\xymatrix{X_{1}\ar[r]^{\theta}\ar[d]_{f_{1}}&X\ar[d]^{f}\\
Y_{1}\ar[r]_{\nu}&Y}$$
Comme $\theta$ est une modification propre, on peut compléter ce diagramme grâce à une résolution spéciale $\pi:\widetilde{X}\rightarrow X$ pour avoir le diagramme commutatif
$$\xymatrix{&\widetilde{X}\ar@/_3pc/[ldd]_{\psi}\ar[ld]_{\alpha}\ar[rd]^{\pi}&\\
X_{1}\ar[rr]^{\theta}\ar[d]_{f_{1}}&&X\ar[d]^{f}\\
Y_{1}\ar[rr]_{\nu}&&Y}$$
Comme les faisceaux $\pi^{*}(\widetilde{\Omega}^{m}_{X})/{\mathcal T}_{\pi}$ et $\nu^{*}(\widetilde{\Omega}^{m}_{Y})/{\mathcal T}_{\nu}$ sont localement libre de rang $m$, on a les isomorphismes
$$\pi^{*}(\widetilde{\Omega}^{m}_{X})/{\mathcal T}_{\pi} \simeq {\rm Det}(\pi).\Omega^{m}_{\widetilde{X}},\,\,\,\nu^{*}(\widetilde{\Omega}^{m}_{Y})/{\mathcal T}_{\nu} \simeq {\rm Det}(\nu).\Omega^{m}_{{Y_{1}}} $$
et $$\psi^{*}(\nu^{*}(\widetilde{\Omega}^{m}_{Y})/{\mathcal T}_{\nu})\simeq {\rm Det}(\psi)\Omega^{m}_{\widetilde{X}}$$
dont on déduit 
$$\pi^{*}(\widetilde{\Omega}^{m}_{X})/{\mathcal T}_{\pi}\simeq{\rm Det}(\pi){\rm Det}(\psi)^{-1}\psi^{*}(\nu^{*}(\widetilde{\Omega}^{m}_{Y})/{\mathcal T}_{\nu})$$
où $\Delta:={\rm Det}(\pi){\rm Det}(\psi)^{-1}$ est encore un faisceau cohérent inversible d'idéaux.\vspace{1mm}

\noindent 
En effet, en notant $\Sigma$ (resp. $\Sigma_{1}$) le lieu singulier ou de non locale liberté de $\Omega^{1}_{X}$ (resp. $\Omega^{1}_{Y}$) et par  ${\rm E}:=\pi^{-1}(\Sigma)$ (resp.  ${\rm E}_{1}:=\nu^{-1}(\Sigma_{1})$, l'inclusion $f^{-1}(\Sigma_{1})\subset \Sigma$ induit, naturellement, l'inclusion $\pi^{-1}(f^{-1}(\Sigma_{1}))\subset {\rm E}$. Mais $\pi^{-1}(f^{-1}(\Sigma_{1}))=\alpha^{-1}(\theta^{-1}((f^{-1}(\Sigma_{1}))))$, qui, puisque le carré est cartésien (changement de base induit par $\nu$), s'écrit aussi
$$\pi^{-1}(f^{-1}(\Sigma_{1}))=\alpha^{-1}(f_{1}^{-1}({\rm E}_{1})=\psi^{-1}({\rm E}_{1})$$
Ainsi, $$\psi^{-1}({\rm E}_{1})\subset {\rm E}$$
De plus, seuls nous intéressent les points de  $\psi^{-1}({\rm E}_{1})$ de codimension $1$ car en codimension deux et plus  ${\rm Det}(\psi)$ coincide avec ${\mathcal O}_{\widetilde X}$ par prolongement analytique. On peut dont supposer que $\displaystyle{{\rm dim}(\psi^{-1}({\rm E}_{1}))={\rm dim}({\rm E})}$ et que, par conséquent, $\psi^{-1}({\rm E}_{1})$ est une réunion finie de composantes irréductibles de ${\rm E}$.  Sans enfreindre la généralité, on se ramène, alors, au cas où $\displaystyle{\psi^{-1}({\rm E}_{1})={\rm E}}$ puisque  ${\rm Det}(\pi)$ et ${\rm Det}(\psi)$ sont isomorphes à ${\mathcal O}_{\widetilde{X}}$ en les points de codimension au moins deux dans $\widetilde{X}$.\vspace{1mm}

\noindent 
Par commutativité du diagramme, on a, alors:
$$f_{*}\pi_{*}(\pi^{*}(\widetilde{\Omega}^{m}_{X})/{\mathcal T}_{\pi})\simeq \nu_{*}\psi_{*}(\Delta.\psi^{*}(\nu^{*}(\widetilde{\Omega}^{m}_{Y})/{\mathcal T}_{\nu})\rightarrow\nu_{*}\psi_{*}\psi^{*}(\nu^{*}(\widetilde{\Omega}^{m}_{Y})/{\mathcal T}_{\nu}) $$
Mais $$\psi_{*}\psi^{*}(\nu^{*}(\widetilde{\Omega}^{m}_{Y})/{\mathcal T}_{\nu})\simeq \nu^{*}(\widetilde{\Omega}^{m}_{Y})/{\mathcal T}_{\nu}\otimes_{{\mathcal O}_{Y_1}}\psi_{*}({\mathcal O}_{\widetilde X})$$
On utilise le morphisme trace dans la composition
$$\psi_{*}({\mathcal O}_{\widetilde X})={f_{1}}_{*}\alpha_{*}{\mathcal O}_{\widetilde X}={f_{1}}_{*}{\mathcal L}^{0}_{X_1}\rightarrow {\mathcal O}_{Y_1}$$
A noter que si $X_1$ est normal ou $\alpha $ à fibres connexes, ce que l'on peut d'ailleurs supposer, $\alpha_{*}{\mathcal O}_{\widetilde X}={\mathcal O}_{X_1}$.\vspace{1mm}

\noindent
Ainsi obtient-on le morphisme trace recherché
$${\mathcal T}^{m}: f_{*}(\overline{\Omega}^{m}_{X})\rightarrow \overline{\Omega}^{m}_{Y}$$ 
qui, par construction, est induit par le morphisme trace $f_{*}{\mathcal L}^{m}_{X}\rightarrow {\mathcal L}^{m}_{Y}$ lui même induit par la trace $f_{*}{\omega}^{m}_{X}\rightarrow {\omega}^{m}_{Y}\,\blacksquare$  

\end{proof} 
\Prop{}{} \label{prop6} Soit $X$ un espace complexe réduit de dimension pure $m$. Alors, le morphisme canonique et injectif
$$\overline{\Omega}^{j}_{X}\rightarrow {\mathcal H}om(\overline{\Omega}^{m-j}_{X}, \overline{\Omega}^{m}_{X})$$
est un isomorphisme.\rm
\begin{proof}
   Vu la description de ces faisceaux en termes de somme de formes holomorphes à coefficients continus localement bornés près du lieu singulier (cf \cite{B5}, {\bf Thm 3.0.2}, p. 58), ils sont stables par cup produit interne
    $$\overline{\Omega}^{j}_{X}\otimes_{{\mathcal O}_{X} }\overline{\Omega}^{k}_{X}\rightarrow \overline{\Omega}^{j+k}_{X}$$
    ce qui permet de définir un morphisme naturel de faisceaux cohérents 
    $$\overline{\Omega}^{j}_{X}\rightarrow {\mathcal H}om(\overline{\Omega}^{m-j}_{X}, \overline{\Omega}^{m}_{X})$$
    nécessairement injectif puisque génériquement bijectif et les faisceaux sont  sans torsion.\vspace{1mm}

    \noindent
    Pour définir une flèche inverse, il suffit de remarquer qu'étant donnée une section $\phi$ du faisceau ${\mathcal H}om(\overline{\Omega}^{m-j}_{X}, \overline{\Omega}^{m}_{X})$ sous faisceau du faisceau $\omega^{j}_{X}$, l'hypothèse $\phi(\alpha) \in \overline{\Omega}^{m}_{X}$ pour tout $\alpha\in \overline{\Omega}^{m-j}_{X})$ permet d'écrire
    $$\phi(\alpha)=\sum_{\vert{I}\vert}b_{I} w_{I}$$
    où $\vert{I}\vert=m$ et $w_{I}$ des $m$-formes holomorphes sans torsion. Dans un système de coordonnées génériques $x=(x_1, x_2,\cdots,x_m)$, cette écriture peut-être  réduite à 
    $$\phi(\alpha)= b_{\alpha} dx$$
    On associera, alors, en particulier, à chaque forme holomorphe $\alpha=dx^{J}$ avec $\vert{J}\vert=m-j$, un coefficient $b_J$ et on pose $$\xi:=\sum_{J}b_{J}dx^{J^c}$$
    avec $J^c$ ensemble d'indice complémentaire de $J$. Comme on le constate c'est une somme de produit de fonctions bornées continues par des formes holomorphes. De plus, les morphismes $\phi$ et cup produit par $\xi$ coincident génériquement et, par suite, globalement puisqu'ils définissent tous deux des sections du faisceau $\omega^{j}_{X}$ qui est sans torsion. D'où un morphisme inverse  
   de  ${\mathcal H}om(\overline{\Omega}^{m-j}_{X}, \overline{\Omega}^{m}_{X})$ vers $\overline{\Omega}^{j}_{X})$ associant à tout morphisme $\phi$ la section associée $\xi\,\blacksquare$
\end{proof}
\cor{}{} \label{co1}Soit $f:X\rightarrow Y$ un morphisme ouvert et fini d'espaces complexes réduits de dimension pure $m$. Alors, il existe des morphismes de type trace ${\mathcal T}^{j}: f_{*}(\overline{\Omega}^{j}_{X})\rightarrow \overline{\Omega}^{j}_{Y}$ naturellement induit par les morphisme trace $\frak{T}^{j}: f_{*}{\omega}^{j}_{X}\rightarrow {\omega}^{j}_{Y}$.\rm
\begin{proof}
Il suffit d'écrire
$$\xymatrix{f_{*}\overline{\Omega}^{j}_{X}\eq[r]\ar[dd]&f_{*}{\mathcal H}om(\overline{\Omega}^{m-j}_{X}, \overline{\Omega}^{m}_{X})\ar[r]&{\mathcal H}om(f_{*}\overline{\Omega}^{m-j}_{X}, f_{*}\overline{\Omega}^{m}_{X})\ar[d]\\
&&{\mathcal H}om(f_{*}\overline{\Omega}^{m-j}_{X}, \overline{\Omega}^{m}_{Y})\ar[d]\\
\overline{\Omega}^{j}_{X}\eq[rr]&&{\mathcal H}om(\overline{\Omega}^{m-j}_{Y}, \overline{\Omega}^{m}_{Y})}$$
où les isomorphismes sont donnés par la \propositionref{prop6} et les deux flèches verticales  par la trace de la \propositionref{prop5} et le pull-back $\,\blacksquare$
\end{proof}
\subsection{\color{blue}{Preuve du \theoremref{T2}}.} Soient $X$ et $S$ deux espaces complexes réduits et $\pi:X\rightarrow S$ un morphisme géométriquement plat à fibres de dimension pure $n$. Pour construire un morphisme de ${\mathcal O}_{S}$-module (cohérents si le morphisme est propre)
$${\mathcal T}^{j}_{\pi}: {\rm I}\!{\rm R}^{n}\pi_{!}\overline{\Omega}^{n+j}_{X}\rightarrow\overline{\Omega}^{j}_{S}$$ 
( en remplaçant $\pi_{!}$ par $\pi_{*}$ dans le cas propre) compatible aux restrictions ouvertes sur $X$ et aux  changements de base admissibles sur $S$, on procède comme  dans \cite{K0} ou \cite{K1} en commençant par exhiber une telle
flèche relativement à une installation locale de $\pi$ puis
procéder à la globalisation grâce à l'exactitude à droite du foncteur ${\rm I}\!{\rm R}^{n}\pi_{!}$ et à la méthode du découpage.\vspace{2mm}

\indent
 {\bf(a) Construction locale.}\vspace{1mm}

\noindent
Soient  $s_{0}$ un point de $S$ et  $x$ un point de
la fibre $\pi^{-1}(s_{0})$.  Considérons une factorisation locale
de $\pi$ en $x$  donnée (abusivement) par
$\xymatrix{X\ar@/_/[rr]_{\pi}\ar[r]^{f}&Y\ar[r]^{q}&S}$ dans
laquelle $f$ est fini, ouvert  et surjectif sur $Y=S\times U$ avec
$U$ polydisque ouvert relativement compact de ${\Bbb C}^{n}$, $S$ un
ouvert de Stein; $q$ et $p$ étant les projections canoniques sur $S$ et $U$ respectivement, on a, pour tout entier $k$, les décompositions  
$$\overline{\Omega}^{k}_{Y}=\bigoplus_{i+j=k} q^{*}(\overline{\Omega}^{i}_{S})\otimes_{{\mathcal O}_{Y}}p^{*}(\Omega^{j}_{U})=\bigoplus_{i+j=k}\overline{\Omega}^{i}_{S}\widehat\otimes_{\Bbb C}\Omega^{j}_{U}$$ 
et, donc, une projection naturelle
$$\overline{\Omega}^{n+q}_{Y}\rightarrow q^{*}(\overline{\Omega}^{q}_{S})\otimes_{{\mathcal O}_{Y}}
\Omega^{n}_{Y/S}$$
Mais le morphisme canonique de ${\mathcal O}_{S}$-modules
$${\overline{\Omega}^{q}_{S}}\otimes_{{\mathcal O}_{S}}{\rm I}\!{\rm R}^{n}q_{!}\Omega^{n}_{Y/S}\rightarrow 
{\rm I}\!{\rm R}^{n}q_{!}(q^{*}(\overline{\Omega}^{q}_{S})\otimes_{{\mathcal O}_{Y}}
\Omega^{n}_{Y/S})$$
est un isomorphisme puisque, pour tout faisceau cohérent ${\mathcal G}$ sur $S$,  les foncteurs
$${\mathcal G}\rightarrow {\mathcal G}\otimes_{{\mathcal O}_{S}}{\rm I}\!{\rm R}^{n}q_{!}\Omega^{n}_{Y/S},\,\, {\mathcal G}\rightarrow{\rm I}\!{\rm R}^{n}q_{!}(q^{*}({\mathcal G})\otimes_{{\mathcal O}_{Y}}
\Omega^{n}_{Y/S})$$
coincident sur les faisceaux localement libres; le cas général s'en déduit en prenant une résolution localement libre à deux termes pour ${\mathcal G}$. Cela se vérifie aussi rapidement en utilisant une formule de type Künneth sur le produit complété.\vspace{1mm}

\noindent
Alors, en appliquant  le foncteur exact à droite ${\rm I}\!{\rm R}^{n}q_{!}$ à la flèche 
$$f_{*}\overline{\Omega}^{n+q}_{X}\rightarrow q^{*}(\overline{\Omega}^{q}_{S})\otimes_{{\mathcal O}_{Y}}\Omega^{n}_{Y/S}$$ 
on obtient le  morphisme
$${\rm I}\!{\rm R}^{n}q_{!}f_{*}\overline{\Omega}^{n+q}_{X}\rightarrow
\overline{\Omega}^{q}_{S}\otimes_{{\mathcal O}_{S}}{\rm I}\!{\rm R}^{n}q_{!}\Omega^{n}_{Y/S}$$
Mais, pour tout ouvert
$S'$ de Stein dans $S$, une formule de type K\"unneth donne
$${\rm I}\!{\rm R}^{n}q_{!}\Omega^{n}_{Y/S}(S')\simeq {\mathcal
O}_{S}(S') {{\widehat{\otimes}}{\atop_{{\Bbb C}}}}{\rm H}^{n}_{c}(U,
\Omega^{n}_{U})$$ qui, eu égard à l'isomorphisme ${\rm
H}^{n}_{c}(U, \Omega^{n}_{U})\simeq {\Bbb C}$ donné par la
dualité de Serre, produit, pour $S$ connexe,  l'isomorphisme bien connu pour un morphisme lisse, ${\rm I}\!{\rm
R}^{n}q_{!}\Omega^{n}_{Y/S}\simeq {\mathcal O}_{S}$. 
La finitude de  $f$ garantissant la dégénérescence de la suite spectrale de Leray, et, donc,  
l'isomorphisme de foncteurs ${\rm I}\!{\rm R}^{n}\pi_{!} \simeq {\rm
I}\!{\rm R}^{n}q_{!}f_{*}$,   permet de mettre en évidence la flèche désirée 
$${\rm I}\!{\rm R}^{n}\pi_{!}\overline{\Omega}^{n+q}_{X}\rightarrow
\overline{\Omega}^{q}_{S}$$
\indent {\bf  (b) Construction globale.}\vspace{1mm}

\noindent
Commençons par remarquer que la construction d'une telle
flèche, qui est toujours de nature locale sur $S$, est aussi  de
nature locale sur $X$ en vertu du  lemme de Reiffen (\cite{R}) impliquant, en particulier, 
l'exactitude à droite du foncteur ${\rm I}\!{\rm R}^{n}\pi_{!}$. Par ailleurs, $\overline{\Omega}^{n+q}_{X}$ étant de nature locale, on a,  pour tout ouvert  $U$ de $X$ muni de l'injection naturelle $j:U\rightarrow X$, 
$$\overline{\Omega}^{n+q}_{X}|_{U}=\overline{\Omega}^{n+q}_{U}$$
Cela étant dit, supposons $X$ paracompact et complètement
paracompact\footnote{Tout ouvert de $X$ est un espace
paracompact; en particulier, pour toute famille paracompactifiante
de supports, l'étendue $\bigcup_{F\in \Phi}F$ est un ouvert
paracompact.}  et considérons un ouvert de Stein (que l'on notera
encore $S$) de $S$ et un recouvrement ouvert localement fini
$(X_{\alpha})_{\alpha\in A}$ de $X$ qui soit $S$-adapté et muni
d'installations locales
 $$\xymatrix{X_{\alpha}\ar[rdd]_{\pi_{\alpha}}\ar[rr]^{\sigma_{\alpha}}\ar[rd]^{f_{\alpha}}&&Z_{\alpha}\ar[ld]_{p_{\alpha}}\ar[ldd]^{r_{\alpha}}\\
&Y_{\alpha}\ar[d]^{q_{\alpha}}&\\&S&}$$ avec $\sigma_{\alpha}$  un
plongement local, $\pi_{\alpha}$ la restriction de $\pi$ à
$X_{\alpha}$, $Y_{\alpha}$ et $Z_{\alpha}$ sont lisses sur $S$,
$p_{\alpha}$,$q_{\alpha}$ et $r_{\alpha}$ lisses.\vspace{1mm}

\noindent Comme
$\overline{\Omega}^{n+q}_{X}|_{X_{\alpha}} = \overline{\Omega}^{n+q}_{X_{\alpha}}$ on déduit,  du cas local précédent, une collection de morphismes
$${\rm I}\!{\rm R}^{n}{\pi_{\alpha}}_{!}\overline{\Omega}^{n+q}_{{X}_{\alpha}}\rightarrow \overline{\Omega}^{q}_{S}$$
Mais l'exactitude à droite du foncteur ${\rm I}\!{\rm
R}^{n}\pi_{!}$  qui assure la
surjectivité de la flèche naturelle
$$ \bigoplus_{\atop\alpha\in A}{\rm I}\!{\rm
R}^{n}{\pi_{\alpha}}_{!} \overline{\Omega}^{n+q}_{{X_{\alpha}}} \rightarrow
{\rm I}\!{\rm R}^{n}{\pi_{!}}{\overline{\Omega}^{n+q}_{X}}$$ permet de
produire, par recollement sur $X$, le morphisme désiré.\vspace{1mm}

\noindent {\bf(i)}  Cette construction qui a été entreprise pour le faisceau ${\mathcal L}^{\bullet}_{X}$ nous donnait (cf \cite{K0}, \cite{K4}) le morphisme
$\displaystyle{{\rm I}\!{\rm R}^{n}{\pi_{!}}{\mathcal L}^{n+q}_{X}\rightarrow {\mathcal L}^{q}_{S}}$
et s'adapte sans aucune modification au faisceau $\omega^{\bullet}_{X}$ pour avoir le morphisme
$\displaystyle{{\rm I}\!{\rm R}^{n}{\pi_{!}}{\omega}^{n+q}_{X}\rightarrow {\omega}^{q}_{S}}$
La commutativité des diagrammes
$$\xymatrix{{\rm I}\!{\rm R}^{n}{\pi_{!}}{\overline{\Omega}^{n+q}_{X}}\ar[d]_{{\mathcal T}^{j}_{\pi,\overline{\Omega} }}\ar[r]&{\rm I}\!{\rm R}^{n}{\pi_{!}}{\mathcal L}^{n+q}_{X}\ar[d]_{{\mathcal T}^{j}_{\pi,{\mathcal L}}} \ar@{=}[rr]&&{\rm I}\!{\rm R}^{n}{\pi_{!}}{\mathcal L}^{n+q}_{X}\ar[r]\ar[d]_{{\mathcal T}^{j}_{\pi,{\mathcal L}}}&{\rm I}\!{\rm R}^{n}{\pi_{!}}{\omega}^{n+q}_{X}\ar[d]_{{\mathcal T}^{j}_{\pi,\omega}}\\
{\Omega}^{q}_{S}\ar[r]&{\mathcal L}^{q}_{S}\ar@{=}[rr]&&{\mathcal L}^{q}_{S}\ar[r]&{\omega}^{q}_{S}}$$
est alors évidente. \vspace{1mm}

\noindent
{\bf(ii)} La compatibilité de la construction aux changements de base est facile à voir car à tout diagramme de changement de base
$$\xymatrix{\widetilde{X}\ar[d]_{\tilde\pi}\ar[r]^{\theta}&X\ar[d]^{\pi}\\
\widetilde{S}\ar[r]_{\nu}&S}$$
sont naturellement associés  des diagrammes  commutatifs (qui coincident tous sur la partie régulière de $\pi$)
$$\xymatrix{\nu^{*}{\rm I}\!{\rm R}^{n}{\pi_{!}}{\mathcal L}^{n+q}_{X}\ar[d]_{\nu^{*}({\mathcal T}^{j}_{\pi,{\mathcal L} })}\ar[r]&{\rm I}\!{\rm R}^{n}{\tilde\pi_{!}}\theta^{*}{\mathcal L}^{n+q}_{X}\ar[r]&{\rm I}\!{\rm R}^{n}{\tilde\pi_{!}}{\mathcal L}^{n+q}_{\widetilde X}\ar[d]^{{\mathcal T}^{j}_{\tilde\pi,{\mathcal L}}}\\
\nu^{*}{\mathcal L}^{q}_{S}\ar[rr]&&{\mathcal L}^{q}_{\widetilde S}}$$
pour tout morphisme $\nu$,
$$\xymatrix{\nu^{*}{\rm I}\!{\rm R}^{n}{\pi_{!}}{\overline{\Omega}^{n+q}_{X}}\ar[d]_{\nu^{*}({\mathcal T}^{j}_{\pi,\overline{\Omega} })}\ar[r]&{\rm I}\!{\rm R}^{n}{\tilde\pi_{!}}\theta^{*}{\overline{\Omega}^{n+q}_{X}}\ar[r]&{\rm I}\!{\rm R}^{n}{\tilde\pi_{!}}{\overline{\Omega}^{n+q}_{\widetilde X}}\ar[d]^{{\mathcal T}^{j}_{\tilde\pi,\overline{\Omega}}}\\
\nu^{*}{\overline{\Omega}}^{q}_{S}\ar[rr]&&{\overline{\Omega}}^{q}_{\widetilde S}}$$
si $\nu$ est un changement de base permis (i.e $\nu^{-1}(S)\subsetneq{\rm Sing}({\widetilde S})$ ),
$$\xymatrix{\nu^{*}{\rm I}\!{\rm R}^{n}{\pi_{!}}{{\omega}^{n+q}_{X}}\ar[d]_{\nu^{*}({\mathcal T}^{j}_{\pi,{\omega} })}\ar[r]&{\rm I}\!{\rm R}^{n}{\tilde\pi_{!}}\theta^{*}{{\omega}^{n+q}_{X}}\ar[r]&{\rm I}\!{\rm R}^{n}{\tilde\pi_{!}}{{\omega}^{n+q}_{\widetilde X}}\ar[d]^{{\mathcal T}^{j}_{\tilde\pi,{\omega}}}\\
\nu^{*}{\omega}^{q}_{S}\ar[rr]&&{\omega}^{q}_{\widetilde S}}$$
pour $\nu$ plat ou $S$ et $\widetilde S$ normaux.\vspace{1mm}

\noindent La compatibilité aux localisations sur $X$ et $S$ est
aisée à vérifier puisque c'est un cas particulier du changement de base plat.\vspace{1mm}

\noindent 
{\bf(iii)} Soient $Z$ et $S$  deux espaces  analytiques complexes de dimension pure avec 
$Z$ dénombrable 
à l'infini et $S$ réduit. Soit $(X_{s})_{s\in S}$ une famille analytique de $n$- cycles de $Z$ 
paramétrée par $S$ dont le graphe $X$ est supposé non entièrement inclus dans le lieu singulier de $S\times Z$ et muni du morphisme $\pi: X\rightarrow S$ (induit par la projection canonique sur $S$). Soient $\Phi$ et  $\Psi$ deux familles de supports avec $\Phi$ paracompactifiante et $\Phi \cap \Psi$ contenue dans la famille des compacts de $Z$, $\forall s\in S$, $X_{s} \in \Psi$ et  $X\cap (S\times \Phi)$ contenu dans la famille des ferm\'es $S$- 
propres.\vspace{1mm}

\noindent  
 L'existence du morphisme d'int\'egration sur les cycles 
$$\overline{\sigma}^{q,0}_{{\Phi},X}:{\rm H}^{n}_{\Phi}(Z, \overline{\Omega}^{n+q}_{Z})\longrightarrow 
{\rm H}^{0}(S, \overline{\Omega}^{q}_{S})$$
est un corollaire presqu'immédiat de ce qui précède. En effet, on met en évidence le diagramme
$$\xymatrix{{\rm H}^{n}_{\Phi}(Z, \overline{\Omega}^{n+q}_{Z})\ar[r]^{\alpha}\ar[d]_{\overline{\sigma}^{q,0}_{{\Phi},X}}&{\rm H}^{n}_{S\times \Phi}(S\times Z, \overline{\Omega}^{n+q}_{S\times Z})\ar[r]^{\beta}&{\rm H}^{n}_{S-propre}(X, \overline{\Omega}^{n+q}_{X})\ar[d]^{\gamma}\\
{\rm H}^{0}(S, \overline{\Omega}^{q}_{S})&&{\rm H}^{0}(S, {\rm I}\!{\rm R}^{n}{\pi_{!}}\overline{\Omega}^{n+q}_{X})\ar[ll]^{{\mathcal T}^{q}_{\pi,\overline{\Omega}}}}$$
où $\alpha$ et $\beta$ sont les morphismes induit par l'image réciproque et la restriction, $\gamma$ est un morphisme latéral de la suite spectrale
$${\rm E}^{i,j}_{2}={\rm H}^{i}(S, {\rm I}\!{\rm R}^{j}{\pi_{!}}\overline{\Omega}^{n+q}_{X})\Longrightarrow\,{\rm H}^{i+j}_{S-propre}(X, \overline{\Omega}^{n+q}_{X})$$
sachant que $${\rm I}\!{\rm R}^{j}{\pi_{!}}\overline{\Omega}^{n+q}_{X}=0,\,\forall\,j>n\,\blacksquare$$\vspace{1mm}

\noindent 
\subsection{\color{blue} {Preuve du \corollaryref{C0}}}: Rappelons les hypothèses:\vspace{1mm}

\noindent $\bullet$ $Z$ et $S$  deux espaces  analytiques complexes de dimension pure avec 
$Z$ dénombrable à l'infini et $S$ réduit, $(X_{s})_{s\in S}$ une famille analytique de $n$- cycles de $Z$ paramétrée par $S$ dont le graphe $X$ est supposé non entièrement inclus dans le lieu singulier de $S\times Z$ et muni du morphisme $\pi: X\rightarrow S$ (induit par la projection canonique sur $S$).\vspace{1mm}

\noindent $\bullet$  $\Psi$ une famille de supports de $Z$ contenant les supports des cycles $X_s$,  $\Theta$  une famille paracompactifiante de supports de $S$  telle que la famille $\widetilde{\Psi}:=(S\times (\Phi\cap \Psi))\cap X$ soit paracompactifiante et contenue dans la famille paracompactifiante des fermés $\Theta$-propre de sorte  que le couple $(\widetilde{\Psi}, \Theta)$ soit adapté à $\pi$.\vspace{1mm}

\noindent 
On montre, alors, l'existence d'un morphisme d'int\'egration sur les cycles 
$$\overline{\sigma}^{q,p}_{{\Phi},X}:{\rm H}^{n+p}_{\Phi}(Z, \overline{\Omega}^{n+q}_{Z})\longrightarrow 
{\rm H}^{p}_{\Theta}(S, \overline{\Omega}^{q}_{S})$$
vérifiant les propriétés fonctorielles induites par celles de ${\mathcal T}^{q}_{\pi, \overline{\Omega}}$.\vspace{2mm}

\indent Selon les notations de \cite{Br}, {\S  IV}, {\bf{p 210-223}}, (cf {\bf{prop.5.3} p.220}), la famille des fermés $\Theta$-propre  est notée et définie par $$\Theta(c):=\{F\in \pi^{-1}(\Theta)/ \,\pi\vert{F}: F\rightarrow S\,{\rm est\,propre}\} $$
D'après cette même proposition, 
$$\widetilde{\Psi}\subset \Theta(c)\Longleftrightarrow\, (\pi(\widetilde{\Psi})\subset \Theta\,{\rm et}\,
\pi\vert{F}: F\rightarrow S\,{\rm est\,propre}\,\forall\,F\in \widetilde{\Psi})$$
 Si $\displaystyle{c(S):=\{F\in X: \pi\vert{F}: F\rightarrow S\,{\rm est\,propre}\} }$ est la famille paracompactifiante des fermés $S$-propres, on a, d'après ({\bf{prop.5.4, 5)} }
 $$\Theta(c(S)):=c(S)(S)\cap \pi^{-1}(\Theta)= c(S)\cap \pi^{-1}(\Theta)=\Theta(c)$$
 De plus, $\Theta(c(S))$ est une famille paracompactifiante puisque $\Theta$ et $c(S)$ le sont ( cf {\bf{prop.5.5}}). Le {\bf{thm 6.1}} nous donne, alors, une suite spectrale de Leray
 $${\rm E}^{i,j}_{2}:={\rm H}^{i}_{\Theta}(S, {\mathcal H}^{j}_{c(S)}(\pi, \overline{\Omega}^{n+q}_{X}))\Longrightarrow\,{\rm H}^{i+j}_{\Theta(c)}(X, \overline{\Omega}^{n+q}_{X})$$
 qui s'écrit aussi
$${\rm E}^{i,j}_{2}:={\rm H}^{i}_{\Theta}(S, {\rm I}\!{\rm R}^{j}\pi_{!}\overline{\Omega}^{n+q}_{X}))\Longrightarrow\,{\rm H}^{i+j}_{\Theta(c)}(X, \overline{\Omega}^{n+q}_{X})$$
 de laquelle résulte le morphisme latéral
 $${\rm H}^{n+p}_{\Theta(c)}(X, \overline{\Omega}^{n+q}_{X})\rightarrow {\rm H}^{p}_{\Theta}(S, {\rm I}\!{\rm R}^{n}{\pi_{!}}\overline{\Omega}^{n+q}_{X})$$
  que l'on compose avec le morphisme induit par ${\mathcal T}^{j}_{\pi,\overline{\Omega}}$ pour en déduire le morphisme $\overline{\sigma}^{q,p}_{{\Phi},X}$ possédant les mêmes propriétés fonctorielles que $\overline{\sigma}^{q,0}_{{\Phi},X}\,\,\blacksquare$
  
\begin{rem}
Si l'on  suppose ${\mathcal T}^{0}_{\pi,\overline{\Omega} }$ construite, les  idées du \theoremref{T'1} permettent de définir 
${\mathcal T}^{j}_{\pi,\overline{\Omega} }$ pour tout $j>0$. En effet, le diagramme 
$$\xymatrix{{\rm I}\!{\rm R}\pi_{!}\overline{\Omega}^{n+q}_{X}[n]\ar[r]\ar[rrd]&{\rm I}\!{\rm R}\pi_{!}{\rm I}\!{\rm R}{\mathcal H}om(\overline{{\Omega}}^{r-q}_{X}, \overline{{\Omega}}^{n+r}_{X})[n]\ar[r]&{\rm I}\!{\rm R}{\mathcal H}om({\rm I}\!{\rm R}\pi_{*}\overline{\Omega}^{r-q}_{X}, {\rm I}\!{\rm R}\pi_{!}\overline{{\Omega}}^{n+r}_{X}[n] )\ar[d]\\
&&{\rm I}\!{\rm R}{\mathcal H}om(\overline{\Omega}^{r-q}_{S}, \overline{\Omega}^{r}_{S} )}$$ 
En prenant la cohomologie de degré $0$ de la flèche oblique et utilisant le \corollaryref{co1}, on obtient le morphisme ${\mathcal T}^{j}_{\pi,\overline{\Omega} }$.
\end{rem}

\noindent 
\subsection{Exemple d'intégration sur une famille analytique de cycles. }\rm
On va construire une famille analytique de cycles paramétrée par un espace analytique faiblement normal pour laquelle l'intégration sur les cycles de classes de cohomologie "holomorphes"  produit des formes méromorphes  avec condition de dépendance intégrale. \vspace{1mm}

\noindent
 Comme de coutume  $\rm{S}ym^2_{0}({\Bbb C}^{2})$ désignera la partie
homogène de degré $2$ de ${\rm Sym}^{2}({\Bbb C}^{2})$
canoniquement  identifiée à ${\Bbb C}^{2}/ \{id,-id\}$ qui est
isomorphe  au cône de ${\Bbb C}^3$. \vspace{1mm}

\noindent Soit $S:=\{(a,b,c)\in {\Bbb C}^{3}: a^{4} + a^{2}bc -c^{2} =0\}$ la surface faiblement normale dont le lieu singulier est la droite $\Sigma:=\{(a,b,c)\in {\Bbb C}^{3}: a=c=0\}$ et l'idéal de ${{\Bbb C}^{6}}$ donné par 
${\mathcal I}:=( u^2 - (a^{2} +bc)t^2, v^2 - a^2, uv -
ct, a^{4} + a^{2}bc -c^{2})$  définissant 
 un certain sous espace que l'on  notera abusivement $X=\{(u,v,a,b,c,t)\in {{\Bbb
C}^{2}}\times S\times D : u^2=(a^{2} +bc)t^2, v^2=a^2, uv=ct\}$, $D$ étant le disque unité. On
désignera par $\pi: X\rightarrow S$ (resp. $f:X\rightarrow S\times
D$) les morphismes induits par  les projections canoniques de
$S\times D\times{\Bbb C}^2\rightarrow S\times D\rightarrow S$ et par
 $F:S\times D\rightarrow \rm{S}ym^2_{0}({\Bbb C}^{2})$,
 l'application analytique
donnée par $(a,b,c,t)\rightarrow((a^{2} +bc)t^2, a^2, ct).$\vspace{2mm}

\noindent Il apparait clairement que  $X$, qui est  le produit
fibré de $S\times D$ par ${\Bbb C}^2$ au dessus de
$\rm{S}ym^2_{0}({\Bbb C}^{2})$, s'identifie, après quotient par
l'involution $(a,b,c,t,u,v)\rightarrow(a,b,c,t,-u,-v)$, au graphe de
l'application $F$. \vspace{1mm}

\noindent Pour $s$ générique, la fibre
$X_{s}:=\pi^{-1}(s)$  est constituée d'un couple de droites en
position générale dont les équations sont:
$$ \left\{
\begin{array}{lr}
u_1={c\over a}t,&u_2=-u_1\\
            \\
            v_1=a,&v_2=-v_1
\end{array}
\right. $$
  Remarquons au passage que $\pi$
n'est pas plat puisque de multiplicité $3$ en l'origine alors
qu'elle vaut  $2$ en les points  génériques. \vspace{2mm}

\noindent
  Comme  les fonctions de Newton fondamentales $u^{2}$, $v^{2}$ et $uv$
ont pour traces $S$-relatives respectives ${\mathcal
T}^{0}_{f}(u^{2})=2(a^{2} +bc)t^2$,  ${\mathcal T}^{0}_{f}(v^{2})=2a^{2}$ et
${\mathcal T}^{0}_{f}(uv)=2ct$, qui sont manifestement holomorphes sur
$S$, on a une famille analytique de revêtements ramifiés de
degré $2$ et de codimension $2$ de $D$ dans $D\times {\Bbb C}^{2}$,
``repéré'' par le morphisme $ F:S\times D\rightarrow
\rm{S}ym^2_{0}({\Bbb C}^{2})$, envoyant chaque couple $(s,t)$ sur
les fonctions symétriques des branches locales
$f_{1}=(u_{1},v_{1})$ et $f_{2}=-f_{1}$. Il s'en suit que $\pi$
définit une famille continue de $1$-cycles de ${\Bbb C}^{3}$ (on
dira plus brièvement que $\pi$ est continûment
géométriquement plat).
\vspace{2mm}

\noindent
L'analyticité de la famille se teste en regardant la trace
$S$-relative des formes de Newton. Dans notre cas,  le faisceau des
$1$- formes de Newton est engendré par les formes $udu$, $vdv$,
$udv$ et $vdu$ dont les traces $S$-relatives sont ${\mathcal
T}^{1}_{f}(udu)= (a^{2} +bc)tdt$, ${\mathcal T}^{1}_{f}(vdv)= 0$ et ${\mathcal
T}^{1}_{f}(udv-vdu)=- 2cdt$ qui sont clairement $S$- holomorphes.
Ainsi $\pi$ définit une famille analytique locale de 1-cycles de
${\Bbb C}^3 $ ( on dira que $\pi$ est analytiquement
géométriquement plat ).\vspace{2mm}

\noindent Remarquons que
l'identification $\rm{S}ym^2_{0}({\Bbb
C}^{2})\simeq\{(x,y,z)\in{\Bbb C}^{3}/ xy=z^{2}\}$, permet de voir
que la forme $udv-vdu$ correspond exactement (par image directe) à
la forme méromorphe (non holomorphe) ${z({dx\over{x}} -
{dy\over{y}}})$ sur le c\^one; elle définit naturellement une
section du faisceau ${\mathcal L}^{1}_{\rm{S}ym^2_{0}({\Bbb
C}^{2})}$.\vspace{2mm}

\noindent
Considérons la forme $udv$ dont la trace absolue donne la forme méromorphe
 $${\mathcal
T}^{1}_{f}(udv)=2{c\over{a}}da$$
qui définit une section du faisceau $\omega^{1}_{S}$. Montrons qu'elle n'est pas section du faisceau $\overline{\Omega}^{1}_{S}$.\vspace{1mm}

\noindent En utilisant les relations différentielles données par les équations de définition, on obtient la relation minimale 
$$\frac{c}{a}da  + \frac{c^{2}}{2bc + 4a^{2}}db =\frac{1}{2}dc$$
avec $\alpha:=\frac{c}{a}da$ et $\beta:=\frac{c^{2}}{2bc + 4a^{2}}db$ toutes deux sections de $\omega^{1}_{S}$. Comme $\xi_{1}:=(\frac{c}{a}da)$ vérifie $\xi_{1}^{\otimes 2}=(a^2 +bc)da^{\otimes 2}$ qui est une équation de dépendance intégrale, la forme $\xi_{2}:=\frac{c^{2}}{2bc + 4a^{2}}db)$ vérifie aussi une équation de dépendance intégrale et qui est
$$\xi^{\otimes 2}_{2} +\xi_{1}\otimes dc + (\xi^{\otimes 2}_{1} - \frac{1}{4}dc^{\otimes 2})=0$$
Comme on le constate la relation différentielle  minimale met en évidence des formes méromorphes vérifiant des équations de dépendance intégrale différentes. \vspace{1mm}

\noindent
D'ailleurs le coefficient de  la forme produit tensoriel $\xi_{1}\otimes \xi_{2}$ n'est pas holomorphe puisqu'il s'écrit, modulo les fonctions holomorphes, à la fonction méromorphe
$$g(a,b,c)=\frac{ac^{3}}{c^{2} +a^{4}}$$
dont on teste l'analyticité en faisant, par exemple,  $a=b$ qui nous donne la courbe 
$${\mathcal C}:=\{(a,c): a^{4} +ca^{3}-c^{2}=0\}$$
qui, pour l'éclatement $(a,c)\rightarrow(uc,c)$, donne la préimage stricte lisse
$$\widetilde{{\mathcal C}}:=\{(u,c): c^{2}u^{3} +c^{2}u^{4}-1=0\}$$
sur laquelle l'image réciproque stricte de $g$ donne la fonction 
$$\widetilde{g}(c,u)=\frac{uc^{2}}{1+u^{4}c^{2}}$$
qui présente manifestement au moins un pôle en $u=\frac{-1}{2}$, $c=\pm 4i$.\vspace{1mm}

\noindent 
En étudiant l'exemple suivant
$$X:=\{(x,y,z,t, \alpha, \beta, \gamma, \delta, t',u,v)\in {\Bbb C}^{8}:xy=zt, xt=z^{2}, yz=t^{2}; \alpha\beta^{3}=\gamma^{3}\delta^{2}, \alpha\delta^{2}=\gamma^{6}, \beta^{3}\gamma^{3}=\delta^{2};  x=u^{3} =\alpha t'^{3},$$
$$\,\,\,\,\,\,\,\,\,\,\,\,y=v^{3} =\beta t'^{3}, z=vu^{2} =\gamma^{3} t'^{3}, t=uv^{2} =\delta^{2} t'^{3}\}$$
avec $$S:=\{(\alpha, \beta, \gamma, \delta)\in {\Bbb C}^{4}: \alpha\beta^{3}=\gamma^{3}\delta^{2}, \alpha\delta^{2}=\gamma^{6}, \beta^{3}\gamma^{3}=\delta^{2}\},$$
on montre que la propriété de  dépendance intégrale est préservée par intégration sur cette famille analytique de cycles.\vspace{2mm}

\phantomsection\addcontentsline{toc}{part}{Références.}

%
\svnid{$Id: S01-intro.tex 272 2020-01-21 14:29:32Z kebekus $}

\end{document}